\documentclass[twoside, 10pt]{article}

\usepackage[hmargin= 0.8in, height=9.3in]{geometry}
\usepackage{amssymb,amsthm, amsfonts,amsmath,mathrsfs, amsxtra, url, hyperref}

\usepackage{fancyhdr}
\fancypagestyle{plain}{
 \fancyhf{}

}

\renewcommand{\theequation}{\thesection.\arabic{equation}}
 \numberwithin{equation}{section}

\newtheorem {thm}{Theorem}[section]
\newtheorem {prop}{Proposition}[section]
\newtheorem {lemm}{Lemma}[section]
\newtheorem {deff}{Definition}[section]
\newtheorem {cor}{Corollary}[section]
\newtheorem {rem}{Remark}[section]

\newtheorem {eg}{Example}[section]

\def\ba{\begin{array}}
\def\ea{\end{array}}
\def\bea{\begin{eqnarray}}
\def\eea{\end{eqnarray}}
\def\beas{\begin{eqnarray*}}
\def\eeas{\end{eqnarray*}}
\def\bi{\begin{itemize}}
\def\ei{\end{itemize}}
\def\bc{\begin{cases}}
\def\ec{\end{cases}}

\def\a{\alpha}
\def\g{\gamma}
\def\d{\delta}
\def\e{\varepsilon}
\def\z{\zeta}
\def\k{\kappa}
\def\l{\lambda}

\def\si{\sigma}

\def\t{\tau}
\def\th{\theta}
\def\o{\omega}
\def\f{\phi}
\def\vf{\varphi}

\def\D{\Delta}
\def\G{\Gamma}
\def\L{\Lambda}
\def\O{\Omega}
\def\F{\Phi}

\def\Th{\Theta}
\def\U{\Upsilon}

\def\bF{{\bf F}}

\def\cA{{\cal A}}
\def\cB{{\cal B}}

\def\cD{{\cal D}}
\def\cE{{\cal E}}
\def\cF{{\cal F}}
\def\cG{{\cal G}}
\def\cH{{\cal H}}
\def\cI{{\cal I}}

\def\cN{{\cal N}}

\def\cS{{\cal S}}
\def\cT{{\cal T}}
\def\cU{{\cal U}}
\def\cV{{\cal V}}

\def\cX{{\cal X}}
\def\cY{{\cal Y}}
\def\cZ{{\cal Z}}

\def\hC{\mathbb{C}}

\def\hE{\mathbb{E}}

\def\hG{\mathbb{G}}
\def\hH{\mathbb{H}}

\def\hL{\mathbb{L}}
\def\hM{\mathbb{M}}
\def\hN{\mathbb{N}}

\def\hR{\mathbb{R}}
\def\hS{\mathbb{S}}

\def\hU{\mathbb{U}}
\def\hV{\mathbb{V}}

\def\hZ{\mathbb{Z}}

\def\sB{\mathscr{B}}

\def\sO{\mathscr{O}}
\def\sP{\mathscr{P}}

\def\sV{\mathscr{V}}

\def\fA{\mathfrak{A}}
\def\fB{\mathfrak{B}}

\def\fD{\mathfrak{D}}

\def\fO{\mathfrak{O}}
\def\fP{\mathfrak{P}}

\def\ff{\mathfrak{f}}

\def\fx{\mathfrak{x}}

\def\fg{\mathfrak{g}}
\def\fw{\mathfrak{w}}

\def\({\textnormal{(}}
\def\){\textnormal{)}}
\def\[{[\neg[}
\def\]{]\neg]}
\def\lan{\langle}
\def\ran{\rangle}

\def\no{\noindent}

\def\ss{\smallskip}
\def\ms{\medskip}

\def\q{\quad}
\def\qq{\qquad}

\def\neg{\negthinspace}
\def\dneg{\neg \neg}
\def\tneg{\neg \neg \neg}

\def\ol{\overline}
\def\ul{\underline}
\def\ua{\mathop{\uparrow}}
\def\da{\mathop{\downarrow}}

\def\wt{\widetilde}
\def\wh{\widehat}

\def\pas{{\hbox{$P-$a.s.}}}

\def\hb{\hbox}
\def\dis{\displaystyle}
\def\cd{\cdot}
\def\cds{\cdots}

\def\fa{\,\forall \,}
\def\pa{\partial}
\def\es{\emptyset}

\def\dfnn{\stackrel{\triangle}{=}}
\def\b1{{\bf 1}}
\def\qed{\hfill $\Box$ \medskip}

\def\essinf{\mathop{\rm essinf}}
\def\esssup{\mathop{\rm esssup}}
\def\liminf{\mathop{\ul{\rm lim}}}
\def\limsup{\mathop{\ol{\rm lim}}}
\newcommand{\esup}[1]{ \underset{#1}{\esssup}\,}
\newcommand{\einf}[1]{ \underset{#1}{\essinf}\,}
\newcommand{\lsup}[1]{ \underset{#1}{\limsup}}
\newcommand{\linf}[1]{ \underset{#1}{\liminf}}
\newcommand{\lmt}[1]{ \underset{#1}{\lim}}
\newcommand{\lmtu}[1]{ \underset{#1}{\lim} \neg \ua \,}
\newcommand{\lmtd}[1]{ \underset{#1}{\lim} \neg \da \,}

\begin{document}

\title{\bf A Weak Dynamic Programming Principle for Zero-Sum Stochastic Differential Games with Unbounded Controls
}

\author{
  Erhan Bayraktar\thanks{ \noindent Department of
  Mathematics, University of Michigan, Ann Arbor, MI 48109; email:
{\tt erhan@umich.edu}.}  \thanks{E. Bayraktar is supported in part by the National Science Foundation under  a Career grant DMS-0955463 and an applied mathematics research grant DMS-1118673.} $\,\,$,
$~~$Song Yao\thanks{
\noindent Department of
  Mathematics, University of Pittsburgh, Pittsburgh, PA 15260; email: {\tt songyao@pitt.edu}. } }
\date{ }

\maketitle

  \begin{abstract}

We analyze a  zero-sum stochastic differential game between two competing players who can  choose unbounded controls.
The payoffs of the game are defined through    backward  stochastic differential equations.
We prove that  each player's  priority  value   satisfies a weak dynamic programming principle
  and thus solves    the associated fully non-linear partial differential equation in the viscosity sense.

\ss  \noindent   {\bf Keywords:}\: Zero-sum stochastic differential games, Elliott-Kalton strategies, weak dynamic programming principle,     backward stochastic differential equations, viscosity solutions,   fully non-linear PDEs.

\end{abstract}

    \tableofcontents

\section{Introduction}

    \setcounter{equation}{0}

\label{sec:intro}

In this paper  we extend the study of Buckdahn and Li  \cite{Buckdahn_Li_1} on a zero-sum  stochastic differential game (SDG),
 whose payoffs  are generated by  backward stochastic differential equations (BSDEs),
 to the case of super-square-integrable  controls (see Remark \ref{rem_super_square_control}).

    Since the seminal paper by Fleming and Souganidis  \cite{Fleming_1989},
    the SDG theory   has    grown rapidly in many aspects
   (see e.g. the references in    \cite{Buckdahn_Li_1}, \cite{Buckdahn_Li_4}).
     Among these developments,  Hamad\`ene et al. \cite{Hamadene_Lepeltier_1995b,  Hamadene_Lepeltier_Peng_1997, EH-03}
   introduced  a  (decoupled) SDE-BSDE system,   with controls only in the drift coefficients,
      to generate the   payoffs    in their studies of saddle point problems of SDGs.
 (For the  evolution and applications of the BSDE theory, see Pardoux and Peng \cite{PP-90},
   El Karoui et al.\,\cite{EPQ-97} and the   references therein.) 
 Later on, \cite{Buckdahn_Li_1} as well as its sequels \cite{Buckdahn_Li_2,Buckdahn_Li_3, Buckdahn_Li_4}
     generalized the SDE-BSDE framework so that     the  two  competing controllers
   can also influence the diffusion coefficient of the state dynamics.
   Unlike   \cite{Fleming_1989},   \cite{Buckdahn_Li_1}     used   a uniform  canonical space
 $\O = \big\{\o  \in  \hC([0,T]; \hR^d ) \neg : \o(0)  =  0 \big\} $
 so that   admissible control processes
  can also depend on the information occurring before the start of the game.
  Such a setting
 allows the authors of \cite{Buckdahn_Li_1} get around a relatively complicated 
  approximation  argument of \cite{Fleming_1989} which was due to a measurability issue (see Remark \ref{rem_measurability_issue}),
  and allows them  to adopt   the notion of  stochastic backward semigroups  and a  BSDE method,
  developed in \cite{Peng_1997, Peng_1992},
 to obtain results similar to \cite{Fleming_1989}:
 the lower and  upper values  of the SDG 
 satisfy a dynamic programming principle and  solve   the associated Hamilton-Jacobi-Bellman-Isaacs equations
 in the viscosity sense.
 However,  \cite{Buckdahn_Li_1}, \cite{Fleming_1989}  as well as some latest advances to the SDG theory
 (e.g. Bouchard et al. \cite{Bouchard_Moreau_Nutz_2012} on stochastic target games,
 Peng and Xu \cite{Peng_XuXm_2010} on SDGs in form of a generalized BSDE with random default time)
 still assume the compactness of  control spaces while Pham and Zhang \cite{Pham_Zhang_2012} on weak formulation of SDGs assumes the boundedness of coefficients in control variables. We are going to address these particular issues.

 In the present paper, since two   players  take super-square-integrable admissible controls over two separable metric spaces  $\hU$ and $\hV$  not necessarily compact, those approximation methods of \cite{Fleming_1989} and \cite{Buckdahn_Li_1} in proving
    the dynamic programming principle are no longer effective. Instead,  we derive a weak
     form of dynamic programming principle in  spirit of Bouchard and Touzi \cite{Bouchard_Touzi_2011} and use it to show that
    each player's priority value  solves the corresponding   fully non-linear PDE in the viscosity sense.
   Vitoria  \cite{PMASCVitoria_2010} has tried to extend the SDG for unbounded controls by proving a weak dynamic programming principle.  However, it still assumed that the control space of the player with priority is compact, see Theorem 75 therein.

 Square-integrable controls were initially considered
 by Krylov \cite[Chapter 6]{Krylov_CDP}, however, for cooperative games
  (i.e. the so called $\sup$ $\sup$ case).
   Browne \cite{Browne_2000} studied a
 specific zero-sum investment game between two small investors who  control the game
  via their square-integrable portfolios. Since the PDEs in this   case have smooth solutions,
    the problem can be solved by   a verification theorem instead of the dynamic programming principle.
     Inspired by the ``tug-of-war" (a discrete-time random turn game, see e.g.
 \cite{PSSW_2009} and \cite{LLPU_1994}),  Atar and Budhiraja  \cite{Atar_Budhiraja_2010} studied
 a zero-sum  stochastic differential game with $\hU=\hV=\{x \in \hR^n: |x|=1 \} \times  [0,\infty) $
  played until the state process exits a given domain.
     As in Chapter 6 of \cite{Krylov_CDP}, the authors approximated such a game with unbounded controls
by a sequence of games with bounded controls which satisfy a dynamic programming principle.
  They showed the equicontinuity of the approximating sequence and thus proved that
  the value function of the game is a unique viscosity solution to the inhomogenous infinity Laplace equation.
 We do not rely on this approximation scheme but directly prove a weak dynamic programming principle
for the game with super-square-integrable controls.

  Following  the probabilistic setting of \cite{Buckdahn_Li_1} (see Remark \ref{rem_measurability_issue}),
  our paper takes the   canonical space
 $\O = \big\{\o  \in  \hC([0,T]; \hR^d ) \neg : \o(0)  =  0 \big\} $, whose
coordinator process $B$ is a Brownian motion under the Wiener measure $P$.
  When  the game starts from time $t \in [0,T]$, under the super-square-integrable controls $\mu \in \cU_t$
  and $\nu \in \cV_t$ selected by player I and II respectively,
  the state  process    $X^{t,\xi,\mu,\nu}$ starting from a random initial state $\xi$ 
  will then evolve according to a  stochastic differential equation (SDE):
    \bea
  X_s= \xi + \int_t^s b \big(r, X_r,\mu_r,\nu_r \big) \, dr
  + \int_t^s \si \big(r, X_r,\mu_r,\nu_r \big) \, dB_r,  \q s \in [t, T] ,  \label{FSDE}
    \eea
 where the drift $b  $ and the diffusion $\si  $ are Lipschitz continuous in $x$
 and have linear growth in $(u,v)$.
 The payoff player I will receive from  player II is determined by the first component  of the unique solution
  $ \big(Y^{t,\xi,\mu,\nu},Z^{t,\xi,\mu,\nu}  \big) $ to the following BSDE:
   \bea   \label{BSDE00}
   Y_s= g \big( X^{t,\xi,\mu,\nu}_T \big)  \neg + \neg  \int_s^T  \neg
    f \big( r, X^{t,\xi,\mu,\nu}_r, Y_r,   Z_r, \mu_r, \nu_r \big)  \, dr
   \neg-\neg \int_s^T \neg Z_r d B_r  , \q    s \in [t,T]    .
          \eea
          Here the generator $f  $  is  Lipschitz continuous in $(y,z)$ and also has linear growth in $(u,v)$.
     When $   g $ and $  f $ are  $2/p-$H\"older continuous
  in $x $ for some $p \in (1,2]$,    $ Y^{t,\xi,\mu,\nu}$ is $p-$integrable.
  As we see from \eqref{FSDE} and \eqref{BSDE00} that the controls $\mu$, $\nu$
  influence the game in two aspects: both affect  \eqref{BSDE00} via the state process $X^{t,\xi,\mu,\nu}$
  and appear directly in the generator $f$ of   \eqref{BSDE00} as parameters.
  In particular, if $f$   is independent of $(y,z)$, $Y$ is in form of the  conditional linear expectation
    of the terminal reward $ g \big( X^{t,\xi,\mu,\nu}_T \big) $ plus the cumulative reward
    $\int_s^T  \neg  f (r, X^{t,\xi,\mu,\nu}_r,   \mu_r, \nu_r)  \, dr $ (cf. \cite{Fleming_1989}).

   When the player (e.g. Player I) with the priority
   chooses firstly  a super-square-integrable control (e.g. $\mu \in \cU_t$),
    its opponent (e.g. Player II) will   select its reacting control
   via a non-anticipative mapping $\beta: \cU_t \to \cV_t$, called \emph{Elliott-Kalton strategy}.
    In particular, using Elliott-Kalton strategies is essential in proving the dynamic programming principle. This phenomenon already appears in the controller-stopper games, i.e. when one of the players is endowed with the right of stopping the game instead of using a control; see \cite{BH11}, which shows that if the stopper acts second it is necessary that the stopper uses non-anticipative strategies in order to prove a dynamic programming principle. This type of phenomenon does not appear (or it is implicitly satisfied) if the controllers only control the drift, see e.g. \cite{OS_CRM} and the references therein, or when there are two stoppers (the so-called Dynkin games), see e.g. \cite{BS12} and the references therein. 

    By  $ w_1 (t,x)  \dfnn \underset{\beta \in \fB_t  }{\essinf} \; \underset{\mu \in \cU_t }{\esssup}
  \; Y^{t,x,\mu, \beta (  \mu )  }_t $ we  denote Player I's  priority  value of the  game starting
   from time $t  $ and state $x$,
  where $\fB_t$ collects  all admissible    strategies for Player II.
   Switching the priority defines Player II's priority value $w_2 (t,x)$.

    Although our setting
     makes the  payoffs $ Y^{t,\xi,\mu, \nu  }_t $        random variables,
   we can show like \cite{Buckdahn_Li_1} that $w_1 (t,x)  $ and $w_2 (t,x)  $ are invariant   under Girsanov
transformation via functions of  the Cameron-Martin space and are thus deterministic, see Lemma \ref{prop_value_constant}.
  To assure  values  $w_1(t,x)$ and $w_2 (t,x)  $ are finite,
      we assume that each player has  some  control neutralizer for coefficients $(b, \si, f)$
      (such an assumption holds for additive controls, see Example \ref{eg_control_neutralizer}),
       and we impose a growth condition on   strategies.  These two  requirements  are also crucial in proving our weak dynamic programming principle.
  When $\hU$ and $\hV$ are compact,  the control neutralizers become futile
and the growth condition  holds automatically for strategies. Thus
 our problem degenerates to   \cite{Buckdahn_Li_1}'s case, see Remark  \ref{rem_BL_case}.

  \ss  Although value functions $w_1 (t,x)$, $w_2(t,x)$ are still $2/p-$H\"older continuous in $x$
  (see Proposition \ref{prop_w_conti}),   they may not   be continuous in $t$.
   Hence we can not follow   \cite{Buckdahn_Li_1}'s approach to
 get  a strong form of  dynamic programming principle for $w_1$ and $w_2$.
 Instead, we prove a weak dynamic  programming principle, say for $w_1$:
           $$
       \underset{\beta \in \fB_t }{\essinf} \; \underset{\mu \in \cU_t }{\esssup}    \;
    Y^{t,x,\mu, \beta (  \mu ) }_t \Big(\t_{\beta,\mu},
    \f \big(\t_{\beta,\mu},  X^{t,x,\mu, \beta (  \mu ) }_{\t_{\beta,\mu}} \big) \Big)
   \neg  \le  \neg  w_1(t,x) \le    \underset{\beta \in \fB_t }{\essinf} \; \underset{\mu \in \cU_t }{\esssup}    \;
    Y^{t,x,\mu, \beta (  \mu ) }_t \Big(\t_{\beta,\mu},
    \wt{\f} \big(\t_{\beta,\mu},  X^{t,x,\mu, \beta (  \mu ) }_{\t_{\beta,\mu}} \big) \Big) ,
    $$
   for any two continuous functions        $    \f      \le     w_1     \le    \wt{\f}   $. Here
        $\t_{\beta,\mu}$ denotes the first existing time of state process $X^{t,x,\mu, \beta (  \mu )}$ from
        the given open ball  $O_{  \d} (t, x)$.

 To prove  the weak dynamic programming principle, we first
 approximate $w_1(t,x) = \underset{\beta \in \fB_t  }{\essinf} \; I(t,x,\beta)  $ from above and
 $ I(t,x,\beta) \dfnn \underset{\mu \in \cU_t }{\esssup}
  \; Y^{t,x,\mu, \beta (  \mu )  }_t $  from below in a probabilistic sense (see Lemma \ref{lem_mu_lmtu}) so that
    we  can  construct  approximately optimal controls/strategies by a pasting technique  similar to
 the one used in \cite{Bouchard_Touzi_2011} and   \cite{PMASCVitoria_2010}.
 Then we make   a series of estimates and eventually obtain  the weak dynamic programming principle
  by using a stochastic backward semigroup  property
  \eqref{eq:p677}, the continuous dependence of  payoff process  on the initial state
    (see Lemma \ref{lem_estimate_Y}) as well as
  the control-neutralizer assumption and the growth condition on strategies.

     Next,  one can deduce from the weak dynamic programming principle and
       the separability of control space $\hU$, $\hV$ that the value functions $w_1$ and $w_2$ are  (discontinuous)
  viscosity solutions  of  the corresponding   fully non-linear PDEs, see Theorem \ref{thm_viscosity}.
  Recently, Krylov \cite{Krylov_2012a} and \cite{Krylov_2012b}  studied the regularity of solutions to
    related  fully non-linear   PDEs: The former obtained
   $C^{1,1}\cap W^{1,2}_{\infty,loc}-$solutions for the case of bounded measurable coefficients;
     while the latter showed the existence of
   $L^p-$viscosity solutions in $C^{1+\a}$  if the fully non-linear Hamiltonian function
    is continuous in gradient variable and
    Lipschitz continuous in Hessian variable.

 The rest of the paper is organized as follows:
 After listing the notations to use, we recall some basic properties of BSDEs in Section \ref{sec:intro}.
 In Section \ref{sec:zs_drgame}, we set up the zero-sum stochastic differential games based on BSDEs and
 present a weak dynamic programming principle
 for priority values of both players defined via Elliott-Kalton strategies.
 With help of the weak dynamic programming principle, we show in Section \ref{sec:PDE} that  the   priority values
 are (discontinuous) viscosity solutions of the corresponding    fully non-linear   PDEs.
 The proofs of our results are deferred to Section \ref{sec:Proofs}.

\subsection{Notation and Preliminaries} \label{subsec:preliminary}


        Let   $ (\hM, \rho_{\overset{}{\hM}}) $ be a generic metric space
   and let  $\sB(\hM)$ be  the Borel $\si-$field of $\hM$. For any $x \in \hM $     and $\d >0$,
            $O_\d(x) \neg \dfnn \neg  \{x' \in \hM \neg  : \rho_{\overset{}{\hM}}(x,x')  < \d \}$
          and $\ol{O}_\d(x) \neg  \dfnn \neg  \{x' \in \hM  \neg : \rho_{\overset{}{\hM}}(x,x')  \le  \d \}$ respectively denote
           the open and closed  ball     centered at $x   $     with radius $\d  $.
                For any function $\phi: \hM \to \hR$,  we define the lower/upper semi-continuous envelopes by
    \beas
      \linf{x' \to x} \phi(x') \dfnn \lmtu{n \to \infty} \underset{x' \in O_{\neg \frac{1}{n}}(x)}{\inf} \phi(x')
     \q \hb{ and }  \q  \lsup {x' \to x} \phi(x') \dfnn \lmtd{n \to \infty} \underset{x' \in O_{\neg \frac{1}{n}}(x)}{\sup} \phi(x')    ,
    \eeas
    where $  \lmtd{n \to \infty} $ (resp. $ \lmtu{n \to \infty} $) denotes the limit of a decreasing (resp. increasing) sequence.

    Fix $d    \in \neg  \hN$ and   a  time horizon $T  \in (0,\infty)  $. We  consider
  the   canonical space $\O \neg \dfnn \neg \big\{\o  \neg  \in  \neg  \hC([0,T]; \hR^d) \neg : \o(0)  \neg = \neg  0 \big\}$
    equipped with Wiener measure $P$, under which    the   canonical process
      $B $        is a  $d-$dimensional Brownian motion.
     Let      $\bF =   \{ \cF_t    \}_{t \in [0,T]}$
   be the   filtration generated by $ B $ and augmented by all $P-$null sets.
   We denote by $\sP$      the   $\bF-$progressively  measurable $\si-$field of $ [0,T] \times \O$.

  Given $t \in [0,T]$, Let $\cS_{t,T}$ collect all $\bF-$stopping times $\t$ with $t \le \t \le T $, \pas ~
   For any   $\t \in \cS_{t,T}$ and $A \in \cF_\t$, we define
   $\[t,\t\[_A  \; \dfnn   \big\{ (r,\o) \in [t,T] \times A \neg : r  <  \t(\o) \big\} $
    and
   $ \[\t,T\]_A \dfnn \big\{ (r,\o) \in [t,T] \times A \neg : r  \ge  \t(\o) \big\}$  for any $A \in \cF_\t$.
 In particular, $ \[t,\t\[ \, \dfnn \[t,\t\[_\O    $ and $ \[\t,T\] \dfnn \[\t,T\]_\O$ are the
 stochastic intervals.

   Let      $\hE$ be a generic Euclidian space. For any $p \in [1,\infty)$ and   $t \in [0,T]$,
      we introduce some   spaces of functions:

  \ss \no 1)  For any sub$-\si-$field  $\cG$  of $\cF_T$,  let
     $\hL^p(\cG,\hE) $ be  the space of all  $\hE-$valued,
$  \cG-$measurable random variables $\xi$ such that $\| \xi \|_{\hL^p(\cG, \hE)}
   \dfnn  \big\{E\big[|\xi|^p\big]\big\}^{1/p  } < \infty$, and let
     $\hL^\infty (\cG,\hE) $ be  the space of all  $\hE-$valued,
$  \cG-$measurable bounded random variables.

\ss \no 2)         $\hC^p_\bF   ([t,T], \hE)$ denotes the space of all $\hE-$valued, $\bF-$adapted processes
  $\{X_s\}_{s \in [t,T]}$ with \pas ~ continuous paths
  such that $\|X\|_{\hC^p_\bF( [t,T],\hE)}
  \dfnn \bigg\{ E \Big[\,\underset{s \in [t,T]}{\sup}|X_s|^p \Big]  \bigg\}^{ 1 / p }<\infty$.

 \ss \no 3)        $ \hH^{p,loc}_\bF   ([t,T], \hE) $  denotes the space of all $\hE-$valued,
 $\bF-$progressively   measurable processes $\{X_s\}_{s \in [t,T]}$ such that
   $ \int_t^T \neg |X_s|^p   \,  ds <\infty $, $P-$a.s.
    For any $\wh{p} \in [1,\infty)$,      $ \hH^{p,\wh{p}}_\bF   ([t,T], \hE) $ denotes  the space of all $\hE-$valued,
 $\bF-$progressively   measurable processes $\{X_s\}_{s \in [t,T]}$ with
   $
   \|X\|_{\hH^{p,\wh{p}}_\bF([t,T], \hE)} \dfnn \left\{  E \left[  \big( \hb{$\int_t^T \neg |X_s|^p   \,  ds$}
  \big)^{ \wh{p} / p } \right] \right\}^{ 1 / \wh{p} }<\infty $.

  \ss \no 4)   We also set  $ \hG^p_\bF([t,T]) \dfnn \hC^p_\bF( [t,T], \hR) \times  \hH^{2,p} _\bF([t,T], \hR^d) $.

 \ss  If $\hE=\hR$, we will drop it from the above notations. Moreover, we will use  the convention $ \inf \es = \infty$.

   \subsection{  Backward Stochastic Differential Equations}

    Given $t \in [0,T]$, a  $t-$parameter set $\big(\eta,f  \big)$ consists of
  a random variable   $\eta \in \hL^0 \big( \cF_T \big)$ and
  a function $f: [t,T] \times \O \times \hR \times \hR^d \to \hR $   that
   is $\sP   \otimes \sB(\hR)  \otimes \sB(\hR^d)/\sB(\hR)-$measurable.
   In particular, $\big(\eta,f \big)$
 is called a $(t,p)-$parameter set for some $p \in [1,\infty)$ if $\eta \in \hL^p \big( \cF_{ T} \big)$.

 \begin{deff}
 Given a $t-$parameter set $\big(\eta,f  \big)$ for some $t \in [0,T]$, a pair $(Y,Z )
 \in \hC^0_{\bF}([t,T]) \times \hH^{2,loc}_{\bF}([t,T],\hR^d)  $ is called a solution of the
   backward stochastic differential equation  on the probability space $(\O,  \cF_T,  P )$ over period $[t,T]$
    with terminal condition $\eta$ and generator $f$
    \big(BSDE$ \big(t, \eta,f   \big)$ for short\big) if it holds $P-$a.s. that
     \bea   \label{BSDE01}
    Y_s = \eta \neg + \neg  \int_s^T  \neg  f  (r,   Y_r, Z_r)  \, dr
   \neg-\neg \int_s^T \neg Z_r d B_r  , \q    s \in [t,T] .
    \eea
\end{deff}

   Analogous to Theorem 4.2 of \cite{BH_Lp_2003}, we have the following well-posedness result of  BSDE \eqref{BSDE01}.

\begin{prop} \label{BSDE_well_posed}
 Given $t \in [0,T]$ and $p \in [1,\infty)$, let  $ \big(\eta,f \big)$ be a $(t,p)-$parameter set
 such that $f$     is  Lipschitz continuous in $(y,z)$: i.e.
   for some $\g > 0 $, it holds for $ds \times dP-$a.s.    $(s, \o) \in [t,T] \times \O $ that
   \beas  
 \big| f (s, \o,y ,z ) - f (s,\o,y', z') \big| \le \g \big( |y-y'| + |z-z'| \big), \q \fa y, y' \in \hR,~ \fa z,z' \in \hR^d .
   \eeas
   If    $    E   \Big[  \big(\int_t^T  \neg \big| f    (s, 0,0 )   \big|  ds \big)^p \Big] \neg < \neg \infty $,
     BSDE \eqref{BSDE01}
        admits a unique solution $ \big(Y ,Z   \big)  \neg \in \neg  \hG^p_\bF\big([t,T]\big) $  that
         satisfies
         \bea   \label{eq:xvx131}
        E \bigg[ \, \underset{s \in [t,T]}{\sup}  |  Y_s     |^p  \Big| \cF_t  \bigg]
    \le  C(T,p,\g) E \bigg[  \,   |    \eta   |^p
   + \Big( \int_t^T  \big| f  (s,0,0 ) \big|  ds \Big)^p \bigg| \cF_t  \bigg], \q   \pas
        \eea

 \end{prop}

\ss  Also,    
 we  have the following  a priori estimate and comparison results  for
  BSDE \eqref{BSDE01}.

 \begin{prop} \label{prop_BSDE_estimate_comparison}
Given  $t \in [0,T]$ and $p \in [1,\infty)$, let $ \big(\eta_i,f_i \big),  i=1,2$  be   two  $(t,p)-$parameter sets
such that $f_1$ is  Lipschitz continuous in $(y,z)$,
  and let  $ \big(Y^i ,Z^i   \big) \in  \hG^p_\bF([t,T])   $,  $i=1,2$ be a solution
   of BSDE$ \big(t,\eta_i,f_i \big)$.

 \no \(1\)  If    $ E \Big[    \Big( \neg \int_t^T   \neg  \big| f_1 (s, Y^2_s ,Z^2_s )
   -    f_2 (s, Y^2_s, Z^2_s) \big| ds    \Big)^{\neg {\wt{p}}}  \Big] < \infty$ for some ${\wt{p}} \neg \in \neg (1,p]$, then
   it holds \pas ~ that
\bea      \label{eq:n211}
     E \bigg[ \underset{s \in [t,T]}{\sup} \big|    Y^1_s \neg - \neg  Y^2_s   \big|^{\wt{p}} \bigg|\cF_t \bigg]
       \neg \le   \neg  C(T,{\wt{p}}, \g )     E \bigg[  \big|   \eta_1   \neg -  \neg  \eta_2    \big|^{\wt{p}}
           \neg +  \neg     \bigg( \neg \int_t^T   \neg  \big| f_1 (s, Y^2_s ,Z^2_s )
 \neg - \neg  f_2 (s, Y^2_s, Z^2_s) \big|   ds    \bigg)^{\neg {\wt{p}}} \bigg|\cF_t \bigg]     .
  \eea

 \no \(2\) If $ \eta_1 \le (\hb{resp.}\,\ge)\, \eta_2    $, \pas ~ and if     $  f_1 (s,  Y^2_s ,Z^2_s )
 \le (\hb{resp.}\,\ge)\, f_2 (s, Y^2_s ,Z^2_s )  $,  $ds \times dP-$a.s. on
  $[t,T] \times \O$,   then  it holds $P-$a.s. that   $ Y^1_s \le (\hb{resp.}\,\ge)\, Y^2_s$  for any $s \in [t,T]$.

\end{prop}

 \section{Stochastic Differential Games with Super-square-integrable Controls} \label{sec:zs_drgame}

    Let  $(\hU, \rho_{\overset{}{\hU}})$ and $(\hV, \rho_{\overset{}{\hV}})$ be two  separable metric spaces.
       For some $u_0 \in \hU$ and $v_0 \in \hV$, we define
  \beas
   [u]_{\overset{}{\hU}} \dfnn \rho_{\overset{}{\hU}}(u,u_0), ~\fa u \in \hU \q \hb{and} \q  [v]_{\overset{}{\hV}} \dfnn \rho_{\overset{}{\hV}}(v,v_0) , ~\fa v \in \hV .
  \eeas

 We shall study  a    zero-sum stochastic differential game  between two players,
   player I and player II,  who  choose super-square-integrable $\hU-$valued controls
    and   $\hV-$valued controls respectively  to compete:

  \begin{deff}   \label{def:controls}
 Given $t \in [0,T]$, an admissible control process $\mu=\{\mu_s\}_{s \in [t,T]}$  for player I   over   period $[t,T]$
 is a $\hU-$valued, $\bF-$progressively measurable process    such that
 $ E  \int_t^T [ \mu_s ]^q_{\overset{}{\hU}} \,ds < \infty$ for some $q > 2$.
 Admissible control processes    for player  II  over   period $[t,T]$ are defined similarly.
 We denote by $\cU_t$ \(resp.\;$\cV_t$\) the set of all admissible controls for player I \(resp.\;II\) over period $[t,T]$.

   \end{deff}

   \begin{rem} \label{rem_super_square_control}
   The reason why we use super-square-integrable controls lies in the fact that in the proof of Proposition \ref{prop_value_constant},  the set of $\hU-$valued \(resp.\,$\hV-$valued\) square integrable processes is
   not closed under Girsanov transformation via functions of  the Cameron-Martin space
    \(see in particular \eqref{eq:xvx319}\).
   \end{rem}

    Clearly,  connecting two $\cU_t-$controls along some $\t \in \cS_{t,T} $ results in a new $\cU_t- $control:

   \begin{lemm}  \label{lem_control_combine}
Let $t \in [0,T]$ and $\t \in \cS_{t,T} $. For any $\mu^1, \mu^2 \in \cU_t$,
 $\mu_s \dfnn  \b1_{\{ s < \t   \}}  \mu^1_s
     \neg + \neg  \b1_{\{ s \ge \t     \}} \mu^2_s $, $ s  \in [t,T] $ defines a  $\cU_t-$control.
  Similarly, for any $\nu^1, \nu^2 \in \cV_t$,
 $\nu_s \dfnn  \b1_{\{ s < \t   \}}  \nu^1_s
     \neg + \neg  \b1_{\{ s \ge \t     \}} \nu^2_s $, $ s  \in [t,T] $ defines a  $\cV_t-$control.

\end{lemm}

    \subsection{Game Setting: A Controlled SDE$-$BSDE  System}

   Our  zero-sum stochastic differential game  is  formulated via a (decoupled) SDE$-$BSDE  system
   with the following parameters: Fix $k \in \hN$, $\g >0$ and $p \in (1,2]$.

  \ss  \no 1)   Let  $b: [0,T] \times  \hR^k  \times \hU \times \hV  \to  \hR^k  $
  be a $\sB([0,T]) \otimes \sB(\hR^k) \otimes \sB(\hU) \otimes \sB(\hV)/\sB(\hR^k)-$measurable
  function and let $\si: [0,T] \times  \hR^k \times \hU \times \hV  \to  \hR^{k \times d} $
  be a $\sB([0,T]) \otimes \sB(\hR^k) \otimes \sB(\hU) \otimes \sB(\hV)/\sB(\hR^{k \times d})-$measurable function  such  that for any $(t,u,v)  \neg \in \neg  [0,T]  \neg \times \neg  \hU  \neg \times \neg  \hV $
    and $x,x' \in      \hR^k $
   \bea
   && \hspace{2.05cm} |b(t,0,u,v)  |     +    |\si(t,0,u,v) |
       \neg \le  \neg  \g \big(1 + [u]_{\overset{}{\hU}} + [v]_{\overset{}{\hV}}  \big)   \label{b_si_linear_growth}   \\
    \hb{and}  && |b(t,x,u,v) \neg - \neg b(t,x',u,v)|    +    |\si(t,x,u,v) \neg - \neg \si(t,x',u,v) |
       \neg \le  \neg  \g |x-x'|      .     \qq     \label{b_si_Lip}
   \eea

     \no 2) Let $  g     : \hR^k  \neg \to \neg  \hR $ be  a $ 2/p-$H\"older continuous function with coefficient $\g$.

  \ss \no 3)  Let   $f: [0,T] \times \hR^k \times \hR  \times \hR^d \times \hU \times \hV \to \hR$  be
  $\sB([0,T]) \otimes \sB(\hR^k)  \otimes \sB(\hR)  \otimes \sB(\hR^d)  \otimes \sB(\hU) \otimes \sB(\hV)/\sB(\hR)-$measurable
  function such that
      for any  $(t,  u,v  ) \in [0,T]   \times \hU \times \hV $ and any $( x,y,z), ( x',y',z' ) \in \hR^k  \times \hR \times \hR^d$
     \bea
   |f(t,0,0,0,u,v)  |   & \dneg  \dneg  \le & \dneg \dneg  \g \Big(1 + [u]^{2/p}_{\overset{}{\hU}}
 + [v]^{2/p}_{\overset{}{\hV}}  \Big)   \label{f_linear_growth} \\
 \hb{and} \q \big|f(t,x,y,z,u,v)-f(t,x',y',z',u,v)\big| & \dneg \dneg\le & \dneg \dneg   \g \big(|x-x'|^{2/p} + |y-y'| + |z-z'|\big)   .  \qq  \label{f_Lip}
 \eea

  For any $\l  \ge  0$,    we let  $c_\l$ denote a generic constant,  depending on $\l$, $T$, $\g $,  $p$ and $|g(0)|$,
        whose form  may vary from line to line. (In particular, $c_0$ stands for a generic constant depending
      on $ T $, $\g $, $p$  and $|g(0)|$.)

 \ss Also, we would like to introduce two control neutralizers $\psi , \wt{\psi}$ for the coefficients: For some $\k > 0$

   \ss \no {\bf (A-u)}  there exist  a function $\psi : [0,T] \times \big( \hU \backslash O_\k(u_0) \big) \to \hV$ that is $ \sB([0,T]) \times \sB \big( \hU \backslash O_\k(u_0) \big) / \sB(\hV) -$measurable and satisfies: for any $(t,x,y,z)  \neg \in \neg  [0,T]
     \neg \times \neg  \hR^k \neg \times \neg  \hR   \neg \times \neg  \hR^d  $ and $u,u' \in  \hU \backslash O_\k(u_0) $
   \beas
 &&  b \big(t,x, u , \psi(t,u)\big) \neg = \neg b \big(t,x, u' , \psi(t,u') \big)  , \q  \si \big(t,x, u , \psi (t,u)\big) \neg = \neg \si \big(t,x, u' , \psi(t,u') \big)  , \\
 &&    f \big(t,x,y,z,u , \psi(t,u)\big)    \neg = \neg f \big(t,x,y,z, u' , \psi(t,u') \big)    \q       \hb{and} \q  [\psi(t,u)]_{\overset{}{\hV}}  \neg \le \neg  \k (1+ [u]_{\overset{}{\hU}})  ;
   \eeas

  \ss \no {\bf (A-v)} and there exists a function $\wt{\psi} : [0,T] \times \big( \hV \backslash O_\k(v_0) \big) \to \hU$ that is $ \sB([0,T]) \times \sB \big( \hV \backslash O_\k(v_0) \big) / \sB(\hU) -$measurable and satisfies: for any $(t,x,y,z)  \neg \in \neg  [0,T]
     \neg \times \neg  \hR^k   \neg \times \neg  \hR \neg \times \neg  \hR^d    $ and $v,v' \in  \hV \backslash O_\k(v_0) $ \beas
 &&  b \big(t,x, \wt{\psi}(t,v),v \big) \neg = \neg b \big(t,x,\wt{\psi}(t,v'),v' \big)  , \q  \si \big(t,x, \wt{\psi}(t,v),v\big) \neg = \neg \si \big(t,x, \wt{\psi}(t,v'),v' \big)  , \\
 &&    f \big(t,x,y,z,\wt{\psi}(t,v),v\big)    \neg = \neg f \big(t,x,y,z, \wt{\psi}(t,v'),v' \big)    \q       \hb{and} \q  [\wt{\psi}(t,v)]_{\overset{}{\hU}}  \neg \le \neg  \k (1 \neg + \neg [v]_{\overset{}{\hV}} ) .
   \eeas

   A typical example   satisfying  both \(A-u\) and \(A-v\) is the additive-control case:

 \begin{eg}  \label{eg_control_neutralizer}
 Let $\hU=\hV = \hR^\ell$ and consider the following coefficients:
 \beas
 && b  (t,x, u , v  ) \neg = \neg b \big(t,x, u+v \big)  , \q  \si  (t,x, u , v  ) \neg = \neg \si \big(t,x, u+v \big)
 \q \hb{and} \\
    && f (t,x, y,z, u , v  ) \neg = \neg f \big(t,x,y,z, u+v \big), \q \fa (t,x,  y,z, u , v) \in [0,T]
     \neg \times \neg  \hR^k   \neg \times \neg  \hR \neg \times \neg  \hR^d \neg \times \neg  \hU \neg \times \neg  \hV .
 \eeas
 Then \(A-u\) and \(A-v\) hold for functions $\psi (u) =-u $  and $\wt{\psi} (v) = -v$  respectively.

 \end{eg}

  Here is another example:

     \begin{eg}  \label{eg_example2}
   Given $\g > 0$,       let  $b_0,\si_0: [0,T] \times  \hR      \to  \hR   $
  be two $\sB([0,T]) \otimes \sB(\hR )  /\sB(\hR )-$measurable
  functions and   let $f_0: [0,T] \times \hR  \times \hR  \times \hR    \to \hR$  be
  $\sB([0,T]) \otimes \sB(\hR )  \otimes \sB(\hR)  \otimes \sB(\hR )   /\sB(\hR)-$measurable
  function such that    for any $ t   \neg \in \neg  [0,T]    $
    and   $( x,y,z), ( x',y',z' ) \in \hR  \times \hR \times \hR$
   \beas
      |b_0(t,x ) \neg - \neg b_0(t,x' )|    +    |\si_0(t,x ) \neg - \neg \si_0(t,x' ) |
    + \big|f_0(t,x,y,z )-f_0(t,x',y',z' )\big|  \neg \le  \neg  \g \big( |x-x'| + |y-y'| + |z-z'|\big)     .
   \eeas
Also,  let $\hU=\hV = \hR $,  $ \k > 0$ and   $\vf: [0,T] \times    \hU \times \hV     \to  \hR   $
  be a jointly continuous  function such that $\vf$ is Lipschitz continuous in $(u,v)$ with coefficient $\g$,
  $ \underset{t \in [0,T]}{\sup} |\vf(t,0,0)| \le \g $, and
 for any $ (t,u,v)  \neg \in \neg  [0,T] \neg \times \neg     \hU  \neg  \times \neg   \hV$  \bea  \label{eq:xqxqx011}
  \underset{|v'| \le \k |u|}{\inf} \vf (t,u,v') \le 0 \le \underset{|v'| \le \k |u|}{\sup} \vf (t,u,v')    \q
  \hb{and} \q
  \underset{|u'| \le \k |v|}{\inf} \vf (t,u',v) \le 0 \le \underset{|u'| \le \k |v|}{\sup} \vf (t,u',v)   .
  \eea
Then $ b(t,x,u,v) \dfnn b_0(t,x) + \vf(t,u,v)  $, $ \si(t,x,u,v) \dfnn \si_0(t,x) + \vf(t,u,v)  $
 and $ f(t,x,y,z,u,v) \dfnn f_0(t,x,y,z) + \vf(t,u,v)  $, $\fa (t,x,y,z,u,v) \in [0,T] \times \hR \times \hR \times \hR
 \times \hU \times \hV$ are the measurable functions satisfying \eqref{b_si_linear_growth}$-$\eqref{f_Lip}
 with  $k=d=1$ and $p=2$. We will show at the beginning of Subsection \ref{subsection:Proofs_S2} that
   \(A-u\) and \(A-v\) hold for these coefficients.

 \end{eg}

  \ms When the game begins at time $t \in [0,T]$,     player I and player II select admissible controls
  $\mu \in \cU_t$ and $\nu \in \cV_t$ respectively. Then the state process starting from  $\xi  \in \hL^2 (\cF_t, \hR^k)$
  will  evolve   according to   SDE \eqref{FSDE} on  the  probability space $\big( \O,  \cF_T,  P \big)$.
      The measurability of functions $b$, $\si$, $\mu$ and $\nu$ implies that
  \beas
      b^{\mu,\nu} (s,\o, x) \neg \dfnn \neg     b \big(s,  x,\mu_s(\o),\nu_s (\o) \big) ,
      \q \fa (s,\o,x) \neg \in  \neg  [t,T]    \neg \times  \neg  \O \neg \times \neg \hR^k
  \eeas
  is $\sP \otimes \sB(\hR^k)  /\sB(\hR^k)-$measurable  and that
    \beas
     \si^{\mu,\nu}   (s,\o, x)  \neg \dfnn \neg         \si \big(s,  x,\mu_s(\o),\nu_s (\o) \big)  ,
     \q \fa (s,\o,x) \neg \in  \neg  [t,T]    \neg \times  \neg  \O \neg \times \neg \hR^k
   \eeas
   is  $\sP \otimes \sB(\hR^k)  /\sB(\hR^{k \times d})-$measurable.
   Also,   \eqref{b_si_Lip}, \eqref{b_si_linear_growth} and H\"older's inequality show that    $b^{\mu,\nu}   $
   , $\si^{\mu,\nu}   $ are Lipschitz continuous in $x$ and   satisfy
  \beas
    E   \left[ \Big(\int_t^T  \neg \big| b^{\mu,\nu}    (s, 0 )   \big| \, ds \Big)^2
   + \Big(\int_t^T  \neg \big| \si^{\mu,\nu}    (s, 0 )   \big| \, ds \Big)^2 \right]
   \le c_0 + c_0 E \neg \int_t^T  \neg \big([ \mu_s ]^2_{\overset{}{\hU}} \neg+\neg [ \nu_s ]^2_{\overset{}{\hV}} \big) ds < \infty .
   \eeas
  Then   it is well-known that
 the SDE \eqref{FSDE} admits      a  unique  solution  $\left\{ X^{t,\xi,\mu,\nu}_s \right\}_{s \in [t,T]}
    \in \hC^2_{\bF }([t,T], \hR^k)$ such that
  \bea
   E \neg \left[  \underset{s \in [t,T]}{\sup}  \big| X^{t,\xi,\mu,\nu}_s  \big|^2 \right]
  & \dneg \le &  \dneg
      c_0    E\big[ |\xi|^2 \big] \neg + \neg c_0   E   \left[ \Big(\int_t^T  \neg \big| b^{\mu,\nu}     (s, 0 )   \big| \, ds \Big)^2
   + \Big(\int_t^T  \neg \big| \si^{\mu,\nu}     (s, 0 )   \big| \, ds \Big)^2 \right] \nonumber  \\
    & \dneg  \le  & \dneg
           c_0  \bigg( 1  \neg + \neg  E \big[|\xi|^2 \big]  \neg + \neg  E  \neg \int_t^T  \dneg \big( [ \mu_s ]^2_{\overset{}{\hU}} \neg + \neg  [ \nu_s ]^2_{\overset{}{\hV}}  \big)  ds  \bigg) < \infty . \q   \label{eq:esti_X_1a}
 \eea

    Given $s \neg \in \neg  [t,T]$, let $ [\mu]^s  $ denote  the restriction of $\mu$ over  period
  $[s,T]  \neg  : i.e., [\mu]^s_r  \neg \dfnn \neg  \mu_r $, $\fa r  \neg \in \neg  [s,T]$.
  Clearly, $[\mu]^s  \neg \in \neg  \cU_s $, similarly, $ \big\{[\nu]^s_r \neg \dfnn \neg  \nu_r \big\}_{ r  \in   [s,T]}    \neg \in \neg  \cV_s $. As
  \beas
    X^{t,\xi,\mu,\nu}_r  &=& X^{t,\xi,\mu,\nu}_s+ \int_s^r b \big(r', X^{t,\xi,\mu,\nu}_{r'},\mu_{r'},\nu_{r'} \big) \, dr'
   + \int_s^r \si \big(r', X^{t,\xi,\mu,\nu}_{r'},\mu_{r'},\nu_{r'} \big) \, dB_{r'} \\
 &=& X^{t,\xi,\mu,\nu}_s+ \int_s^r b \big(r', X^{t,\xi,\mu,\nu}_{r'},[\mu]^s_{r'}, [\nu]^s_{r'} \big) \, dr'
   + \int_s^r \si \big(r', X^{t,\xi,\mu,\nu}_{r'},[\mu]^s_{r'},[\nu]^s_{r'} \big) \, dB_{r'} , \q   r \in [s,T] ,
 \eeas
 we see that $\big\{X^{t,\xi,\mu,\nu}_r\big\}_{r \in [s,T]} \in \hC^2_{\bF }([s,T], \hR^k)$ solves \eqref{FSDE} with the parameters
   $   \big(s, X^{t,\xi,\mu,\nu}_s, [\mu]^s , [\nu]^s \big)$. To wit,  it holds \pas ~ that
   \bea  \label{eq:xxa603}
   X^{t,\xi,\mu,\nu}_r = X^{s, X^{t,\xi,\mu,\nu}_s, [\mu]^s , [\nu]^s}_r , \q  \fa r \in [s,T].
   \eea




 Moreover, the state process  depends on controls in the following way:

 \begin{lemm} \label{lemm_control_compare}
 Given $ t   \in [0,T]$, let $  \xi \in   \hL^2(\cF_t, \hR^k)$ and $( \mu,  \nu ), \big( \wt{\mu}, \wt{\nu} \big) \in \cU_t \neg \times \neg  \cV_t $.
 If      $ \big(  \mu,  \nu \big)  \neg = \neg  \big( \wt{\mu}, \wt{\nu} \big)  $, $dr  \neg \times \neg  d P-$a.s.
  on $\[t,\t\[ \, \cup \, \[\t,T\]_A$ for some $\t \neg \in \neg   \cS_{t,T}$ and $A \in \cF_\t$, then it holds \pas ~ that
 \bea  \label{eq:p611}
  \b1_A X^{t,\xi,\mu,\nu}_s   + \b1_{A^c} X^{t,\xi,\mu,\nu}_{\t \land s}
  = \b1_A  X^{t,\xi,\wt{\mu}, \wt{\nu}}_s   + \b1_{A^c} X^{t,\xi,\wt{\mu}, \wt{\nu}}_{\t \land s}  ,
  \q   \fa   s \in [t, T ] .
 \eea
 \end{lemm}

  Now,  let $\Th$ stand for the quadruplet $(t,\xi,\mu,\nu)$.
  Given  $\t \in \cS_{t,T}$,   the measurability of $ ( f ,  X^\Th , \mu, \nu ) $  and  \eqref{f_Lip} imply that
   \beas
    \q   f^\Th_\t   (s,\o,y,z ) \dfnn \b1_{\{    s < \t(\o)\}}  f  \Big(s,  X^\Th_s(\o) , y,  z,
        \mu_s (\o),   \nu_s (\o) \Big) ,
    \q   \fa (s,\o,y,z  ) \in   [t,T] \times \O \times \hR \times \hR^d
    \eeas
    is a    $\sP   \otimes \sB(\hR)  \otimes \sB(\hR^d)/\sB(\hR)-$measurable function
      that is Lipschitz continuous in $(y,z)$ with coefficient $\g$.
  And  one can deduce from   \eqref{f_linear_growth}, \eqref{f_Lip} and  H\"older's inequality    that
 \bea    \label{eq:s031}
      E   \bigg[      \Big(\int_t^T  \neg \big| f^\Th_\t   (s,0, 0 )   \big|  ds \Big)^p\bigg]
   \neg  \le  \neg  c_0    \neg + \neg  c_0   E \bigg[   \underset{s \in [t,T]}{\sup}   \big| X^\Th_s \big|^2
   \neg + \neg \int_t^T \dneg \big( [ \mu_s ]^2_{\overset{}{\hU}} \neg+\neg [ \nu_s ]^2_{\overset{}{\hV}}  \big)  ds \bigg]
    \neg < \neg  \infty .   \q
   \eea
    Thus, for any $\eta \in \hL^p(\cF_{\t}) $,  Proposition \ref{BSDE_well_posed} shows that  the BSDE$(t,\eta,  f^\Th_\t )$
 admits a unique solution $ \big(Y^\Th (\t,\eta),  Z^\Th (\t,\eta)  \big) \\ \neg \in \neg \hG^p_\bF\big([t,T]\big) $, which has the following estimate as a consequence of \eqref{eq:n211}.

     \begin{cor} \label{cor_estimate_Y0}

 Let     $ t \neg \in \neg  [0,T]$, $\xi  \neg  \in \neg  \hL^2(\cF_t,\hR^k)$,  $( \mu,\nu)  \neg \in  \neg   \cU_t  \neg \times \neg  \cV_t$ and $\t  \neg \in \neg  \cS_{t,T}$.
 Given  $\eta_1,  \eta_2  \neg \in \neg  \hL^p(\cF_{\t}) $, it holds for any $\wt{p}  \neg \in \neg  (1,p]$ that
   \bea
      E \bigg[ \, \underset{s \in [t,T]}{\sup} \Big|  Y^{t,\xi,\mu,\nu}_s   (\t, \eta_1  )
     - Y^{t, \xi, \mu, \nu}_s   (\t,  \eta_2  ) \Big|^{\wt{p}} \bigg|\cF_t \bigg]
    \le c_{\wt{p}}       E \big[     | \eta_1 -  \eta_2  |^{\wt{p}}   \big|\cF_t    \big] , \q  \pas     \label{eq:xvx071}
     \eea

 \end{cor}

 \ss      Given another  stopping time  $\z \neg \in \neg  \cS_{t,T}  $ with $\z  \neg \le \neg  \t  $, $P-$a.s.,  one can easily show that
 $  \Big\{  \neg  \Big( Y^\Th_{\z \land  s}   (\t,\eta),  \b1_{\{s < \z \}}Z^\Th_s   (\t,\eta)   \Big)
   \neg  \Big\}_{s \in [t,T]} \\
  \in \hG^p_\bF\big([t,T]\big)$
 solves the BSDE$\big(t,  Y^\Th_\z  (\t,\eta),f^\Th_\z    \big)$. To wit, we have
 \bea
  \Big(Y^\Th_s \big(\z,  Y^\Th_\z  (\t,\eta)\big) ,  Z^\Th_s \big(\z,  Y^\Th_\z  (\t,\eta)\big)  \Big)
    =  \Big( Y^\Th_{\z \land  s}   (\t,\eta), \b1_{\{s < \z \}}Z^\Th_s   (\t,\eta)   \Big), \q s \in [t,T].   \label{eq:p677}
 \eea
 In particular, when $\z =\t$,
 \bea   \label{eq:xux313}
  \Big(Y^\Th_s (\t,\eta) ,  Z^\Th_s (\t,\eta)   \Big)
    \neg = \neg  \Big( Y^\Th_{\t \land  s}   (\t,\eta), \b1_{\{s < \t \}}Z^\Th_s   (\t,\eta)    \Big), \q  s \in [t,T].
 \eea

   On the other hand, if $\t \in \cS_{s,T}$ for some   $s \neg \in \neg  [t,T]$,
 letting $\Th^s \dfnn \big( s, X^\Th_s, [\mu]^s , [\nu]^s \big)$,  we can deduce from \eqref{eq:xxa603} that
   $\big\{\big(Y^\Th_r (\t,\eta)  ,Z^\Th_r (\t,\eta)   \big) \big\}_{r \in [s,T]} \neg \in \neg \hG^p_\bF\big([s,T]\big) $ solves
    the following BSDE$\big(s,\eta,  f^{\Th^s}_\t \big)$:
   \beas
       Y_s &=& \eta  \neg + \neg  \int_r^T  \neg \b1_{\{r' < \t \}}
        f (r', X^\Th_{r'}, Y_{r'},   Z_{r'}, \mu_{r'}, \nu_{r'})  \, dr'
     \neg-\neg \int_r^T \neg Z_{r'} d B_{r'}             \\
  & = & \eta  \neg + \neg  \int_r^T  \neg \b1_{\{r' < \t \}} f (r',X^{\Th^s}_{r'}, Y_{r'},   Z_{r'}, [\mu]^s_{r'}, [\nu]^s_{r'})  \, dr'   \neg-\neg \int_r^T \neg Z_{r'} d B_{r'}  , \q    r \in [s,T]    .
    \eeas
  Hence,  it holds  \pas ~ that
 \bea  \label{eq:xxa605}
 Y^\Th_r  (\t,\eta ) = Y^{\Th^s}_r (\t,\eta), \q \fa  r \in [s,T].
 \eea

 The $ 2/p-$H\"older continuity of functions $ g$ and \eqref{eq:esti_X_1a} show that $ g\big( X^\Th_T \big) \in \hL^p(\cF_T)$.
 Set   $  J(\Th) \dfnn  Y^\Th_t   \big(T,g \big(X^\Th_T \big)  \big)      $
 From \eqref{eq:n211} and the standard estimate  of SDE \eqref{FSDE}, we can deduce the following a priori estimate:
  \begin{lemm} \label{lem_estimate_Y}
 Let     $\, t \in [0,T]$   and $( \mu,\nu) \in   \cU_t \times \cV_t$.
 Given   $\xi_1, \xi_2 \in \hL^2(\cF_t,\hR^k)$, it holds for any $ \wt{p} \in (1, p] $ that
  \bea
      E \bigg[ \, \underset{s \in [t,T]}{\sup} \Big|  Y^{t,\xi_1,\mu,\nu}_s  \big(T, g\big(X^{t,\xi_1,\mu,\nu}_T \big) \big)
     - Y^{t,\xi_2,\mu,\nu}_s  \big(T, g\big(X^{t,\xi_2,\mu,\nu}_T \big) \big) \Big|^{\wt{p}} \bigg|\cF_t \bigg]
    \le c_{\wt{p}}   | \xi_1 \neg - \neg  \xi_2  |^{\frac{2 \wt{p}}{p}} , \q  \pas     \label{eq:s025}
     \eea

\end{lemm}

  \subsection{Definition of  the Value Functions and a Weak Dynamic Programming Principle}

\ss  Now, we are ready to introduce   values of the zero-sum stochastic differential games via the following
  version of Elliott$-$Kalton strategies (or non-anticipative strategies).  

\begin{deff}   \label{def:strategy}
  Given   $t  \neg \in \neg  [0,T]$,    an admissible strategy $\a$
  for player I  over  period   $[t,T]$ is a mapping $\a \neg: \cV_t \neg \to \neg \cU_t $ satisfying:
     \(i\)  There exists a $C_\a \neg > \neg  0$ such that  for  any $\nu  \neg \in \neg  \cV_t$
   $   
    \big[ \big( \a (\nu)\big)_s   \big]_{\overset{}{\hU}}    \neg   \le  \neg   \k  \neg + \neg  C_\a [\nu_s ]_{\overset{}{\hV}} $,
    $ ds  \neg \times \neg  dP- $a.s.,
   where $\k$ is the constant that appears in \(A-u\) and \(A-v\); \,
    \(ii\)  For any    $\nu^1, \nu^2  \neg \in \neg  \cV_t$, $ \t  \neg \in \neg    \cS_{t,T}$
  and $ A   \neg \in \neg  \cF_\t  $, if
 $\nu^1  \neg = \neg   \nu^2 $,  $ds  \neg \times \neg  d P -$a.s.  on $\[t,\t\[  \, \cup  \,   \[\t, T\]_A   $,
   then $\a (\nu^1) \neg = \neg \a (\nu^2)$,   $ds  \neg \times \neg  d P -$a.s.
    on $\[t,\t\[  \,   \cup  \,   \[\t, T\]_A   $.

 Admissible strategies $\beta: \cU_t \to \cV_t$  for player II  over period $[t,T]$ are defined similarly.    The collection of all  admissible strategies   for player I \(resp.\;II\) over   period $[t,T]$ is denoted by $\fA_t$ \big(resp.\;$\fB_t $\big).
\end{deff}

 \begin{rem} \label{rem_nonanticipate}
  The condition \(ii\) of Definition \ref{def:strategy} is called the nonanticipativity of strategies.
   It is said in  \cite[line 4 of page 456]{Buckdahn_Li_1} that
   ``From the nonanticipativity of $\beta_2$ we have $\beta_2(u^\e_2) = \sum_{j \ge 1} \b1_{\D_j} \beta_2(u^2_j)$, $\cds$".
   What actually used in this equality is not the  nonanticipativity of $\beta_2$ as defined in Definition 3.2 therein, but
   the  requirement:
   \bea  \label{eq:xvx531}
  \hb{For any    $u, \wt{u}  \neg \in \neg  \cU_{t+\d,T}$
  and $ A   \neg \in \neg  \cF_{t+\d}  $, if
 $u    \neg = \neg   \wt{u} $   on $  [t+\d , T]  \neg \times \neg  A   $,
   then $\beta_2 (u )  \neg = \neg  \beta_2(\wt{u}) $
    on $ [t \neg + \neg \d , T]  \neg \times \neg  A    $.}
   \eea
   Since $\beta_2$ is a restriction of strategy $\beta \in \cB_{t,T}$ over period $[t \neg + \neg \d,T]$,
   \eqref{eq:xvx531} entails the following condition on $\beta$.
  \beas
  && \hb{For any    $u, \wt{u}  \neg \in \neg  \cU_{t,T}$, any $ s \in [t,T]$
  and any $ A   \neg \in \neg  \cF_s  $, if
 $u    \neg = \neg   \wt{u} $   on $ \big([t,s) \times \O\big) \cup \big( [s , T]  \neg \times \neg  A \big)   $,}  \\
 && \hb{then $\beta (u )  \neg = \neg  \beta (\wt{u}) $
    on $  \big([t,s) \times \O\big) \cup \big( [s , T]  \neg \times \neg  A \big)    $.}
   \eeas
 which is exactly a simple    version of our nonanticipativity condition on strategies with $\t = s$.

 \end{rem}

      For any $(t,x) \in  [0,T] \times  \hR^k $,    
   we define
  \beas
   && w_1 (t,x) \dfnn \underset{\beta \in \fB_t }{\essinf} \; \underset{\mu \in \cU_t}{\esssup} \;
    J \big(t,x,\mu, \beta (\mu)\big)   =  \underset{\beta \in \fB_t }{\essinf} \; \underset{\mu \in \cU_t}{\esssup} \; Y^{t,x,\mu, \beta (\mu) }_t   \Big(T,g \Big(X^{t,x,\mu, \beta ( \mu ) }_T \Big)  \Big) \qq \\
  \hb{and} &&   w_2 (t,x )  \dfnn   \underset{\a \in  \fA_t}{\esssup} \, \underset{\nu \in \cV_t}{\essinf} \; J \big(t,x,\a (  \nu ) ,\nu \big)   =  \underset{\a \in  \fA_t}{\esssup} \, \underset{\nu \in \cV_t}{\essinf} \; Y^{t,x,\a (  \nu ) ,\nu}_t  \Big(T, g \Big(X^{t,x,\a (  \nu ) ,\nu}_T  \Big) \Big)
  \eeas
as    player I's and player II's   {\it priority  values}   of the zero-sum stochastic differential  game that starts   from time $t  $ with initial  state $x$.

 \begin{rem}
 When  $f$ is independent of $(y,z) $,  $w_1$ and $w_2$ are in form of
  \beas
   && \hspace{-5mm} w_1 (t,x) \dfnn \underset{\beta \in \fB_t }{\essinf} \; \underset{\mu \in \cU_t}{\esssup} \;
      E \bigg[ g\big(X^{t,x,\mu, \beta (\mu)}_T \big)
    + \int_t^T f \big(s, X^{t,x,\mu, \beta (\mu)}_s, \mu_s , (\beta (\mu))_s \big) ds \bigg| \cF_t \bigg]    \qq \\
  \hb{and} && \hspace{-5mm} w_2 (t,x )  \dfnn   \underset{\a \in  \fA_t}{\esssup} \, \underset{\nu \in \cV_t}{\essinf} \;
   E \bigg[ g\big(X^{t,x,\a (  \nu ) ,\nu}_T \big)
    + \int_t^T f \big(s, X^{t,x,\a (  \nu ) ,\nu}_s, ( \a (  \nu ))_s ,\nu_s \big) ds \bigg| \cF_t \bigg] ,
     ~ \;   \fa   (t,x) \in  [0,T] \times  \hR^k .
  \eeas
 \end{rem}

 \begin{rem}   \label{rem_BL_case}
 When $U$ and $V$ are compact \(say $\hU = \ol{O}_\k(u_0)$ and $\hV = \ol{O}_\k(v_0)$\),
 Assumptions \(A-u\), \(A-v\) 
 are no longer needed, and the integrability condition
 in Definition \ref{def:controls} as well as  the condition \(i\) in Definition
 \ref{def:strategy} hold automatically.  Thus
  our game problem degenerates to the case of  \cite{Buckdahn_Li_1}.
 \end{rem}

  Let us    review  some basic properties of the essential extrema for the later use \big(see e.g.
\cite[Proposition VI-\b1-1]{Neveu_1975} or \cite[Theorem A.32]{Follmer_Schied_2004}\big):

 \begin{lemm}
 \label{lem_ess}
Let $\{\xi_i\}_{i \in \cI}$,   $\{\eta_i\}_{i \in \cI}$ be two
classes of $\cF_T-$measurable random variables with the same index set
$\cI$.

\ss \no  (1)  If $\xi_i \le (=)\; \eta_i$, \pas~ holds for all $i \in
\cI$,
 then $\underset{i \in \cI}{\esssup}\, \xi_i \le (=)\; \underset{i \in \cI}{\esssup}\, \eta_i$, \pas

\ss \no  (2)  For any $A \in \cF_T$, it holds \pas~that
 $\;\underset{i \in \cI}{\esssup}\, \big( \b1_A \xi_i + \b1_{A^c} \eta_i
\big) = \b1_A \,\underset{i \in \cI}{\esssup}\, \xi_i + \b1_{A^c}\,
\underset{i \in \cI}{\esssup}\, \eta_i$. In particular, $\underset{i
\in \cI}{\esssup}\, \big( \b1_A \xi_i  \big) =  \b1_A\, \underset{i
\in \cI}{\esssup}\, \xi_i $, \pas

\ss \no (3) For any $\cF_T-$measurable random variable $\eta $ and any
 $\l>0$, we have $\underset{i \in \cI}{\esssup}\,  (\l \xi_i + \eta )
 = \l \, \underset{i \in \cI}{\esssup}\, \xi_i + \eta $, \pas

\ss \no  (1)-(3) also hold when we replace $\,\underset{i \in
\cI}{\esssup}\,$  by $\,\underset{i \in \cI}{\essinf}\,$.
\end{lemm}

 The  values $w_1$, $w_2$ are bounded as follows:
  \begin{prop} \label{prop_value_bounds}
  For any $(t,x) \in [0,T] \times \hR^k$, it holds \pas ~ that $|w_1 (t,x)| + |w_2 (t,x)| \le c_\k  + c_0 |x|^{2/p}$.

  \end{prop}

  Similar to  Proposition 3.1 of \cite{Buckdahn_Li_1}, the following result allows us to regard
    $ w_1  $ and $ w_2  $ as     deterministic functions on $[0,T] \times \hR^k$:
  \begin{prop} \label{prop_value_constant}
 Let $i=1,2$.  For any $(t,x) \in [0,T] \times \hR^k$, it holds \pas ~ that  $ w_i(t,x)=  E[w_i(t,x)]  $.

  \end{prop}

 Moreover,  as a consequence of \eqref{eq:s025},  $w_1$ and $w_2$ are $ 2/p-$H\"older continuous in $x$:

 \begin{prop}  \label{prop_w_conti}
 For any $t \in [0,T] $ and $ x_1,x_2 \in \hR^k$,
 $   \big| w_1(t,x_1) \neg - \neg  w_1(t,x_2) \big|    +\big| w_2(t,x_1) \neg - \neg  w_2(t,x_2) \big|
   \le c_0 | x_1 \neg - \neg x_2 |^{2/p} $.

 \end{prop}

  However, the values  $w_1$,  $w_2$  are generally not 
  continuous in $t$ unless   $\hU$, $\hV$ are compact. 

 \begin{rem}  \label{rem_measurability_issue}

   When trying to directly prove the dynamic programming principle, 
 \cite{Fleming_1989} encountered a measurability issue:
 The pasted strategies for approximation may not be progressively measurable,  see  page 299 therein. So
    they  first proved that the value functions
 are unique viscosity solutions to the associated Bellman-Isaacs equations
 by a   time-discretization approach \(assuming that the limiting Isaacs equation has a comparison principle\), which relies on the following regularity of the approximating  values $v_\pi$
   \beas
   |v_\pi(t,x)- v_\pi(t',x')| \le C \big( |t-t'|^{1/2} + |x-x'|\big), \q \fa (t,x), (t',x') \in [0,T] \times \hR^k
   \eeas
 with a uniform coefficient $C>0$  for all partitions $\pi$  of $[0,T]$.
 Since our value functions $w_1$, $w_2$ may  not be $1/2-$H\"older continuous in $t$, 
 this method seems not suitable for our problem.   Hence, we  adopt Buckdahn and Li's probability setting.

   \end{rem}


   \ss   The following  weak dynamic programming principle for value functions $w_1$, $w_2$ is the main result of the
   paper:

\begin{thm}   \label{thm_DPP}

  \no 1\)     Given   $ t \neg \in \neg [0,T)$,  let $\f, \wt{\f} \neg : [t,T]   \neg    \times  \neg
     \hR^k  \neg \to \neg  \hR $   be two  continuous functions such that
     $    \f (s,x)    \le     w_1 (s,x)   \le    \wt{\f} (s,x) $,
     $   (s,x)  \neg \in \neg  [t,T]  \neg \times \neg  \hR^k$. Then    for any
      $x  \neg \in \neg     \hR^k$ and $\d  \neg \in \neg  (0, T-t) $, it holds \pas ~ that
           $$
       \underset{\beta \in \fB_t }{\essinf} \; \underset{\mu \in \cU_t }{\esssup}    \;
    Y^{t,x,\mu, \beta (  \mu ) }_t \Big(\t_{\beta,\mu},
    \f \big(\t_{\beta,\mu},  X^{t,x,\mu, \beta (  \mu ) }_{\t_{\beta,\mu}} \big) \Big)
   \neg  \le  \neg  w_1(t,x) \le    \underset{\beta \in \fB_t }{\essinf} \; \underset{\mu \in \cU_t }{\esssup}    \;
    Y^{t,x,\mu, \beta (  \mu ) }_t \Big(\t_{\beta,\mu},
    \wt{\f} \big(\t_{\beta,\mu},  X^{t,x,\mu, \beta (  \mu ) }_{\t_{\beta,\mu}} \big) \Big) ,
    $$
 where    $\t_{\beta,\mu} \dfnn  \inf \big\{s  \neg \in \neg  ( t , T  ] \neg  : \big( s,X^{t,x,\mu, \beta (  \mu ) }_s \big)
  \neg  \notin   \neg   O_{  \d} (t, x) \big\}   $.

  \ss \no 2\) Given    $ t \neg \in \neg [0,T)$,   let $\f, \wt{\f} \neg : [t,T]     \times
     \hR^k  \neg \to \neg  \hR $   be two   continuous functions   such that    $   \f (s,x) \le w_2 (s,x)
      \neg \le \neg  \wt{\f} (s,x)      $, $   (s,x)
        \neg \in \neg [t,T] \neg \times \neg \hR^k$.   Then   for any
      $ x   \neg \in \neg     \hR^k$ and $\d  \neg \in \neg  (0, T-t) $,    it holds \pas ~ that
           $$
  \underset{\a \in \fA_t }{\esssup} \; \underset{\nu \in \cV_t }{\essinf}    \;
    Y^{t,x,  \a ( \nu ), \nu }_t \Big(\t_{\a,\nu},
     \f \big(\t_{\a,\nu},  X^{t,x, \a ( \nu ), \nu }_{\t_{\a,\nu}} \big) \Big) \neg \le \neg w_2(t,x) \neg \le \neg    \underset{\a \in \fA_t }{\esssup} \; \underset{\nu \in \cV_t }{\essinf}    \;
    Y^{t,x,  \a ( \nu ), \nu }_t \Big(\t_{\a,\nu},
     \wt{\f} \big(\t_{\a,\nu},  X^{t,x, \a ( \nu ), \nu }_{\t_{\a,\nu}} \big) \Big) ,
    $$
where    $\t_{\a,\nu} \dfnn  \inf \big\{s  \neg \in \neg  ( t , T  ] \neg  : \big( s,X^{t,x,\a(\nu),   \nu  }_s \big)
  \neg  \notin   \neg   O_{  \d} (t, x) \big\}   $.

 \end{thm}

 The significance of such a weak dynamic programming principle  lies in the following fact:
 Since $w_i$, $i=1,2$  may not be continuous in $t$,
 $ w_i \big(\t_{\beta,\mu},  X^{t,x,\mu, \beta (  \mu ) }_{\t_{\beta,\mu}} \big) $ may not be $\cF_{\t_{\beta,\mu}}-$measurable. Then $Y^{t,x,\mu, \beta (  \mu ) }_t \Big(\t_{\beta,\mu},
    w_i \big(\t_{\beta,\mu},  X^{t,x,\mu, \beta (  \mu ) }_{\t_{\beta,\mu}} \big) \Big)$ and thus
    the strong dynamic programming principle  may not be well-defined.

   \section{Viscosity Solutions of Related Fully Non-linear PDEs}

\label{sec:PDE}

\ms In this section,  we show that  the   priority  values
   are (discontinuous) viscosity solutions to the following partial differential equation  with a
fully non-linear Hamiltonian $H$:
  \bea \label{eq:PDE}
 \hspace{-0.3cm}     - \frac{\pa }{\pa t} w(t,x)
        -   H \big(t,x, w(t,x),  D_x w(t,x), D^2_x w(t,x)\big)
      \neg = \neg  0 ,    \q    \fa (t,x)
 \neg \in  \neg  (0,T)  \neg \times  \neg  \hR^k .    ~ \;       
\eea

 \if{0}
 \bea \label{eq:PDE}
   \min \Big\{ (w-\ul{l})(t,x), ~ \max \Big\{ - \neg \frac{\pa }{\pa t} w(t,x)
    - H \big(t,\xi, w(t,x), D_x w(t,x), D^2_x w(t,x)\big),
    (w-\ol{l})(t,x)  \Big\} \Big\}= 0 , ~ \fa (t,x) \in (0,T) \times \hR^k  .       
\eea
 \fi

\begin{deff}  \label{def:viscosity_solution}
 Let  us denote by $\hS_k$  the set of all  $\hR^{k \times k}-$valued symmetric matrices and let $H \dneg : [0,T] \neg \times   \neg  \hR^k  \neg  \times  \neg   \hR  \neg  \times  \neg    \hR^k
  \neg  \times    \hS_k  \neg \to \neg  [-\infty, \infty]$.
  An upper \(resp.\,lower\) semicontinuous function $w  \neg :  [0, T]  \neg \times \neg  \hR^k  \neg \to \neg  \hR$    is called  a viscosity subsolution \(resp. supersolution\) of \eqref{eq:PDE} 
  if   for any $(t_0,x_0, \vf )  \neg \in \neg  (0,T)  \neg \times \neg  \hR^k  \neg \times \neg  \hC^{1,2}\big([0,T]  \neg \times \neg  \hR^k\big)$
  such that 
  $w  \neg - \neg  \vf$ attains a strict local maximum $0$ \(resp.\;strict local minimum $0$\)  at $(t_0,x_0)$, we have
 \beas
        -  \frac{\pa }{\pa t} \vf (t_0,x_0)
     -  H \big(t_0,x_0, \vf (t_0,x_0),  D_x \vf (t_0,x_0), D^2_x \vf (t_0,x_0)\big)   \neg \le \neg   (\hb{resp.}   \ge )  \;  0.
 \eeas

  \end{deff}

  \ss    For any $(t,x,y,z,\G, u,v ) \neg \in \neg  [0,T]  \neg \times \neg  \hR^k \neg \times  \neg  \hR
   \neg \times  \neg  \hR^d  \neg \times \neg \hS_k \neg \times \neg  \hU  \neg \times \neg  \hV   $,
    set 
      \beas
    H  (t,x,y,z,\G,u,v) \dfnn \frac12 trace\big(\si \si^T(t,x,u,v) \,\G\big)+ z \cd b(t,x,u,v)
    +  f\big(t,x,y,z \cd \si(t,x,u,v), u,v \big) .
    \eeas
 We  consider   the following Hamiltonian functions:
  \beas
    \ul{H}_1(\Xi)   & \dneg \dfnn &  \dneg         \underset{u \in \hU  }{\sup} \;
 \linf{ \Xi'  \to \Xi }   \;    \underset{v \in \sO_u}{\inf}  \;  H (\Xi',u,v) , \q
 \ol{H}_1(\Xi)  \dfnn  \lmtd{n \to \infty}  \;     \underset{u \in \hU  }{\sup} \;
 \underset{v \in \sO^n_u  }{\inf}   \;  \lsup{  u' \to u} \; \underset{\Xi'     \in O_{\neg \frac{1}{n}} (\Xi)}{\sup} \;  H (\Xi',u',v) ,  \\
 \hb{and}  \q  \ol{H}_2(\Xi)  &  \dneg \dfnn &  \dneg     \underset{v \in \hV  }{\inf}   \;
 \lsup{ \Xi'  \to \Xi }    \;    \underset{u \in \sO_v}{\sup} \;  H (\Xi',u,v)  ,  \q
\ul{H}_2(\Xi)   \dfnn      \lmtu{n \to \infty}  \;     \underset{v \in \hV  }{\inf}  \;  \underset{u \in \sO^n_v  }{\sup}  \;
 \linf{  v' \to v} \; \underset{ \Xi'     \in O_{\neg \frac{1}{n}} (\Xi)}{\inf} \;  H (\Xi',u,v'),   \q
    \eeas
 where $\Xi = (t,x,y,z,\G)$,   $ \sO^n_u \dfnn \big\{ v \in \hV:  [v]_{\overset{}{\hV}}
  \le \k + n  [u]_{\overset{}{\hU}}  \big\}$,
    $ \sO^n_v \dfnn \big\{ u \in \hU:  [u]_{\overset{}{\hU}} \le \k + n   [v]_{\overset{}{\hV}}   \big\}$,
 $\sO_u \dfnn \underset{n \in \hN}{\cup} \sO^n_u = \b1_{\{u = u_0\}} \ol{O}_\k (v_0) + \b1_{\{u \ne u_0\}} \hV $
 and $\sO_v \dfnn \underset{n \in \hN}{\cup} \sO^n_v = \b1_{\{v = v_0\}} \ol{O}_\k (u_0) + \b1_{\{v \ne v_0\}} \hU $.

 \begin{rem}
 When $U$ and $V$ are compact \(say $\hU = \ol{O}_\k(u_0)$ and $\hV = \ol{O}_\k(v_0)$\),
 it holds for any $(u,v) \in \hU \times \hV$ and $n \in \hN$ that $ \big( \sO^n_u, \sO^n_v \big) = (\hV,\hU)$.
 If  further assuming as \cite{Buckdahn_Li_1} that for any $(x,y,z) \in \hR^k \times \hR \times \hR^d$,
 $b(\cd, x, \cd, \cd)$, $\si(\cd, x, \cd, \cd)$,$f(\cd, x, y,z,\cd, \cd)$ are all continuous in $(t,u,v)$,
 one can deduce from \eqref{b_si_linear_growth}$-$\eqref{f_Lip} that the continuity of $H (\Xi ,u,v)$   in
 $\Xi$ is uniform in $(u,v)$. It follows that
  \beas
 \ul{H}_1(\Xi)   =         \underset{u \in \hU  }{\sup} \;
 \linf{ \Xi'  \to \Xi }   \;    \underset{v \in \hV}{\inf}  \;  H (\Xi',u,v)
 =       \underset{u \in \hU  }{\sup} \;
     \underset{v \in \hV}{\inf}  \; \linf{ \Xi'  \to \Xi }   \;  H (\Xi',u,v)
     =  \underset{u \in \hU  }{\sup} \;
     \underset{v \in \hV}{\inf}  \;    H (\Xi,u,v) ,
 \eeas
 and that
 \beas
 \ol{H}_1(\Xi) & = & \lmtd{n \to \infty}  \;     \underset{u \in \hU  }{\sup} \;
 \underset{v \in \hV  }{\inf}   \;  \lsup{  u' \to u} \; \underset{\Xi'     \in O_{\neg \frac{1}{n}} (\Xi)}{\sup} \;  H (\Xi',u',v)  =      \underset{u \in \hU  }{\sup} \;
 \underset{v \in \hV  }{\inf}   \;  \lsup{  u' \to u} \; \lmtd{n \to \infty}  \;  \underset{\Xi'     \in O_{\neg \frac{1}{n}} (\Xi)}{\sup} \;  H (\Xi',u',v) \\
 & = & \underset{u \in \hU  }{\sup} \; \underset{v \in \hV  }{\inf}   \;  \lsup{  u' \to u} \;    H (\Xi,u',v)
 = \underset{u \in \hU  }{\sup} \; \underset{v \in \hV  }{\inf}   \;       H (\Xi,u,v)  = \ul{H}_1(\Xi)   .
 \eeas
  Similarly, $  \ul{H}_2(\Xi) =  \ol{H}_2(\Xi) = \underset{v \in \hV  }{\inf}   \;
      \underset{u \in \hU}{\sup} \;  H (\Xi,u,v)   $.

 \end{rem}

 \ss  For $i=1,2$,     Proposition \ref{prop_w_conti} implies that
 \beas
   ~ \;     \ul{w}_i (t,x )  \neg  \dfnn   \neg    \linf{ t' \to t}  w_i \,  (t',x )
    \neg  =  \neg    \linf{ (t',x') \to (t,x) }  w_i   (t',x')
   ~ \hb{and}  ~
    \ol{w}_i (t,x )   \neg   \dfnn   \neg    \lsup{ t' \to t}  w_i   (t',x )
   \neg  =  \neg  \lsup{ (t',x') \to (t,x) }  w_i  (t',x') , ~ \fa (t,x) \neg \in \neg  [0,T]  \neg \times \neg  \hR^k  .
 \eeas
  In fact,   $\ul{w}_i$  is the largest lower semicontinuous function below  $w_i$ (known as the lower    semicontinuous envelope  of $w_i$)
  while      $\ol{w}_i$ is the smallest upper semicontinuous function above $w_i$ (known as the upper    semicontinuous envelope  of $w_i$).

\begin{thm} \label{thm_viscosity}

    For $i=1,2$, $\ul{w}_i $ \(resp.\;$ \ol{w}_i $\)  is a  viscosity supersolution
 \(resp.\;subsolution\)  of   \eqref{eq:PDE} with the fully non-linear
 Hamiltonian $ \ul{H}_i $ \(resp.\;$ \ol{H}_i $\).
   \if{0}
     1\)     $\ul{w}_1 $ \(resp.\;$ \ol{w}_1 $\)  is a  viscosity supersolution
 \(resp.\;subsolution\)  of   \eqref{eq:PDE} with the fully non-linear
 Hamiltonian $ \ul{H}_1 $ \(resp.\;$ \ol{H}_1 $\).

  \no    2\)         $\ol{w}_2 $ \(resp.\;$ \ul{w}_2 $\)
   is a  viscosity subsolution \(resp.\;supersolution\) of  \eqref{eq:PDE} with the fully non-linear Hamiltonian $\ol{H}_2$    \(resp.\;$ \ul{H}_2 $\).
  \fi

\end{thm}

 Since there is no regularity, even semi-continuity,  in the fully non-linear
 Hamiltonian functions $ \ul{H}_i $ and  $ \ol{H}_i$, this  existence result
 of viscosity solutions to the fully non-linear  PDEs \eqref{eq:PDE} is quite general. In general, a comparison result for the PDEs that we analyze may not hold since  it is not clear whether $ \ul{H}_i $ =$ \ol{H}_i$ unless the control spaces are compact.

 \begin{rem}
  Given $i=1,2$ and  $x \in \hR^k$, although $w_i(T,x) =   g(x) $, it is possible that  neither
   $\ul{w}_i (T,x)$ nor $\ol{w}_i (T,x)$ equals to $g(x)$ since  
   $w_i$  may not be continuous in $t$. This phenomenon already appears in stochastic control problems with unbounded control; see e.g. \cite{BS13}.
 \end{rem}

 \section{Proofs}

  \label{sec:Proofs}



   \subsection{Proofs of  the results in Section \ref{sec:intro}}


    \ss \no {\bf Proof of Proposition \ref{BSDE_well_posed}:}
Set $\ff(s,\o,y,z ) \neg \dfnn \neg  \b1_{\{s \ge t\}} f (s,\o,y,z ) $, $\fa (s,\o,y,z )
 \neg \in \neg  [0,T] \times \O \times \hR \times \hR^d$.
  Clearly, $\ff$ is also a $\sP   \neg  \otimes \neg  \sB(\hR)   \neg \otimes \neg  \sB(\hR^d)/\sB(\hR)-$measurable
   function  Lipschitz continuous in $(y,z)$.
  As $    E   \Big[  \big(\int_0^T  \neg \big| \ff    (s, 0,0 )   \big| \, ds \big)^p \Big]
  = E   \Big[  \big(\int_t^T  \neg \big| f    (s, 0,0 )   \big| \, ds \big)^p \Big] < \infty $,
     Theorem 4.2 of \cite{BH_Lp_2003} shows that the BSDE
   \bea   \label{BSDE02}
    Y_s = \eta \neg + \neg  \int_s^T  \neg  \ff  (r,   Y_r, Z_r)  \, dr
   \neg-\neg \int_s^T \neg Z_r d B_r  , \q    s \in [0,T] .
    \eea
  admits a unique solution $ \big(Y ,Z   \big) \in \hG^p_\bF\big([0,T]\big) $.
 In particular, $\big\{ (Y_s ,Z_s  ) \big\}_{s \in [t,T]} \in \hG^p_\bF\big([t,T]\big) $ solves
 \eqref{BSDE01}.

 \ss  Suppose that  $(Y', Z') $ is another solution of \eqref{BSDE01} in $ \hG^p_\bF\big([t,T]\big) $.
     Let $(\wt{Y}', \wt{Z}') \in \hG^p_\bF([0,t])$ be the unique solution of the following BSDE with zero generator:
   \beas
   \wt{Y}'_s =  Y'_t - \int_s^t \wt{Z}'_r d B_r      , \q    s \in [0,t] .
   \eeas
   Actually, $\wt{Y}'_s = E[ Y'_t |\cF_s]$.  Then $(\cY', \cZ' )
   \dfnn  \big\{ \big( \b1_{\{s < t\}}\wt{Y}'_s \neg + \neg  \b1_{\{s  \ge  t\}} Y'_s,
   \b1_{\{s < t\}}\wt{Z}'_s  \neg + \neg  \b1_{\{s  \ge  t\}} Z'_s  \big) \big\}_{s \in [0,T]} \in \hG^p_\bF([0,T]) $
  also solves    BSDE \eqref{BSDE02}.
   So $(\cY', \cZ' ) = (Y,Z)$. In particular, $(Y'_s ,Z'_s ) = (Y_s, Z_s)$, $\fa s \in [t,T]$.

      \ss
        Given $A \in \cF_t$, multiplying $\b1_A$ to both sides of \eqref{BSDE01} yields that
      \beas
      \b1_A Y_s= \b1_A  \eta    \neg + \neg  \int_s^T  \neg
      \b1_A f (r,   \b1_A  Y_r,   \b1_A  Z_r )  \, dr
   \neg-\neg \int_s^T \neg  \b1_A  Z_r d B_r  , \q    s \in [t,T] .
     \eeas
   Let $(Y^A, Z^A) \in \hG^p_\bF([0,t])$ be the unique solution of the following BSDE with zero generator:
   \beas
   Y^A_s = \b1_A  Y_t - \int_s^t Z^A_r d B_r      , \q    s \in [0,t] .
   \eeas
    Then $ \big(\cY^A, \cZ^A \big)
   \dfnn  \big\{ \big( \b1_{\{s < t\}}Y^A_s \neg + \neg  \b1_{\{s  \ge  t\}}\b1_A Y_s,
   \b1_{\{s < t\}}Z^A_s  \neg + \neg  \b1_{\{s  \ge  t\}}\b1_A Z_s  \big) \big\}_{s \in [0,T]} \in \hG^p_\bF([0,T]) $
   solves the  BSDE
    \beas
    \cY^A_s= \b1_A \eta   \neg + \neg  \int_s^T  \neg
    f_A (r,   \cY^A_r,   \cZ^A_r )  \, dr
   \neg-\neg \int_s^T \neg \cZ^A_r d B_r  , \q    s \in [0,T] ,
    \eeas
   where $f_A (r,\o,y,z ) \dfnn \b1_{\{r  \ge  t\}} \b1_{\{ \o \in A \}} f (r,\o,y,z) $.
   Since
 $ \{\b1_{\{r  \ge  t   \} \cap A}   \}_{r \in [0,T]}$ is a right-continuous $\bF-$adapted process,
   the measurability and Lipschitz continuity of $ f $ imply that    $f_A$ is also
   a    $\sP   \otimes \sB(\hR)  \otimes \sB(\hR^d)/\sB(\hR)-$measurable function   Lipschitz continuous in $(y,z)$.
 Since    $     E \bigg[   \Big( \int_0^T \neg  \big| f_A (s,0, 0)  \big| ds  \Big)^p  \bigg]
  \le   E \bigg[   \Big( \int_t^T \neg  \big| f (s,0, 0)  \big| ds  \Big)^p  \bigg]
    < \infty $,
      applying       Proposition 3.2 of \cite{BH_Lp_2003}
  yields that
   \beas
      E \Big[ \b1_A \underset{s \in [t,T]}{\sup}  \big|   Y_s     \big|^p   \Big]
    &  \tneg  \le & \tneg   E \Big[ \, \underset{s \in [0,T]}{\sup}  \big|  \cY^A_s     \big|^p   \Big]
     \neg  \le  \neg  C(T,p,\g) E \bigg[ \b1_A  | \eta   |^p
     \neg  + \neg \Big( \int_0^T \neg  \big| f_A (s,0, 0)  \big| ds  \Big)^p  \bigg]   \\
 &  \tneg  = & \tneg  C(T,p,\g)  E \bigg[ \b1_A  |  \eta   |^p
  \neg + \neg   \b1_A    \Big( \int_t^T \neg  \big|  f (s,0, 0)  \big| ds  \Big)^p  \bigg]    .
     \eeas
         Letting $A$ vary in $\cF_t$ yields \eqref{eq:xvx131}.   \qed

       \ss \no {\bf Proof of Proposition \ref{prop_BSDE_estimate_comparison}:}
   (1)  Set $(\wh{Y}, \wh{Z}) \dfnn
      \big( Y^1  - Y^2, Z^1  - Z^2  \big) $, which  solves the  BSDE
     \bea   \label{eq:xvx017a}
      \wh{Y}_s= \eta_1 -  \eta_2   + \neg  \int_s^T  \neg
    \wh{f} (r,   \wh{Y}_r,   \wh{Z}_r )  \, dr
   \neg-\neg \int_s^T \neg \wh{Z}_r d B_r  , \q    s \in [t,T] ,
     \eea
     where $\wh{f}(r,\o,y,z) \neg \dfnn \neg  f_1 \big(r, \o,     y  +   Y^2_r    (\o),
     z  +   Z^2_r     (\o)   \big)      \neg - \neg f_2 \big(r, \o, Y^2_r    (\o),        Z^2_r    (\o) \big) $.
     Clearly, $ \wh{f} $ is a $\sP   \otimes \sB(\hR)  \otimes \sB(\hR^d)/\sB(\hR)-$ measurable function
       Lipschitz continuous in $(y,z)$.
     \if{0}
      Given $A \in \cF_t$, multiplying $\b1_A$ to both sides of \eqref{eq:xvx017a} yields that
      \beas
      \b1_A \wh{Y}_s= \b1_A \big( \eta_1 - \eta_2 \big)   \neg + \neg  \int_s^T  \neg
      \b1_A \wh{f} (r,   \b1_A  \wh{Y}_r,   \b1_A  \wh{Z}_r )  \, dr
   \neg-\neg \int_s^T \neg  \b1_A  \wh{Z}_r d B_r  , \q    s \in [t,T] ,
     \eeas

      \ss  Let $(\cY^A, \cZ^A) \in \hG^p_\bF([0,t])$ be the unique solution of the following BSDE with zero generator:
   \beas
   \cY^A_s = \b1_A \wh{Y}_t - \int_s^t \cZ^A_r d B_r      , \q    s \in [0,t] .
   \eeas
    Then $ \big(\wh{\cY}^A, \wh{\cZ}^A \big)
   \dfnn  \big\{ \big( \b1_{\{s < t\}}\cY^A_s \neg + \neg  \b1_{\{s  \ge  t\}}\b1_A \wh{Y}_s,
   \b1_{\{s < t\}}\cZ^A_s  \neg + \neg  \b1_{\{s  \ge  t\}}\b1_A \wh{Z}_s  \big) \big\}_{s \in [0,T]} \in \hG^p_\bF([0,T]) $
   solves the  BSDE
    \beas
    \wh{\cY}^A_s= \b1_A (  \eta_1 -  \eta_2  )   \neg + \neg  \int_s^T  \neg
    \wh{f}_A (r,   \wh{\cY}^A_r,   \wh{\cZ}^A_r )  \, dr
   \neg-\neg \int_s^T \neg \wh{\cZ}^A_r d B_r  , \q    s \in [0,T] ,
    \eeas
   where $\wh{f}_A (r,\o,y,z ) \dfnn \b1_{\{r  \ge  t\}} \b1_{\{ \o \in A \}} \wh{f} (r,\o,y,z) $.
   Since
 $ \{\b1_{\{r  \ge  t   \} \cap A}   \}_{r \in [0,T]}$ is a right-continuous $\bF-$adapted process,
   the measurability and Lipschitz continuity of $ f_1 $ imply that    $\wh{f}_A$ is also
   a    $\sP   \otimes \sB(\hR)  \otimes \sB(\hR^d)/\sB(\hR)-$measurable function   Lipschitz continuous in $(y,z)$.

  \ss    Since $ \hG^p_\bF\big([t,T]\big) \subset \hG^{\wt{p}}_{\bF } \big([t,T]\big)$
  by H\"older's inequality and since
  \beas
  \q   E \bigg[   \Big( \int_0^T \neg  \big| \wh{f}_A (s,0, 0)  \big| ds  \Big)^{\wt{p}}  \bigg]
  \le E \bigg[   \Big( \int_t^T \neg  \big| \wh{f} (s,0, 0)  \big| ds  \Big)^{\wt{p}}  \bigg]
  \le E \bigg[   \Big( \int_t^T \neg  \big| f_1 \big(s,    Y^2_s     ,      Z^2_s         \big)
      \neg - \neg f_2 \big(s,   Y^2_s ,  Z^2_s  \big)  \big| ds \Big)^{\wt{p}}    \bigg] < \infty ,
  \eeas
   applying       Proposition 3.2 of \cite{BH_Lp_2003} with $ \wt{p} $
  yields that
   \beas
      E \Big[ \b1_A \underset{s \in [t,T]}{\sup}  \big|   \wh{Y}_s     \big|^{\wt{p}}   \Big]
    &  \tneg  \le & \tneg   E \Big[ \, \underset{s \in [0,T]}{\sup}  \big|  \wh{\cY}^A_s     \big|^{\wt{p}}   \Big]
     \neg  \le  \neg  C(T,\wt{p},\g) E \bigg[ \b1_A  | \eta_1  \neg - \neg  \eta_2  |^{\wt{p}}
     \neg  + \neg \Big( \int_0^T \neg  \big| \wh{f}_A (s,0, 0)  \big| ds  \Big)^{\wt{p}}  \bigg]   \\
 &  \tneg  = & \tneg  C(T,\wt{p},\g)  E \bigg[ \b1_A  |  \eta_1  \neg - \neg  \eta_2  |^{\wt{p}}
  \neg + \neg   \b1_A    \Big( \int_t^T \neg  \big|  \wh{f} (s,0, 0)  \big| ds  \Big)^{\wt{p}}  \bigg]    .
     \eeas
         Letting $A$ vary in $\cF_t$ yields that
      \bea
        E \bigg[   \underset{s \in [t,T]}{\sup}   |   \wh{Y}_s      |^{\wt{p}}  \Big|\cF_t \bigg]
     & \le &  C(T,\wt{p},\g) E \bigg[    |  \eta_1 - \eta_2  |^{\wt{p}} +  \Big(  \int_t^T \neg  \big|  \wh{f} (s,0, 0)  \big|   ds \Big)^{\wt{p}} \Big|\cF_t  \bigg]  \label{eq:xvx023}  \\
     &=& C(T,\wt{p},\g) E \bigg[    |  \eta_1 - \eta_2  |^{\wt{p}}
     +  \Big(  \int_t^T \neg   \big| f_1 \big(s,    Y^2_s     ,      Z^2_s         \big)
      \neg - \neg f_2 \big(s,   Y^2_s ,  Z^2_s  \big)  \big| ds \Big)^{\wt{p}} \Big|\cF_t  \bigg] , \qq  \pas
      \nonumber
      \eea
      \fi
  Suppose that  $ E \Big[  \big( \int_t^T  \big| \wh{f}  (s,0,0 ) \big|  ds \big)^{\wt{p}} \Big]
      = E \Big[    \big( \neg \int_t^T   \neg  \big| f_1 (s, Y^2_s ,Z^2_s )
   -    f_2 (s, Y^2_s, Z^2_s) \big| ds    \big)^{\neg {\wt{p}}}  \Big] < \infty$ for some ${\wt{p}} \neg \in \neg (1,p]$.
   Since  $ \hG^p_\bF\big([t,T]\big) \subset \hG^{\wt{p}}_{\bF } \big([t,T]\big)$ by H\"older's inequality,
   applying Proposition \ref{BSDE_well_posed} with $p = \wt{p}$ shows that $(\wh{Y}, \wh{Z}) $ is the unique solution of
   BSDE$\big(t, \eta_1 \neg - \neg \eta_2, \wh{f}\big)$ in $\hG^{\wt{p}}_{\bF } \big([t,T]\big)$ satisfying
           \beas 
       E \bigg[ \, \underset{s \in [t,T]}{\sup}  |  \wh{Y}_s     |^{\wt{p}}  \Big| \cF_t  \bigg]
    \le  C(T,\wt{p},\g) E \bigg[  \,   |    \eta_1 \neg - \neg \eta_2   |^{\wt{p}}
   + \Big( \int_t^T  \big| \wh{f}  (s,0,0 ) \big|  ds \Big)^{\wt{p}} \bigg| \cF_t  \bigg] , \q   \pas ,
        \eeas
which is exactly \eqref{eq:n211}.

 \ss  \no (2)    Next, suppose that $ \eta_1 \neg \le  \neg  (\hb{resp.}\ge)  \eta_2  $, \pas ~
   and that   $ \d f_s  \neg \dfnn \neg  f_1 (s,  Y^2_s ,Z^2_s )  \neg - \neg  f_2 (s, Y^2_s ,Z^2_s )
  \neg \le \neg  (\hb{resp.}\ge) \, 0  $,  $ds  \neg \times \neg  dP-$a.s. on
  $[t,T] \times \O$.    By \eqref{f_Lip},
      \beas
      \mathfrak{a}_s \dfnn \b1_{\{\wh{Y}_s \ne 0 \}} \frac{ f_1 \big(s,
         Y^1_s     ,         Z^1_s         \big)
         - f_1 \big(s,  Y^2_s     ,   Z^1_s      \big) }{\wh{Y}_s}
         \in [-\g,\g ],     \q s \in [t,T]
     \eeas
     defines an $\bF-$progressively measurable bounded process.    For $i = 1,\cds, d$,
     analogous to process $\mathfrak{a}$
     \beas
        \mathfrak{b}^i_s & \tneg  \dfnn  & \tneg  \b1_{\{Z^{1,i}_s \ne Z^{2,i}_s \}} \frac{ 1 }{Z^{1,i}_s- Z^{2,i}_s}
         \Big( f_1 \big(s,
         Y^2_s     , (  Z^{2,1}_s , \cds  \neg , Z^{2,i-1}_s, Z^{1,i}_s, \cds  \neg , Z^{1,n}_s  )    \big) \\
        & & \hspace{3.5cm}  -    f_1 \big(s,         Y^2_s     ,    ( Z^{2,1}_s , \cds  \neg , Z^{2,i}_s, Z^{1,i+1}_s, \cds  \neg , Z^{1,n}_s   )   \big) \Big)  \neg  \in  \neg  [-\g,\g ], \q s  \neg \in \neg  [t,T]
       \eeas
       also defines an $\bF-$progressively measurable bounded process.

  \ss  Then  we can alternatively express \eqref{eq:xvx017a} as
      \beas
      \wh{Y}_s= \eta_1 - \eta_2  \neg + \neg  \int_s^T  \neg
  ( \mathfrak{a}_r      \wh{Y}_r + \mathfrak{b}_r \neg \cd \neg   \wh{Z}_r + \d f_r ) \, dr
   \neg-\neg \int_s^T \neg \wh{Z}_r d B_r  , \q    s \in [t,T] .
     \eeas

          Define $Q_s\dfnn \exp \big\{ \int_t^s \mathfrak{a}_r d r-\frac{1}{2}\int_t^s
|\mathfrak{b}_r|^2 d r +\int_t^s \mathfrak{b}_r d B_r \big\}$, $s \in [t,T]$. Applying integration by parts  yields that
 \bea
 Q_s \wh{Y}_s &\tneg  \dneg  =  & \tneg  \dneg  Q_T \wh{Y}_T \neg + \dneg
 \int_s^T Q_r \big( \mathfrak{a}_r \wh{Y}_r
 \neg + \neg    \mathfrak{b}_r \neg  \cd \neg \wh{Z}_r \neg + \neg  \d f_r  \big) d r
    \neg  - \dneg \int_s^T \neg   Q_r \wh{Z}_r d B_r
  \neg - \dneg \int_s^T\neg \wh{Y}_r Q_r \mathfrak{a}_r d r  \neg -\dneg \int_s^T\neg \wh{Y}_r Q_r \mathfrak{b}_r d B_r
   \neg -\dneg \int_s^T\neg   Q_r \mathfrak{b}_r \neg \cd \neg \wh{Z}_r  d r  \nonumber \\
& \tneg  \dneg  = & \tneg  \dneg  Q_T (\eta_1 - \eta_2) \neg + \dneg
 \int_s^T Q_r      \d f_r    d r \neg - \dneg    \int_s^T
  Q_r ( \wh{Z}_r  \neg + \neg  \wh{Y}_r   \mathfrak{b}_r )  d B_r , \q    \pas    \label{eq:xvx077}
 \eea
  One can deduce from the Burkholder-Davis-Gundy inequality and H\"older's inequality that
 \bea
 E    \bigg[ \underset{s \in [t,T]}{\sup} \Big| \int_t^s \neg
  Q_r ( \wh{Z}_r  \neg + \neg  \wh{Y}_r   \mathfrak{b}_r  )  d B_r \Big| \bigg]
  & \tneg  \tneg \le & \tneg \tneg   c_0 E     \bigg[   \Big( \int_t^T  \neg
  Q^2_r |   \wh{Y}_r   \mathfrak{b}_r  \neg + \neg  \wh{Z}_r |^2  d r \Big)^{\frac12} \bigg]
    \neg \le  \neg  c_0 E     \bigg[  \underset{s \in [t,T]}{\sup} | Q_r |
    \bigg\{ \underset{s \in [t,T]}{\sup} | \wh{Y}_r |
   \neg  +  \neg  \Big( \int_t^T  \neg | \wh{Z}_r  |^2  d r \Big)^{\frac12}    \bigg\} \bigg]  \nonumber  \\
 & \tneg  \tneg  \le & \tneg  \tneg   c_0 \bigg( E     \Big[  \underset{s \in [t,T]}{\sup} | Q_r |^{\wh{p}} \Big]
 \bigg)^{\neg 1 /  \wh{p}} \neg \Big( \big\| \wh{Y}_r\big\|_{\hC^p_\bF( [t,T] )    } + \big\| \wh{Z}_r\big\|_{\hH^{2,p} _\bF([t,T], \hR^d)} \Big) ,   \label{eq:xvx075}
    \eea
 where $\wh{p} = \frac{p}{p-1} $.   Also,  Doob's martingale inequality implies that
    \beas
    E     \Big[  \underset{s \in [t,T]}{\sup} | Q_r |^{\wh{p}} \Big]
  & \tneg \dneg \le & \tneg \dneg   c_0 E     \Big[     | Q_T |^{\wh{p}} \Big]
     \neg  =  \neg  c_0 E     \bigg[     \exp \bigg\{ \wh{p} \neg  \int_t^T   \neg  \mathfrak{a}_r d r
      \neg +   \neg      \frac{\wh{p}^{\,2} \neg - \neg 1}{2 }  \neg   \int_t^T   \neg
|\mathfrak{b}_r|^2 d r   \neg -   \neg \frac{\wh{p}^{\,2}}{2 } \neg \int_t^T  \neg
|\mathfrak{b}_r|^2 d r   \neg +  \neg  \wh{p} \neg \int_t^T   \neg  \mathfrak{b}_r d B_r \bigg\} \bigg] \\
 & \tneg  \dneg \le &  \dneg \tneg c_0 \exp\Big\{\wh{p} \g T \neg + \neg \frac{\wh{p}^{\,2} \neg - \neg 1}{2}\g^2 T \Big\}
 E     \bigg[     \exp \bigg\{     \neg -   \neg \frac{\wh{p}^{\,2}}{2 } \neg \int_t^T  \neg
|\mathfrak{b}_r|^2 d r   \neg +  \neg  \wh{p} \neg \int_t^T   \neg  \mathfrak{b}_r d B_r \bigg\} \bigg]
  \neg = \neg  c_0 \exp\Big\{\wh{p} \g T \neg + \neg \frac{\wh{p}^{\,2} \neg - \neg 1}{2}\g^2 T \Big\} ,
    \eeas
    which together with \eqref{eq:xvx075} shows that
  $\big\{ \int_t^s \neg
  Q_r (  \wh{Y}_r   \mathfrak{b}_r  \neg + \neg  \wh{Z}_r)  d B_r \big\}_{s \in [t,T]} $  is a uniformly integrable martingale. Then for any $s \in [t,T]$, taking $E[\cd |\cF_s]$ in \eqref{eq:xvx077} yields that \pas
  \beas
   Q_s \wh{Y}_s = E \bigg[  Q_T (\eta_1 - \eta_2)  + \int_s^T Q_r \d f_r  d r \bigg| \cF_s \bigg] \le (\hb{resp.}\,\ge)\, 0 ,
       \q \hb{thus} \q  \wh{Y}_s \le (\hb{resp.}\,\ge)\,  0 .
   \eeas
  By the continuity of process $\wh{Y}$, it holds $P-$a.s. that   $ Y^1_s \le (\hb{resp.}\,\ge)\, Y^2_s$  for any $s \in [t,T]$.  \qed

\subsection{Proofs of the Results in Section \ref{sec:zs_drgame}}

  \label{subsection:Proofs_S2}

 \ss \no {\bf Proof of Example \ref{eg_example2}:}
  For any $(t,u) \in [0,T] \times \hU$, the continuity of $\vf$ and
 \eqref{eq:xqxqx011} show that   $   \{v  \in  [-   \k|u|,\k|u| ] :  \vf (t,u, v)=0 \}$
 is a non-empty closed set. So we can define $\sV(t,u) \dfnn \min  \{v  \in  [-   \k|u|,\k|u| ] :  \vf (t,u, v)=0 \} $.

 Given $n \in \hN$,  for any $i =0,\cds \neg , 2^n - 1  $ and $j \in \hZ$,
 we set $t^n_i = i 2^{-n}T$, $u^n_j = j 2^{-n} $
 and $\psi^n_{i,j} \dfnn \underset{(t,u) \in \cD^n_{i,j}}{\inf} \sV (t,u) \in [-\k - \k |u|, \k + \k |u|]  $ with
  \beas
  \cD^n_{i,j} \dfnn \left\{
  \ba{ll}
   \, \tneg   [t^n_i, t^n_{i+1}) \times [u^n_j, u^n_{j+1})  , \q & \hb{if } i < 2^n-1, \\
   \, \tneg   [t^n_i, T] \times  [u^n_j, u^n_{j+1}) , \q & \hb{if } i = 2^n-1 .
  \ea
  \right.
  \eeas
 Clearly,
   $\dis \psi_n (t,u) \neg \dfnn \neg  \sum^{2^n-1}_{i =0 } \sum_{j \in \hZ} \psi^n_{i,j} \b1_{\{(t,u) \in \cD^n_{i,j}\}}
    \neg \in \neg  [-\k  \neg - \neg  \k |u|, \k  \neg + \neg  \k |u|]  $,
  $\fa (t,u)  \neg \in \neg  [0,T]  \neg \times \neg  \hU$ defines a $\sB([0,T]) \otimes \sB(\hU )  /\sB(\hR )$ $-$measurable function.
  As $ \psi_n   \neg  \le \neg  \psi_{n+1} $, 
  the function $\psi (t,u) \neg  \dfnn \neg  \lmtu{n \to \infty}  \psi_n (t,u)  \neg \in  \neg
  [-\k  \neg - \neg  \k |u|, \k  \neg + \neg  \k |u|]  $,
  $\fa (t,u)  \neg \in \neg  [0,T]  \neg \times \neg  \hU$ is also
   $\sB([0,T])  \neg \otimes \neg  \sB(\hU )  /\sB(\hR )-$measurable.

 \ss  Now, let $(t,u) \in [0,T] \times \hU$ and $\e>0$. By the continuity of $\vf$ in
 $t$, there exists a $\d \in (0,\e/3\g)$ such that
  \bea  \label{eq:xqxqx014}
     \big|\vf\big(s,u, \psi(t,u)\big)- \vf\big(t,u, \psi(t,u)\big)\big|   < \e/3  ,  \q
   \fa s \in [t-\d,t+\d]\cap [0,T]  .
   \eea
  For any $n > \log_2  ( 1\vee T ) - \log_2  ( \d )$,
    $(t,u) \in \cD^n_{i,j} $ for some  $(i,j) \in \{0,\cds \neg ,2^n-1 \} \times \hZ$, and we can find
    $(t',u') \in \cD^n_{i,j}$ such that $  \sV(t',u') \le \psi^n_{i,j} + \d $. Then
    \eqref{eq:xqxqx014} and the Lipschitz continuity of $\vf$ in $(u,v)$ show that
  \beas
  && \hspace{-0.8cm} \big| \vf(t,u, \psi(t,u) )  \big| =
  \big| \vf(t,u, \psi(t,u) ) - \vf(t',u',\sV (t',u') ) \big| \\
  && \le  \big| \vf(t,u, \psi(t,u) ) \neg - \neg  \vf(t' ,u ,\psi (t,u) ) \big|  \neg + \neg
  \big| \vf(t',u, \psi(t,u) )  \neg - \neg  \vf(t' ,u ,\psi_n(t,u) ) \big|  \neg + \neg
  \big| \vf(t',u, \psi^n_{i,j} )  \neg - \neg  \vf(t',u',\sV (t',u') ) \big| \\
  && \le \e /3 + \g \big|\psi(t,u) - \psi_n(t,u) \big|
  + \g \big(|u-u'|+ |\psi^n_{i,j}- \sV (t',u')| \big)  \\
  && \le \e + \g \big|\psi(t,u) - \psi_n(t,u) \big|  + 2 \g \d
  \le  \e   + \g \big|\psi(t,u) - \psi_n(t,u) \big| .
  \eeas
  Letting $n \to \infty$ yields that $  \big| \vf(t,u, \psi(t,u) )  \big| \le  \e$. Then as
  $\e \to 0$, we obtain that $ \vf(t,u, \psi(t,u) ) = 0 $.

 \ss  Similarly, we can construct a measurable function $\wt{\psi}$ on $ [0,T] \times \hV $ such that
  \beas
  \vf \big(t, \wt{\psi}(t,v), v \big) = 0
 \q \hb{and} \q  |\wt{\psi}(t,v)| \le \k (1+|v|) , \q \fa (t,v) \in [0,T] \times \hV.
  \eeas
 Hence   \(A-u\) and \(A-v\) are satisfied.   \qed

  \no {\bf Proof of Lemma \ref{lem_control_combine}:} It suffices to prove for $\cU_t-$controls.
  Let $s \in [t,T]$ and $U \in \sB \big(\hU  \big)  $.
     Since   $ \[t, \t  \[ , \[\t , T \] \neg \in \neg  \sP $,
      we see that both  $\cD_1  \neg \dfnn   \neg    \[t, \t \[
      \, \cap \, ( [t,s] \times \O)   $
      and $ \cD_2  \neg \dfnn   \neg   \[\t, T \]
      \, \cap \, ( [t,s] \times \O) $ belong to $  \sB\big([t,s]\big)   \otimes    \cF_s $.
              It then follows   that
          \beas
      && \hspace{-1cm}   \big\{ (r,  \o) \neg \in \neg  [t,s] \times \O :   \mu_r (\o) \neg \in \neg U \big\}
        =      \big\{ (r,  \o) \neg \in \neg  \cD_1 :   \mu^1_r (\o) \neg \in \neg U \big\}
        \cup  \big\{ (r,  \o) \neg \in \neg   \cD_2 :   \mu^2_r (\o) \neg \in \neg U \big\}  \\
         &  &      =
      \Big( \cD_1 \cap \big\{ (r,  \o) \neg \in \neg  [t,s] \times \O :   \mu^1_r (\o) \neg \in \neg U \big\} \Big)
      \cup \Big( \cD_2 \cap \big\{ (r,  \o) \neg \in \neg  [t,s] \times \O  :   \mu^2_r (\o) \neg \in \neg U \big\} \Big)
      \neg   \in \neg  \sB\big([t,s]\big) \otimes \cF_s  ,
       \eeas
   which shows that  the process $\mu$ is   $\bF-$progressively  measurable.

  For $i=1,2$, suppose that
   $E \neg \int_t^T  [     \mu^i_s    ]^{q_i}_{\overset{}{\hU}} ds   < \infty$ for some $q_i > 2$.
   One can deduce  that
   $
     E \neg \int_t^T \neg  [  \mu_r   ]^{q_1 \land q_2}_{\overset{}{\hU}} dr
     \le E \neg \int_t^T \neg  [     \mu^1_r      ]^{q_1 \land q_2}_{\overset{}{\hU}} dr
       + E \neg \int_t^T \neg  [      \mu^2_r     ]^{q_1 \land q_2}_{\overset{}{\hU}} dr < \infty  $.
      Thus $\mu \in \cU_t$.     \qed

\ss \no {\bf Proof of Lemma \ref{lemm_control_compare}:}
  Both  $ \big\{   X^{t,\xi,\mu,\nu}_{\t \land  s}\big\}_{s \in [t,T]}$
  and $  \big\{   X^{t,\xi,\wt{\mu}, \wt{\nu}}_{\t \land  s}\big\}_{s \in [t,T]}$
     satisfy the same SDE:
    \bea \label{eq:p605}
  X_s=    \xi  + \int_t^s     b^{\mu,\nu}_{\t}  (r, X_r   ) \, dr
  + \int_t^s      \si^{\mu,\nu}_{\t} (r, X_r ) \, dB_r,  \q s \in [t, T] ,
    \eea
 where  $  b^{\mu,\nu}_{ \t} (r,\o, x ) \dfnn   \b1_{\{r < \t (\o) \}}  b^{\mu,\nu}  (r, \o,x )
 $ and $  \si^{\mu,\nu}_{ \t} (r,\o, x ) \dfnn   \b1_{\{r < \t (\o) \}}   \si^{\mu,\nu}  (r, \o,x )
 $,   $ \fa (r,\o,x) \in [t,T] \times \O \times \hR^k$. Like $b^{\mu,\nu}$ and $\si^{\mu,\nu}$,
  $b^{\mu,\nu}_{ \t}$ is a $\sP \otimes \sB(\hR^k)  /\sB(\hR^k)-$measurable function
  and  $\si^{\mu,\nu}_{ \t}$ is a $\sP \otimes \sB(\hR^k)  /\sB(\hR^{k \times d})-$measurable function that
  is Lipschitz continuous in $(y,z)$ with coefficient $\g$ and satisfies
  \beas
      E   \left[ \Big(\int_t^T  \neg \big| b^{\mu,\nu}_\t    (s, 0 )   \big| \, ds \Big)^2
   + \Big(\int_t^T  \neg \big| \si^{\mu,\nu}_\t    (s, 0 )   \big| \, ds \Big)^2 \right]
     < \infty .
     \eeas
   Thus          \eqref{eq:p605} has a unique solution. It then holds \pas ~ that
     \bea  \label{eq:xxd137}
         X^{t,\xi,\mu,\nu}_{\t \land  s}   =    X^{t,\xi,\wt{\mu}, \wt{\nu}}_{\t \land  s},    \q
     \fa  s    \in    [t,T ]   .
     \eea

    One can   deduce that
 \beas
 \q  X^{t,\xi,\mu,\nu}_{  s} - X^{t,\xi,\mu,\nu}_{\t \land s  }
\neg = \neg   X^{t,\xi,\mu,\nu}_{\t \vee s} - X^{t,\xi,\mu,\nu}_{\t  }
 \neg  = \neg  \int_\t^{\t \vee s}  \neg  b \big(r, X^{t,\xi,\mu,\nu}_r,\mu_r,\nu_r  \big) dr  \neg + \neg
   \int_\t^{\t \vee s}  \neg  \si \big(r, X^{t,\xi,\mu,\nu}_r,\mu_r,\nu_r  \big) d B_r ,   \q s \in [t, T].
 \eeas
 Multiplying $\b1_A$ on both sides yields that
  \beas
     \cX_s & \tneg \dfnn & \tneg  \b1_A \big( X^{t,\xi,\mu,\nu}_{  s} -    X^{t,\xi,\mu,\nu}_{\t \land s  } \big)
   =   \int_\t^{\t \vee s} \b1_A b \big(r, X^{t,\xi,\mu,\nu}_r,\mu_r,\nu_r  \big) dr +
     \int_\t^{\t \vee s} \b1_A \si \big(r, X^{t,\xi,\mu,\nu}_r,\mu_r,\nu_r  \big) d B_r  \\
&  \tneg  = & \tneg    \int_t^s \b1_{\{r  \ge  \t\}} \b1_A b \big(r, \cX_r +   X^{t,\xi,\mu,\nu}_{\t \land r }, \mu_r,\nu_r  \big) dr +
     \int_t^s \b1_{\{r  \ge  \t\}} \b1_A \si \big(r,\cX_r +   X^{t,\xi,\mu,\nu}_{\t \land r }, \mu_r,\nu_r  \big) d B_r
    ,   \q s \in [t, T].
 \eeas
 Similarly, we see from \eqref{eq:xxd137} that
 \beas
 \wt{\cX}_s & \tneg \dneg \dfnn & \tneg  \dneg    \b1_A \big( X^{t,\xi,\wt{\mu}, \wt{\nu}}_s  \neg  -  \neg    X^{t,\xi,\wt{\mu}, \wt{\nu}}_{\t \land s} \big)
  \neg  =   \dneg  \int_t^s  \neg  \b1_{\{r  \ge  \t\}} \b1_A b \big(r, \wt{\cX}_r \neg + \neg   X^{t,\xi,\wt{\mu}, \wt{\nu}}_{\t \land r },\wt{\mu}_r,\wt{\nu}_r  \big) dr  \neg  +  \dneg     \int_t^s  \neg  \b1_{\{r  \ge  \t\}} \b1_A \si \big(r, \wt{\cX}_r \neg + \neg   X^{t,\xi,\wt{\mu}, \wt{\nu}}_{\t \land r },\wt{\mu}_r,\wt{\nu}_r  \big) d B_r \\
  &  \tneg  \dneg  = &  \tneg  \dneg    \int_t^s \b1_{\{r  \ge  \t\}} \b1_A b \big(r,  \wt{\cX}_r \neg + \neg   X^{t,\xi,\mu,\nu}_{\t \land r  },\mu_r,\nu_r  \big) dr +     \int_t^s \b1_{\{r  \ge  \t\}} \b1_A \si \big(r, \wt{\cX}_r \neg + \neg   X^{t,\xi,\mu,\nu}_{\t \land r },\mu_r,\nu_r  \big) d B_r     ,   \q s \in [t, T] .
 \eeas
  To wit,   $\cX, \wt{\cX} \in \hC^2_\bF ([t,T], \hR^k)$ satisfy the same SDE:
 \bea  \label{eq:xxd143}
  X_s=      \int_t^s     \wh{b}   (r, X_r   ) \, dr
  + \int_t^s      \wh{\si}  (r, X_r ) \, dB_r,  \q s \in [t, T] ,
    \eea
where  $  \wh{b}  (r,\o, x ) \dfnn   \b1_{\{r  \ge  \t (\o) \}} \b1_{\{\o \in A \}}  b   \big(r,  x+ X^{t,\xi,\mu,\nu}_{\t \land r }(\o), \mu_r (\o),\nu_r (\o) \big)
 $ and $  \wh{\si}  (r,\o, x ) \dfnn   \b1_{\{r  \ge  \t (\o) \}} \b1_{\{\o \in A \}}  \si^{\mu,\nu}  \big(r,  x+ X^{t,\xi,\mu,\nu}_{\t \land r }(\o), \mu_r (\o),\nu_r (\o) \big)
 $,   $ \fa (r,\o,x) \in [t,T] \times \O \times \hR^k$.  The measurability of functions $b$,
   $X^{t,\xi,\mu,\nu}  $,  $\mu$ and $\nu$ implies that the mapping $(r,\o,x) \to b   \big(r, \o,x+ X^{t,\xi,\mu,\nu}_{\t \land r }(\o), \mu_r (\o),\nu_r (\o) \big) $  is $\sP \otimes \sB(\hR^k)  /\sB(\hR^k)-$measurable. Clearly,
 $ \{\b1_{\{r  \ge  \t   \} \cap A}   \}_{r \in [t,T]}$ is a right-continuous $\bF-$adapted process.
 Thus $\wh{b}$ is also $\sP \otimes \sB(\hR^k)  /\sB(\hR^k)-$measurable. Similarly, $\wh{\si}$ is
 $\sP \otimes \sB(\hR^k)  /\sB(\hR^{k \times d})-$measurable. By \eqref{b_si_Lip}, both $\wh{b}$ and $\wh{\si}$
 are Lipschitz continuous in $x$. Since
\beas
    E   \left[ \Big(\int_t^T  \neg \big| \wh{b}    (r, 0 )   \big| \, dr \Big)^2
   + \Big(\int_t^T  \neg \big| \wh{\si}     (r, 0 )   \big| \, dr \Big)^2 \right]
   \le c_0 + c_0 E \neg \left[     \big| X^{t,\xi,\mu,\nu}_\t  \big|^2 \right] + c_0 E \neg \int_t^T  \neg \big([ \mu_r ]^2_{\overset{}{\hU}} \neg+\neg [ \nu_r ]^2_{\overset{}{\hV}} \big) dr
   < \infty
   \eeas
 by \eqref{b_si_linear_growth}, \eqref{b_si_Lip} and H\"older's inequality,
 the SDE \eqref{eq:xxd143} admits      a  unique  solution. Hence,   $P \big(\cX_s = \wt{\cX}_s,\; \fa s \in [t,T]\big)=1$,  which together with \eqref{eq:xxd137} proves \eqref{eq:p611}.  \qed

   \ss \no {\bf Proof of Lemma \ref{lem_estimate_Y}:}
     For $i=1,2$,  let  $\Th_i  \neg  \dfnn  \neg  (t,\xi_i, \mu, \nu ) $ and set
     $ (Y^i,Z^i) \dfnn \Big(  Y^{\Th_i }  \big(T, g\big(X^{\Th_i }_T \big) \big) ,
    Z^{\Th_i }  \big(T, g\big(X^{\Th_i }_T  \big) \big) \Big) $.
    Given $\wt{p} \in (1,p]$, \eqref{f_Lip} and H\"older's inequality show that
      \beas
       E \bigg[    \bigg(  \int_t^T \neg  \big|      f^{\Th_1}_T \big(r,     Y^2_r    ,    Z^2_r      \big)
        \neg  -  \neg   f^{\Th_2}_T \big(r,   Y^2_r     ,   Z^2_r  \big)  \big|  ds  \bigg)^{\wt{p}}      \bigg]
         \neg  \le \neg  c_{\wt{p}} E \bigg[   \underset{s \in [t,T]}{\sup}  \big| X^{\Th_1}_s  \neg - \neg  X^{\Th_2 }_s  \big|^{\frac{2 \wt{p}}{p}}      \bigg]
          \neg \le \neg  c_{\wt{p}} \bigg\{ E \bigg[   \underset{s \in [t,T]}{\sup}  \big| X^{\Th_1}_s  \neg - \neg  X^{\Th_2 }_s  \big|^2      \bigg] \bigg\}^{\frac{ \wt{p}}{p}} \neg  < \neg  \infty .
                 \eeas
    Then   we can deduce from \eqref{eq:n211}   that
      \beas
      E \Big[   \underset{s \in [t,T]}{\sup}  \big|   Y^1_s  \neg - \neg Y^2_s  \big|^{\wt{p}}  \big|\cF_t \Big]
    & \tneg \le & \tneg  c_{\wt{p}} E \bigg[   \big| g\big(X^{\Th_1}_T \big)  \neg - \neg  g\big(X^{ \Th_2 }_T \big) \big|^{\wt{p}}
      \neg + \neg   \int_t^T \neg  \big|      f^{\Th_1}_T \big(r,     Y^2_r    ,    Z^2_r      \big)
        \neg  -  \neg   f^{\Th_2}_T \big(r,   Y^2_r     ,   Z^2_r  \big)  \big|^{\wt{p}} ds   \bigg| \cF_t    \bigg]  \\
     &  \tneg \le & \tneg  c_{\wt{p}} E \bigg[   \underset{s \in [t,T]}{\sup}  \big| X^{\Th_1}_s  \neg - \neg  X^{\Th_2 }_s  \big|^{\frac{2 \wt{p}}{p}}  \Big| \cF_t    \bigg]   , \q   \pas
      \eeas
     Then         a standard a priori estimate  of SDEs   (see e.g.  \cite[pg. 166-168]{Ikeda_Watanabe_1989}
     and \cite[pg. 289-290]{Kara_Shr_BMSC}) leads to that
     \beas
   \hspace{2.3cm}    E \bigg[   \underset{s \in [t,T]}{\sup}  \big|  Y^1_s  \neg - \neg Y^2_s  \big|^{\wt{p}}  \Big|\cF_t \bigg]
     \le c_{\wt{p}} E \bigg[   \underset{s \in [t,T]}{\sup}  \big| X^{\Th_1}_s \neg - \neg X^{\Th_2 }_s  \big|^{\frac{2 \wt{p}}{p}}
      \Big| \cF_t    \bigg]
     \le c_{\wt{p}}   | \xi_1 \neg - \neg \xi_2 |^{\frac{2 \wt{p}}{p}} , \q  \pas   \hspace{2.3cm}  \hb{\qed}
       \eeas

  \ss \no {\bf Proof of Proposition \ref{prop_value_bounds}:}  Given $\beta \neg \in \neg \fB_t$,
 \eqref{eq:xvx131}  and H\"older's inequality  imply   that
 \bea \label{eq:xvx051}
  \big| J \big(t,x,u_0, \beta (u_0)\big) \big|^p
   & \le & E \bigg[ \, \underset{s \in [t,T]}{\sup} \Big|  Y^{t,x,u_0, \beta (u_0)}_s  \big(T, g\big(X^{t,x,u_0, \beta (u_0)}_T \big) \big) \Big|^p  \Big| \cF_t  \bigg]  \nonumber   \\
   & \le & c_0 E \bigg[   \big|    g\big(X^{t,x,u_0, \beta (u_0)}_T \big)   \big|^p
   + \int_t^T  \big| f^{t,x,u_0, \beta (u_0)}_T (s,0,0 ) \big|^p ds \Big| \cF_t  \bigg] , \q \pas
   \eea
 Since $ \big[ \big( \beta (u_0)\big)_s   \big]_{\overset{}{\hV}}   \neg \le \neg \k $, $ds   \times   dP-$a.s.,
     the  $2/p-$H\"older continuity of $g$,  \eqref{f_linear_growth}, \eqref{f_Lip}
    as well as a conditional-expectation version of \eqref{eq:esti_X_1a} show that  \pas
  \bea
 && \hspace{-1cm}  \big| J \big(t,x,u_0, \beta (u_0)\big) \big|^p
   \le c_0 + c_0 E \bigg[   \big|     X^{t,x,u_0, \beta (u_0)}_T     \big|^2
   + \int_t^T \big( \big| X^{t,x,u_0, \beta (u_0)}_s \big|^2 + \big[ \big(\beta (u_0)\big)_s \big]^2_{\overset{}{\hV}} \big) ds \Big| \cF_t  \bigg]  \nonumber \\
  && \le c_\k + c_0 E \bigg[  \underset{s \in [t,T]}{\sup} \big|     X^{t,x,u_0, \beta (u_0)}_s     \big|^2 \Big| \cF_t  \bigg]
   \le c_\k + c_0   |x|^2  + c_0 E \bigg[     \int_t^T   \big[ \big(\beta (u_0)\big)_s \big]^2_{\overset{}{\hV}}  ds \Big| \cF_t  \bigg] \le c_\k + c_0   |x|^2 .    \label{eq:xvx053}
  \eea
 So it follows that
    \beas
  w_1 (t,x)  \ge  \underset{\beta \in \fB_t }{\essinf} \;
    J \big(t,x,u_0, \beta (u_0)\big)  \ge   - c_\k - c_0   |x|^{2/p}  , \q \pas
    \eeas

  We extensively set $\psi(t,u ) \neg \dfnn \neg  v_0$, $\fa (t,u)  \neg \in \neg  [0,T]  \neg \times \neg  O_\k(u_0)$,  then it is  $ \sB([0,T])
  \neg \times \neg \sB(\hU) / \sB(\hV) -$measurable.
    For any $\mu \in \cU_t$,    the measurability   of function $\psi$ and process $\mu$
    implies  that
     \bea   \label{def_beta_psi}
       \big(\beta_\psi (\mu) \big)_s \neg \dfnn \neg  \psi (s, \mu_s)  ,  \q  s  \neg \in \neg  [t,T]
     \eea
    defines  a $\hV-$valued, $\bF-$progressively measurable process,  and we see from (A-u) that
       $ [\big(\beta_\psi (\mu) \big)_s]_{\overset{}{\hV}}  \neg \le \neg  \k  \neg + \neg \k  [\mu_s]_{\overset{}{\hU}} $,
       $\fa s  \neg \in \neg  [t,T]$.    So  $\beta_\psi (\mu)  \neg \in \neg  \cV_t$.
        Let    $\mu^1, \mu^2 \neg  \in \neg  \cU_t $ such that  $\mu^1  \neg = \neg   \mu^2 $,
$ds  \neg \times \neg  d P -$a.s. on $ \[t,\t\[ \, \cup \, \[\t, T\]_A   $  for some  $ \t \neg \in  \neg   \cS_{t,T}$
  and $ A  \neg  \in  \neg  \cF_\t  $.   It  clearly holds  $ds  \neg \times \neg  d P -$a.s. on $ \[t,\t\[ \, \cup \, \[\t, T\]_A   $ that
   \beas
  \big(\beta_\psi (\mu^1) \big)_s = \psi (s,\mu^1_s) = \psi (s,\mu^2_s) =  \big(\beta_\psi (\mu^2) \big)_s \, .
   \eeas
   Hence, $\beta_\psi \in \fB_t $.

   Fix a $u_\sharp \in \pa O_\k(u_0)$. For any  $\mu \in  \cU_t$, similar to \eqref{eq:xvx051} and \eqref{eq:xvx053},
   we can deduce   that \pas
 \bea
  \big| J \big(t,x,\mu, \beta_\psi (\mu)\big) \big|^p
   & \tneg  \le  & \tneg  c_0 E \bigg[   \big|    g\big(X^{t,x,\mu, \beta_\psi (\mu)}_T \big)   \big|^p
   + \int_t^T  \big| f^{t,x,\mu, \beta_\psi (\mu)}_T (s,0,0 ) \big|^p ds \Big| \cF_t  \bigg] \nonumber \\
    & \tneg  \le  & \tneg  c_0 + c_0 E \bigg[   \big|     X^{t,x,\mu, \beta_\psi (\mu)}_T     \big|^2
       + \int_t^T \Big( \b1_{\{\mu_s \in  O_\k(u_0) \}} \big| f (s, X^{t,x,\mu, \beta_\psi (\mu)}_s, 0,0, \mu_s, v_0  ) \big|^p \nonumber   \\
  & \tneg  &  +  \b1_{\{\mu_s \notin   O_\k(u_0) \}} \big| f \big(s, X^{t,x,\mu, \beta_\psi (\mu)}_s, 0,0, u_\sharp, \psi (s,u_\sharp) \big) \big|^p \Big) ds \Big| \cF_t  \bigg] \label{eq:xvx055} \\
   & \tneg   \le & \tneg c_\k + c_0 E \bigg[  \underset{s \in [t,T]}{\sup} \big|     X^{t,x,\mu, \beta_\psi (\mu)}_s     \big|^2
     \Big| \cF_t  \bigg]  \nonumber  \\
     & \tneg  \le  & \tneg  c_\k + c_0 |x|^2
     \neg + \neg c_0   E   \bigg[ \Big(\int_t^T  \neg \big| b \big(s, 0, \mu_s, (\beta_\psi (\mu))_s \big)   \big| \, ds \Big)^2
   + \Big(\int_t^T  \neg \big| \si  \big(s, 0, \mu_s, (\beta_\psi (\mu))_s \big)  \big| \, ds \Big)^2 \Big|\cF_t \bigg] ,   \nonumber
     \eea
     where we used a conditional-expectation version of \eqref{eq:esti_X_1a}   in the last inequality.
  Then  an  analogous decomposition and estimation to \eqref{eq:xvx055} leads to that
    $    \big| J \big(t,x,\mu, \beta_\psi (\mu)\big) \big|^p \le c_\k + c_0 |x|^2 $,  \pas    ~
    It follows that
  \beas
      w_1 (t,x) \le   \underset{\mu \in \cU_t }{\esssup} \;
    J \big(t,x,\mu, \beta_\psi (\mu)\big) \le  c_\k + c_0 |x|^{2/p} , \q \pas
   \eeas
   Similarly, one has $ |w_2 (t,x)| \le c_\k  + c_0 |x|^{2/p}  $, \pas \qed

 \ss \no {\bf Proof of Proposition \ref{prop_value_constant}:}  Let  $\cH$ denote the Cameron-Martin
 space of all absolutely continuous functions $h \in \O$ whose derivative $\dot{h}$ belongs to
 $\hL^2([0, T],\hR^d)$. For any $h \in \cH$, we define    $\cT_h(\o) \dfnn \o + h$, $\fa  \o \in \O  $.
 Clearly,  $\cT_h : \O \to \O$ is a bijection and its law is given by
 $ P_h \dfnn P \circ \cT^{-1}_h = \exp \big\{\int_0^T  \dot{h}_s dB_s
 - \frac12 \int_0^T  |  \dot{h}_s  |^2 ds \big\} P$.
     \if{0}
 The Brownian Motion with respect to $ P_h $ is $B_t - \lan B, \int \dot{h}_s dB_s \ran_t = B_t -h(t) $:
 For any $\cD \in sB(\O )$,
 \beas
 P_h \big\{\o \in \O : B(\o)-h \in \cD \big\}
 &=& P   \big\{\o \in \O : B \big( \cT_h(\o) \big)-h \in \cD \big\}
 = P   \big\{\o \in \O : B \big( \o + h \big)-h \in \cD \big\}    \\
 &=& P   \big\{\o \in \O :  B( \o ) + h  - h \in \cD \big\}
 = P   \big\{\o \in \O :   B( \o )   \in \cD \big\} .
 \eeas
    \fi
    Fix $(t,x)  \in [0,T] \times \hR^k $   and set
   $\cH_t \dfnn \{h \in \cH: h(s) = h(s \land t), \, \fa s \in [0,T] \}$.

 \ss \no  {\bf a)}
 { \it Let $h \neg \in \neg  \cH_t $.   We first show that
 \bea   \label{eq:xvx333}
    \big( \mu(\cT_h), \nu(\cT_h) \big)      \in     \cU_t \times \cV_t , \q \fa (\mu,\nu) \in \cU_t  \times \cV_t  .
  \eea }

  Fix  $\mu  \neg \in \neg  \cU_t$. Given  $s  \neg \in \neg  [t,T]$, we   set
  $ \U^h_s(\cD)  \neg \dfnn \neg   \big\{(r, \o)  \neg \in \neg  [t, s]  \neg \times \neg  \O  \neg  : \big(r, \cT_h (\o)\big) \neg \in \neg  \cD \big\}$
  for any $\cD  \neg \subset \neg  [t,s]  \neg \times \neg  \O$.
 As the mapping
 \bea   \label{eq:xvx314}
  \hb{$\cT_h  \neg = \neg  B  \neg + \neg  h $    is $\cF_s / \cF_s -$measurable,}
  \eea
 \if{0}
 (Clearly,    the identity mapping $B : \O \to \O    $
 and the  constant mapping  $ h  : \O \to \O  $  are    both $\cF_s / \cF_s -$measurable.)
  \fi
   it holds for any  $\cE  \neg \in \neg  \sB \big([t,s]\big)$ and $A  \neg \in \neg  \cF_{s}$ that
    \beas 
    \U^h_s(\cE \times A) =  \big\{(r, \o) \in [t, s] \times \O: \big(r, \cT_h (\o)\big) \in \cE \times A \big\}
    =  \big( \cE \cap [t, s] \big) \times \cT^{\,-1}_h(A) \in \sB\big([t,s]\big) \otimes \cF_{s} .
    \eeas
     So  $ \cE \times A \in \L^h_{s} \dfnn \big\{ \cD \subset [t,s] \times \O : \U^h_s ( \cD ) \in \sB\big([t,s]\big) \otimes \cF_{s} \big\}$. In particular, $\es \times \es \in \L^h_{s}$ and $[t,s] \times \O \in \L^h_{s}  $.
              For any $\cD \in \L^h_{s}$ and $\{\cD_n\}_{n \in \hN} \subset \L^h_{s}$,   one can deduce that
             \beas
  \U^h_s \big( ( [t,s] \times \O) \backslash \cD \big) & \tneg \dneg   = & \tneg  \dneg  \big\{(r, \o) \in [t, s] \times \O: \big(r, \cT_h (\o)\big) \in ( [t,s] \times \O) \backslash \cD \big\}    \nonumber \\
    & \tneg  \dneg  =&   \tneg \dneg  ( [t,s] \times \O) \backslash  \big\{(r, \o) \in [t, s] \times \O: \big(r, \cT_h (\o)\big) \in   \cD \big\}   =   ( [t,s] \times \O) \backslash  \U^h_s(\cD)  \in \sB\big([t,s]\big) \otimes \cF_{s} , \\
 \hb{and } \;      \U^h_s  \Big( \underset{n \in \hN}{\cup} \cD_n \Big)
     & \tneg \dneg  = & \tneg  \neg   \big\{(r, \o) \in [t, s] \times \O: \big(r, \cT_h (\o)\big)  \in \underset{n \in \hN}{\cup} \cD_n  \Big\}  \nonumber \\
    & \tneg \dneg  = & \tneg \neg  \underset{n \in \hN}{\cup}  \big\{(r, \o) \in [t, s] \times \O: \big(r, \cT_h (\o)\big) \in    \cD_n  \big\}
 =  \underset{n \in \hN}{\cup} \U^h_s ( \cD_n )   \in \sB\big([t,s]\big) \otimes \cF_{s}  ,
     \eeas
     i.e.  $ ( [t,s] \times \O) \backslash \cD $, $ \underset{n \in \hN}{\cup} \cD_n \in \L^h_s$.
     Thus    $\L^h_{s}$ is   a $\si-$field of $[t,s] \times \O$.
      It follows that
      \bea   \label{eq:xvx211}
       \sB([t,s]) \otimes \cF_{s} = \si \big\{\cE \times A:    \cE \in \sB \big([t,s]\big), \,    A \in \cF_{s} \big\} \subset  \L^h_{s}  .
       \eea
         Given    $U \in \sB(\hU) $, the $\bF-$progressive  measurability of $\mu$ and \eqref{eq:xvx211} show that
  \beas
 \cD_U \dfnn \big\{ (r,  \o) \in [t,s] \times \O: \, \mu_r( \o ) \in U \big\} \in \sB([t,s]) \otimes \cF_s \subset \L^h_{s}.
  \eeas
  That is
  \bea   \label{eq:xvx323}
  \big\{(r, \o)  \neg \in \neg  [t, s]  \neg \times \neg  \O  \neg  :
     \mu_r \big( \cT_h (\o) \big) \in U      \big\} = \big\{(r, \o)  \neg \in \neg  [t, s]  \neg \times \neg  \O  \neg  :
     \big(r, \cT_h (\o)\big)  \in  \cD_U       \big\} =
   \U^h_s \big(  \cD_U \big) \in \sB([t,s]) \otimes \cF_s ,
  \eea
  which shows  the $\bF-$progressive measurability of  process $\mu(\cT_h)$.

 \ss Suppose that $E \int_t^T    [\mu_s]^q_{\overset{}{\hU}} ds < \infty$ for some $q > 2$.
  Then   one can deduce  from H\"older's inequality  that for any    $\wt{q} \in (2, q)  $
  \bea
 && \hspace{-0.7cm}  E \int_t^T    \big[  \mu_s (\cT_h)\big]^{\wt{q}}_{\overset{}{\hU}} \, ds
   =  E_{P_h} \int_t^T    \big[  \mu_s  \big]^{\wt{q}}_{\overset{}{\hU}} \, ds
   =  E  \bigg[ \exp \bigg\{\int_0^T  \dot{h}_s dB_s
 - \frac12 \int_0^T  |  \dot{h}_s  |^2 ds \bigg\} \int_t^T    \big[  \mu_s  \big]^{\wt{q}}_{\overset{}{\hU}} \, ds \bigg] \nonumber \\
 &&  \le T^{\frac{q-\wt{q}}{q}} \exp\bigg\{\frac{\wt{q}}{2(q-\wt{q})} \int_0^T  |  \dot{h}_s  |^2 ds \bigg\}  E  \bigg[ \exp \bigg\{\int_0^T  \dot{h}_s dB_s
 - \frac{q}{2(q-\wt{q})} \int_0^T  |  \dot{h}_s  |^2 ds \bigg\} \Big( \int_t^T    \big[  \mu_s  \big]^q_{\overset{}{\hU}} \, ds \Big)^{ \frac{ \wt{q}}{q} } \bigg]   \nonumber \\
 &&  \le T^{\frac{q -\wt{q}}{q}} \exp\bigg\{\frac{\wt{q}}{2(q \neg - \neg \wt{q})}  \neg  \int_0^T  \neg  |  \dot{h}_s  |^2 ds \bigg\} \Bigg( E  \bigg[ \exp \bigg\{\frac{q}{q \neg - \neg \wt{q}} \int_0^T  \neg  \dot{h}_s dB_s
  \neg - \neg  \frac{q^2}{2(q \neg - \neg \wt{q})^2} \neg \int_0^T  \neg  |  \dot{h}_s  |^2 ds \bigg\}  \bigg] \Bigg)^{ \neg  \frac{q-\wt{q}}{q}}   \bigg( E  \neg  \int_t^T   \neg   \big[  \mu_s  \big]^{q}_{\overset{}{\hU}} \, ds \bigg)^{ \frac{ \wt{q}}{q}}   \nonumber \\
 &&  = T^{\frac{q-\wt{q}}{q}} \exp\bigg\{\frac{\wt{q}}{2(q-\wt{q})} \int_0^T  |  \dot{h}_s  |^2 ds \bigg\}
  \bigg( E  \neg   \int_t^T   \neg   \big[  \mu_s  \big]^{q}_{\overset{}{\hU}} \, ds      \bigg)^{  \frac{ \wt{q}}{q} }
  < \infty .   \label{eq:xvx319}
  \eea
  Hence, $\mu (\cT_h) \in \cU_t$.
  Similarly, $ \nu(\cT_h) \in  \cV_t $ for any $\nu \in \cV_t$.

  \ss \no {\bf b)} { \it We next show that
  \bea   \label{eq:xvx413}
 J(t, x, \mu,\nu)(\cT_h)    =  J \big(t, x, \mu(\cT_h), \nu(\cT_h)\big),  \q  \pas
 \eea }

 Let $\{\F_s\}_{s \in [t,T]}$ be an $ \hR^{k \times d}-$valued, $\bF-$progressively measurable process and set
 $M_s \dfnn \int_t^s \F_r dB_r$, $ s \in [t,T]$.
 We know that (see e.g. Problem 3.2.27 of    \cite{Kara_Shr_BMSC}, which is proved on page 228 therein)
    there exists a sequence of $ \hR^{k \times d}-$valued, $ \bF -$simple processes
   $  \Big\{\F^n_s = \sum^{\ell_n}_{i=1} \xi^n_i \, \b1_{ \big\{s \in (t^n_i, t^n_{i+1}] \big\} } ,
  \,   s \in [t,T] \Big\}_{n \in \hN}$ \big(where $t=t^n_1< \cds< t^n_{\ell_n+1}=T$
   and $\xi^n_i  \in   \cF_{  t^n_i}$  for $i=1, \cds, \ell_n$\big) such that
      \beas
     P \neg  - \neg  \lmt{n \to \infty} \int_t^T     trace\Big\{ \big(  \F^n_r -  \F_r  \big)
     \big(    \F^n_r   -   \F_r  \big)^T  \Big\} \,   ds =0 \q  \hb{and} \q
       P  \neg - \neg   \lmt{n \to \infty} \, \underset{s \in [t,T]}{\sup} \left|  M^n_s  - M_s  \right| = 0   ,
    \eeas
   where $  M^n_s \dfnn \int_t^s \F^n_r dB_s = \sum^{\ell_n}_{i=1}  \xi^n_i  \big(  B_{s \land t^n_{i+1}} -  B_{s \land t^n_i} \big)  $.  By the equivalence of $P_h$ to $ P $,  one has
     \bea
     P_h  \neg  - \neg  \lmt{n \to \infty} \int_t^T     trace\Big\{ \big(  \F^n_r -  \F_r  \big)
     \big(    \F^n_r   -   \F_r  \big)^T  \Big\} \,   ds  & \tneg \dneg = & \tneg \dneg
      P_h  \neg - \neg   \lmt{n \to \infty} \, \underset{s \in [t,T]}{\sup} \left|  M^n_s  - M_s  \right| =0 ,
     \nonumber \\
   \hb{or } \,   P \neg  - \neg  \lmt{n \to \infty} \int_t^T   \neg  trace\Big\{ \big(  \F^n_r (\cT_h)
   \neg - \neg  \F_r  (\cT_h) \big)
     \big(    \F^n_r (\cT_h)  \neg  -  \neg   \F_r (\cT_h) \big)^T  \Big\} \,   ds    & \tneg \dneg  = & \tneg \dneg
  P  \neg - \neg   \lmt{n \to \infty} \, \underset{s \in [t,T]}{\sup} \big|  M^n_s (\cT_h)
  \neg - \neg M_s (\cT_h) \big|  \neg = \neg 0 . \qq \q
    \label{eq:xvx011}
    \eea
    Applying Proposition 3.2.26 of    \cite{Kara_Shr_BMSC} yields that
    \bea   \label{eq:p417}
      0  =     P   \neg - \neg   \lmt{n \to \infty}  \,  \underset{s \in [t,T]}{\sup}
        \Big|  \int_t^s \F^n_r (\cT_h) dB_r
         - \int_t^s \F_r (\cT_h) dB_r  \Big|    .
    \eea
    As $h \in \cH_t $, one can deduce that
    \beas
     M^n_s (\cT_h) & \tneg  = &  \tneg     \Big( \sum^{ \ell_n}_{i=1}    \xi^{\,n}_{\,i}
     \big(B_{s   \land   t^{\,n}_{\, i+1} }
       \neg - \neg  B_{ s   \land   t^{\,n}_{\, i} } \big) \Big)(\cT_h)
       =       \sum^{ \ell_n}_{i=1}    \xi^{\,n}_{\,i} (\cT_h)
     \big(B_{s   \land   t^{\,n}_{\, i+1} } (\cT_h)
       \neg - \neg  B_{ s   \land   t^{\,n}_{\, i} }  (\cT_h)  \big)  \\
     &  \tneg  =  &   \tneg       \sum^{ \ell_n}_{i=1}    \xi^{\,n}_{\,i} (\cT_h)
     \big(B_{s   \land   t^{\,n}_{\, i+1} }   \neg - \neg h(s   \land   t^{\,n}_{\, i+1})
       \neg - \neg  B_{ s   \land   t^{\,n}_{\, i} }    + h(s   \land   t^{\,n}_{\, i}) \big)
     =    \int_t^s \F^n_r (\cT_h) dB_r   , \q \fa s \in [t,T]   ,
    \eeas
    which together with \eqref{eq:xvx011} and \eqref{eq:p417} leads to  that \pas
     \bea  \label{eq:xvx317}
     \int_t^s \F_r (\cT_h) dB_r = M_s(\cT_h) = \Big( \int_t^s \F_r dB_r \Big) (\cT_h)      ,   \q   s \in [t,T]   .
     \eea

    Let   $(\mu,\nu) \neg \in \neg  \cU_t  \neg \times \neg  \cV_t $ and set $\Th \neg = \neg (t,x,\mu,\nu)$. By \eqref{eq:xvx314}, the process $X^\Th (\cT_h)$    is $\bF-$adapted, and the  equivalence of $P_h$ to $ P $ implies that $X^\Th (\cT_h)$ has \pas ~ continuous  paths.
    Suppose that $E  \neg  \int_t^T    [\mu_s]^q_{\overset{}{\hU}} ds
     \neg + \neg  E  \neg  \int_t^T    [\nu_s]^q_{\overset{}{\hV}} ds  \neg < \neg  \infty$
    for some $q  \neg > \neg  2$.
    A standard   estimate  of SDEs   (see e.g.  \cite[pg. 166-168]{Ikeda_Watanabe_1989}
     and \cite[pg. 289-290]{Kara_Shr_BMSC}) shows that
 \bea
   E \neg \left[  \underset{s \in [t,T]}{\sup}  \big| X^\Th_s  \big|^q \right]
  & \dneg \le &  \dneg
      c_q      |x|^q   \neg + \neg c_q   E   \left[ \Big(\int_t^T  \neg \big| b^{\mu,\nu}     (s, 0 )   \big| \, ds \Big)^q
   + \Big(\int_t^T  \neg \big| \si^{\mu,\nu}     (s, 0 )   \big| \, ds \Big)^q \right] \nonumber  \\
    & \dneg  \le  & \dneg
           c_q  \bigg( 1  \neg + \neg   |x|^q    \neg + \neg  E  \neg \int_t^T  \dneg \big( [ \mu_s ]^q_{\overset{}{\hU}} \neg + \neg  [ \nu_s ]^q_{\overset{}{\hV}}  \big)  ds  \bigg) < \infty . \label{eq:xvx321}
 \eea
 Similar to \eqref{eq:xvx319}, one can deduce that   $ E   \bigg[  \underset{s \in [t,T]}{\sup}  \big| X^\Th_s (\cT_h) \big|^{\wt{q}} \bigg] \neg < \neg  \infty$ for any $ \wt{q }  \neg \in \neg  [2,q) $. In particular,
   $ X^\Th (\cT_h)  \neg \in \neg  \hC^2_\bF([t,T],\hR^k) $.
  It follows from   \eqref{eq:xvx317} that
 \beas
   X^\Th_s (\cT_h) & \tneg  \neg = &  \neg \tneg  x \neg +  \neg \int_t^s b \big(r, X^\Th_r(\cT_h),\mu_r(\cT_h),\nu_r(\cT_h) \big) \, dr  \neg + \neg
    \bigg( \int_t^s \si(r, X^\Th_r,\mu_r,\nu_r) \, dB_r \bigg) (\cT_h) \\
    & \tneg \neg  = & \tneg  \neg  x \neg  + \neg  \int_t^s b \big(r, X^\Th_r(\cT_h),\mu_r(\cT_h),\nu_r(\cT_h) \big) \, dr  \neg + \neg
      \int_t^s \si \big(r, X^\Th_r(\cT_h),\mu_r(\cT_h),\nu_r(\cT_h) \big) \, dB_r   ,  \q s \in [t, T] .
   \eeas
 Thus  the uniqueness of  SDE \eqref{FSDE} with parameters $\Th_h = \big(t, x, \mu (\cT_h),\nu (\cT_h) \big)$
 shows that
  \bea    \label{eq:xvx331}
 X^{\Th_h}_s = X^\Th_s (\cT_h)  , \q \fa s \in [t, T]  .
  \eea

  Let $\big(\wh{Y},\wh{Z}\big) = \big(Y^\Th \big(T, g (X^\Th_T) \big) , Z^\Th \big(T, g (X^\Th_T) \big)  \big)$.
  Analogous to   $X^\Th (\cT_h)$, $\wh{Y} (\cT_h)$ is an $\bF-$adapted continuous process.  And using
  the similar arguments that leads to \eqref{eq:xvx323}, we see that the process $\wh{Z} (\cT_h)$ is $\bF-$progressively measurable.           By \eqref{eq:xvx321},
  $ g \big( X^\Th_T   \big) \in \hL^{\frac{pq}{2}} (\cF_T) 
   $,   and a similar argument to \eqref{eq:s031} yields that
  \beas
      E    \bigg[      \Big(\int_t^T  \neg \big| f^\Th_T   (s,0, 0 )   \big|  ds \Big)^{\frac{pq}{2}} \bigg]
   \neg  \le  \neg  c_q    \neg + \neg  c_q   E \bigg[   \underset{s \in [t,T]}{\sup}   \big| X^\Th_s \big|^q
   \neg + \neg \int_t^T \dneg \big( [ \mu_s ]^q_{\overset{}{\hU}} \neg+\neg [ \nu_s ]^q_{\overset{}{\hV}}  \big)  ds \bigg]
    \neg < \neg  \infty .
   \eeas
   Then we know from Proposition \ref{BSDE_well_posed}   that the unique solution $\big(\wh{Y},\wh{Z}\big)$
      of BSDE$\big(t,  g (X^\Th_T), f^\Th_T \big)$ in $\hG^p_\bF([t,T])$ actually belongs to  $\hG^{\frac{pq}{2}}_\bF([t,T])$.
   Similar to \eqref{eq:xvx319}, one can deduce that
     $
      E   \bigg[  \underset{s \in [t,T]}{\sup}  \big| \wh{Y}_s (\cT_h) \big|^{\wt{q}}
      + \bigg( \int_t^T  \big| \wh{Z}_s (\cT_h) \big|^2 ds \bigg)^{\wt{q}/2} \, \bigg] \neg < \neg  \infty  $
      for any $ \wt{q }  \neg \in \neg  \big[p,\frac{pq}{2} \big)$.
            In particular,
   $  \big( \wh{Y} (\cT_h), \wh{Z} (\cT_h) \big)  \neg \in \neg  \hG^p_\bF([t,T]) $.

 \ss   Applying \eqref{eq:xvx317} again, we can deduce from \eqref{eq:xvx331} that
  \beas
   \wh{Y}_s (\cT_h) & \tneg  \neg = &  \neg \tneg  g \big(X^\Th_T (\cT_h) \big)  \neg +  \neg \int_s^T f \big(r, X^\Th_r(\cT_h),\wh{Y}_r (\cT_h),\wh{Z}_r (\cT_h),\mu_r(\cT_h),\nu_r(\cT_h) \big) \, dr  \neg - \neg
    \bigg( \int_s^T  \wh{Z}_r \, dB_r \bigg) (\cT_h) \\
    & \tneg \neg  = & \tneg  \neg  g \big(X^{\Th_h}_T \big)  \neg +  \neg \int_s^T f \big(r, X^{\Th_h}_r ,\wh{Y}_r (\cT_h),\wh{Z}_r (\cT_h),\mu_r(\cT_h),\nu_r(\cT_h) \big) \, dr  \neg - \neg
      \int_t^s \wh{Z}_r (\cT_h)   \, dB_r   ,  \q s \in [t, T] .
   \eeas
   Thus  the uniqueness of  BSDE$\big(t,  g (X^{\Th_h}_T), f^{\Th_h}_T \big)$ implies that  \pas
   \beas
   Y^{\Th_h}_s \big(T, g \big(X^{\Th_h}_T \big) \big)
 =   \wh{Y}_s (\cT_h)   , \q s \in [t, T]  .
  \eeas
 \if{0}
 \beas
 \Big( Y^{\Th_h}_s \big(T, g \big(X^{\Th_h}_T \big) \big), Z^{\Th_h}_s \big(T, g \big(X^{\Th_h}_T \big) \big) \Big)
 = \big( \wh{Y}_s (\cT_h) , \wh{Z}_s (\cT_h) \big) , \q s \in [t, T]  .
  \eeas
\beas
 P \Big( Y^{t,x,\mu,\nu}_s (\cT_h) = Y^{t,x,\mu(\cT_h),\nu(\cT_h)}_s , \; \fa s \in [t, T]\Big)=1; \\
 Z^{t,x,\mu,\nu}_s (\cT_h) = Z^{t,x,\mu(\cT_h),\nu(\cT_h)}_s , ds \times dP-a.s. \hb{ on } [t, T] \times  \O;
\eeas
 \fi
 In particular,
 \beas
 J(t, x, \mu,\nu)(\cT_h) = \wh{Y}_t (\cT_h) = Y^{\Th_h}_t \big(T, g \big(X^{\Th_h}_T \big) \big)
  =  J \big(t, x, \mu(\cT_h), \nu(\cT_h)\big),  \q  \pas
 \eeas

 \ss \no {\bf c)} { \it Now, we show that $ w_1(t, x)(\cT_h)
   \neg  = \neg  w_1(t, x) $,     \pas }

 \ss Let $\beta \in \fB_t$ and  define
  \beas
    \beta_h(\mu) \dfnn \beta \big(\mu(\cT_{-h})\big)(\cT_h)   , \q  \fa \mu \in \cU_t  .
   \eeas
   Similar to \eqref{eq:xvx333},     $\mu(\cT_{-h}) \in \cU_t$ as $-h $ also belongs to $\cH$.
    It follows that $\beta \big(\mu(\cT_{-h})\big) \in \cV_t $. Using  \eqref{eq:xvx333} again
    shows  that $ \beta_h(\mu) = \beta \big(\mu(\cT_{-h})\big)(\cT_h) \in \cV_t $.
   Since
   $    \big[ \big(\beta \big(\mu(\cT_{-h})\big) \big)_s \big]_{\overset{}{\hV}} \le \k+ C_\beta  [   \mu_s(\cT_{-h})  ]_{\overset{}{\hU}} $, $   ds \times dP -$a.s.,  the equivalence of $P_h$ to $ P $ shows that
   $\big[ \big(\beta \big(\mu(\cT_{-h})\big) \big)_s \big]_{\overset{}{\hV}}
      \le     \k+ C_\beta  [   \mu_s(\cT_{-h})  ]_{\overset{}{\hU}} $, $ ds \times dP_h -$a.s., or
     \beas
        \big[ \big(\beta_h(\mu)\big)_s \big]_{\overset{}{\hV}}
     =     \big[ \big(\beta \big(\mu(\cT_{-h})\big)(\cT_h)\big)_s \big]_{\overset{}{\hV}}
     \le \k+ C_\beta  [  \mu_s ]_{\overset{}{\hU}} \, , \q ds \times dP -a.s.
   \eeas

   Let    $\mu^1, \mu^2  \neg \in \neg  \cU_t$ such that
 $\mu^1  \neg = \neg   \mu^2   $,  $ds \times d P -$a.s.  on $ \[t,\t\[ \, \cup \, \[\t, T\]_A   $
 for some $ \t \in   \cS_{t,T}$  and $ A  \in \cF_\t  $. By the equivalence of $P_{-h}$ to $ P $,
 $\mu^1  \neg = \neg   \mu^2   $,  $ds \times d P_{-h} -$a.s.  on $ \[t,\t\[ \, \cup \, \[\t, T\]_A   $,
 or  $\mu^1 (\cT_{-h}) \neg = \neg   \mu^2 (\cT_{-h}) $,  $ds \times d P -$a.s.
 on $ \[t,\t (\cT_{-h}) \[ \, \cup \, \[\t (\cT_{-h}) , T\]_{\cT_h (A)  }   $.
 Given $s \in [t,T]$, similar to \eqref{eq:xvx314}, $\cT_{-h}  $  is also $\cF_s / \cF_s -$measurable.   It follows that
  \beas
  \{ \t (\cT_{-h}) \le s \} & \tneg \dneg = & \tneg \dneg  \big\{\o:  \cT_{-h} (\o) \in \{ \t   \le s \} \big\}=
  \cT^{-1}_{-h}\big(\{ \t   \le s \}  \big) \in \cF_s \\
 \hb{and }  \cT_h (A) \cap \{ \t (\cT_{-h}) \le s \} & \tneg \dneg = & \tneg \dneg
 \cT^{-1}_{-h} (A) \cap \cT^{-1}_{-h}\big(\{ \t   \le s \}  \big) = \cT^{-1}_{-h}\big( A \cap \{ \t   \le s \}  \big)
 \in \cF_s ,
   \eeas
   which shows that $\t (\cT_{-h})  $ is an $\bF-$stopping time and $ \cT_h (A) \in \cF_{\t (\cT_{-h})} $.
   As $t \le  \t \le T$, \pas, the equivalence of $P_{-h}$ to $ P $ shows that $t \le  \t \le T$, $P_{-h}-$a.s., or
   $t \le  \t (\cT_{-h}) \le T$, \pas~ So   $\t (\cT_{-h}) \in \cS_{t,T}$, and we see from
    Definition \ref{def:strategy} that  $\beta \big( \mu^1 (\cT_{-h}) \big)
   \neg = \neg  \beta \big(  \mu^2 (\cT_{-h}) \big) $,  $ds \times d P -$a.s.
 on $ \[t,\t (\cT_{-h}) \[ \, \cup \, \[\t (\cT_{-h}) , T\]_{\cT_h (A)  }   $.
   The equivalence of $P_h$ to $ P $ then shows that
 $\beta \big( \mu^1 (\cT_{-h}) \big)
   \neg = \neg  \beta \big(  \mu^2 (\cT_{-h}) \big) $,  $ds \times d P_h -$a.s.
 on $ \[t,\t (\cT_{-h}) \[ \, \cup \, \[\t (\cT_{-h}) , T\]_{\cT_h (A)  }   $, or
 $ \beta_h(\mu^1) = \beta \big( \mu^1 (\cT_{-h}) \big) (\cT_h)
   \neg = \neg  \beta \big(  \mu^2 (\cT_{-h}) \big) (\cT_h) =  \beta_h(\mu^2)  $,  $ds \times d P -$a.s.
 on $ \[t,\t   \[ \, \cup \, \[\t   , T\]_A     $. Hence,  $\beta_h \in \fB_t $.

 \ss Set $I(t, x, \beta) \dfnn \esup{\mu\in\cU_t} J\big(t, x, \mu, \beta(\mu)\big)$.
 For any $ \mu \in \cU_t$, as $I(t, x, \beta)
 \ge J\big(t, x,  \mu , \beta(  \mu )\big)$, \pas, ~ the equivalence of $P_h$ to $ P $ shows that
   $I(t, x, \beta)
 \ge J\big(t, x,  \mu , \beta(  \mu )\big)$, $P_h-$a.s., or
 \bea   \label{eq:xvx411}
  I(t, x, \beta)(\cT_h) \ge J\big(t, x, \mu, \beta(\mu)\big)(\cT_h) , \q  \pas
  \eea
    Let  $\xi$ be another  random variable   such that
 $\xi \ge J\big(t, x, \mu, \beta(\mu)\big)(\cT_h)$, \pas, or  $\xi (\cT_{-h})  \ge J\big(t, x, \mu, \beta(\mu)\big) $,
 $P_h-$a.s.  for any $ \mu \in \cU_t$. By the equivalence of $P_h$ to $ P $, it holds for any $ \mu \in \cU_t$ that
 $\xi (\cT_{-h})  \ge J\big(t, x, \mu, \beta(\mu)\big) $, \pas ~ Taking essential supremum over
  $ \mu \in \cU_t$ yields that  $\xi (\cT_{-h}) \ge I(t,x,\beta) $, \pas ~ or
 $\xi  \ge I(t,x,\beta) (\cT_h) $, $P_{-h}-$a.s. Then it follows from the   equivalence of $P_{-h}$ to $ P $ that
 $\xi  \ge I(t,x,\beta) (\cT_h) $, \pas, which together with \eqref{eq:xvx411} implies that
 \bea  \label{eq:xvx415}
     \esup{\mu\in\cU_t} \Big( J\big(t, x, \mu, \beta(\mu)\big)(\cT_h) \Big)
     = I(t,x,\beta) (\cT_h) = \Big( \esup{\mu\in\cU_t}  J\big(t, x, \mu, \beta(\mu)\big)  \Big) (\cT_h) ,  \q   \pas
  \eea
  Similarly,     $
     \einf{\beta \in \fB_t} \Big( I (t, x,  \beta )(\cT_h) \Big)
     \neg  = \neg \Big(  \einf{\beta \in \fB_t}  I  (t, x,   \beta )  \Big) (\cT_h)   $, \pas,  which together \eqref{eq:xvx413}
       and \eqref{eq:xvx415} yields that
 \bea
    w_1(t, x)(\cT_h) & \tneg \neg = & \tneg \neg  \Big(  \einf{\beta \in \fB_t}  I  (t, x,   \beta )  \Big) (\cT_h)
 =  \einf{\beta \in \fB_t} \Big( I (t, x,  \beta )(\cT_h) \Big)
 =\einf{\beta\in\fB_t} \esup{\mu\in\cU_t} \Big( J\big(t, x, \mu, \beta(\mu)\big)(\cT_h) \Big) \nonumber  \\
 &  \tneg \neg  = & \tneg   \neg \einf{\beta\in\fB_t} \esup{\mu\in\cU_t} J\big(t, x, \mu(\cT_h), \beta_h(\mu(\cT_h))\big)
   \neg = \neg  \einf{\beta\in\fB_t} \esup{\mu\in\cU_t} J\big(t, x, \mu, \beta_h(\mu)\big)  \nonumber \\
 &  \tneg \neg  =&  \tneg \neg   \einf{\beta\in\fB_t} \esup{\mu\in\cU_t} J\big(t, x, \mu, \beta(\mu)\big)
   \neg  = \neg  w_1(t, x) ,   \q \pas  ,  \label{eq:xvx417}
\eea
 where we  used the facts that $\big\{\mu(\cT_h):  \mu  \in \cU_t \big\} = \cU_t $ and $\big\{\beta_h : \beta \in \fB_t \big\} = \fB_t$.

   \ss \no {\bf d)}  As an $\cF_t-$measurable random variable, $w_1(t, x)$ only depends  on the restriction of
 $\o \in \O  $ to the time interval $[0, t]$. So \eqref{eq:xvx417} holds even for any $h \in \cH$.
 Then an application of Lemma 3.4 of \cite{Buckdahn_Li_1} yields that $ w_1(t,x)=  E[w_1(t,x)]  $, \pas~Similarly,
 one can deduce that $ w_2(t,x)=  E[w_2 (t,x)]  $, \pas \qed

    \ss \no {\bf Proof of Proposition \ref{prop_w_conti}:}
  Let $t \in [0,T]$ and  $x_1, x_2 \in   \hR^k $.    For any $(\beta,\mu) \in \fB_t \times \cU_t  $,
   \eqref{eq:s025} implies that
  \beas   
   \Big|  J \big( t,x_1,\mu,\beta (  \mu ) \big)
     - J \big( t,x_2,\mu,\beta (  \mu ) \big)   \Big|^p
      \le c_0 | x_1 \neg - \neg x_2 |^2 , \q \pas
  \eeas
  which leads to that
   \beas
    J \big( t,x_2,\mu,\beta (  \mu ) \big)
  \neg  - \neg  c_0 | x_1 \neg - \neg x_2 |^{2/p}
    \neg  \le \neg  J \big( t,x_1,\mu,\beta (  \mu ) \big)
         \neg  \le \neg  J \big( t,x_2,\mu,\beta (  \mu ) \big)
     \neg + \neg  c_0 | x_1 \neg - \neg x_2 |^{2/p} , \q \pas    
  \eeas
   Taking essential supremum over $\mu \in \cU_{t }$ and then taking essential infimum over $\beta \in \fB_{t }$ yield that
  \beas
   w_1(t,x_2)    -    c_0 | x_1 \neg - \neg x_2 |^{2/p}  \le    w_1(t,x_1)
            \le    w_1(t,x_2)   +    c_0 | x_1 \neg - \neg x_2 |^{2/p}    .
  \eeas
  So $ \big|   w_1(t,x_1) \neg - \neg  w_1(t,x_2) \big| \le c_0 | x_1 \neg - \neg x_2 |^{2/p} $.
  Similarly, one has  $ \big|   w_2(t,x_1) \neg - \neg  w_2(t,x_2) \big|
  \le c_0 | x_1 \neg - \neg x_2 |^{2/p} $.  \qed

  \subsection{Proof of  the Weak Dynamic Programming Principle}

  \label{subsection:DPP}

  To prove the weak dynamic programming principle (Theorem \ref{thm_DPP}), we begin with two auxiliary result. The first one
 shows that the pasting of  state processes (resp. payoff processes) is
 exactly the state process (resp. payoff process) with the pasted controls.

  \begin{lemm}   \label{lem_partition}
 Given $t \in [0,T]$, let $\{A_i\}^n_{i=1} \subset \cF_t$ be a partition of $\O$.
 For any   $ \big\{ (\xi_i,\mu^i,\nu^i) \big\}^n_{i=0} \subset   \hL^2 (\cF_t, \hR^k) \times \cU_t \times \cV_t $, if
 $  \xi_0 =   \sum^n_{i=1} \b1_{A_i} \xi_i $, \pas ~ and if $  ( \mu^0 ,\nu^0 ) =  \big( \sum^n_{i=1} \b1_{A_i} \mu^i, \sum^n_{i=1} \b1_{A_i} \nu^i \big)$, $ds \times dP-$a.s., then  it holds \pas ~ that
 \bea    \label{eq:xxa271}
 X^{t,\xi_0,\mu^0,\nu^0}_s   = \sum^n_{i=1} \b1_{A_i} X^{t,\xi_i, \mu^i, \nu^i }_s , \q \fa s \in [t,T] .
 \eea
  Moreover, for any     $ \{ ( \t_i , \eta_i ) \}^n_{i=0} \subset  \cS_{t,T} \times \hL^p (\cF_T )$ such that
  each $ \eta_i$ is  $ \cF_{\t_i}-$measurable, if $   \t_0  \neg = \neg   \sum^n_{i=1} \b1_{A_i} \t_i    $, \pas ~
   and   if $   \eta_0  \neg = \neg    \sum^n_{i=1} \b1_{A_i} \eta_i  $, \pas, then  it holds \pas ~ that
 \bea   \label{eq:xxa272}
 Y^{t,\xi_0,\mu^0,\nu^0}_s \big(\t_0, \eta_0 \big)   = \sum^n_{i=1} \b1_{A_i} Y^{t,\xi_i, \mu^i, \nu^i }_s
  \big(\t_i, \eta_i \big) , \q \fa s \in [t,T] .
 \eea
   In particular, one has
 \bea   \label{eq:xxa273}
 J( t,\xi_0,\mu^0,\nu^0)   = \sum^n_{i=1} \b1_{A_i} J( t,\xi_i, \mu^i, \nu^i ) ,  \q  \pas
 \eea

 \end{lemm}

\ss \no {\bf Proof:} Let $\big( X^i,   Y^i,  Z^i  \big) = \big( X^{t,\xi_i, \mu^i, \nu^i},   Y^{t,\xi_i, \mu^i, \nu^i}  (\t_i, \eta_i),  Z^{t,\xi_i, \mu^i, \nu^i}  (\t_i, \eta_i)  \big)$ for  $i=0,\cds \neg , n$.  We define
     \beas
      \big(\cX,\cY,\cZ   \big) \dfnn \sum^n_{i=1} \b1_{A_i}
   \big( X^i,      Y^i,  Z^i   \big) \in  \hC^2_\bF([t,T], \hR^k ) \times \hG^p_\bF([t,T]) .
   \eeas
   For any  $s \in [t, T]$ and $i=1,\cds \neg ,n$,  multiplying $\b1_{A_i}$ to  SDE \eqref{FSDE} with parameters $(t,\xi_i, \mu^i, \nu^i)$, we can deduce that
         \bea
  \b1_{A_i} X^i_s &=& \b1_{A_i} \xi_i + \b1_{A_i} \int_t^s  b(r, X^i_r,\mu^i_r,\nu^i_r ) \, dr + \b1_{A_i} \int_t^s \si(r, X^i_r,\mu^i_r,\nu^i_r) \, dB_r  \nonumber \\
   &=& \b1_{A_i} \xi_i +  \int_t^s \b1_{A_i} b(r, X^i_r,\mu^i_r,\nu^i_r ) \, dr +   \int_t^s \b1_{A_i} \si(r, X^i_r,\mu^i_r,\nu^i_r) \, dB_r  \nonumber \\
   &=& \b1_{A_i} \xi_i +  \int_t^s \b1_{A_i} b(r, \cX_r,\mu^0_r,\nu^0_r ) \, dr +   \int_t^s \b1_{A_i} \si(r, \cX_r,\mu^0_r,\nu^0_r) \, dB_r   \nonumber  \\
    &=& \b1_{A_i} \xi_i + \b1_{A_i} \int_t^s  b(r, \cX_r,\mu^0_r,\nu^0_r ) \, dr +  \b1_{A_i} \int_t^s \si(r, \cX_r,\mu^0_r,\nu^0_r) \, dB_r  ,  \q \pas   \label{eq:xxa251}
    \eea
  Adding them up over $ i \in \{ 1,\cds \neg ,n \}$ and using the continuity of process $\cX$ show that  \pas
  \beas
  \cX_s  = \xi_0 + \int_t^s  b(r, \cX_r,\mu^0_r,\nu^0_r ) \, dr +    \int_t^s \si(r, \cX_r,\mu^0_r,\nu^0_r) \, dB_r , \q  s \in [t,T].
  \eeas
  So $\cX = X^{t,\xi_0,\mu^0,\nu^0}  $, i.e. \eqref{eq:xxa271}.

 \ss Next, for any  $s \neg \in \neg  [t, T]$ and $i \neg = \neg 1,\cds \neg ,n$,  similar to \eqref{eq:xxa251}, multiplying $\b1_{A_i}$ to   BSDE$\Big(t, \eta_i, f^{t,\xi_i,\mu^i,\nu^i}_{\t_i} \Big)$  yields that
   \beas
      \b1_{A_i} Y^i_s &=& \b1_{A_i} \eta_i  \neg + \neg  \b1_{A_i} \int_s^T  \neg \b1_{\{r < \t_i \}} f (r, X^i_r, Y^i_r,  Z^i_r, \mu^i_r, \nu^i_r)  \, dr
     \neg-\neg \b1_{A_i} \int_s^T \neg Z^i_r d B_r                 \\
    &=& \b1_{A_i} \eta_i  \neg + \neg  \b1_{A_i} \int_s^T  \neg \b1_{\{r < \t_0 \}} f (r, \cX_r,  \cY_r,   \cZ_r, \mu^0_r, \nu^0_r)  \, dr
   \neg-\neg \b1_{A_i} \int_s^T \neg \cZ_r d B_r , \q \pas
    \eeas
     Adding them up and using the continuity of process $ \cY   $, we obtain that  \pas
   \beas
        \cY_s =     \eta_0  \neg + \neg    \int_s^T  \neg \b1_{\{r < \t_0 \}}
        f (r, X^{t,\xi_0,\mu^0,\nu^0}_r,  \cY_r,   \cZ_r, \mu^0_r, \nu^0_r)  \, dr
   \neg-\neg   \int_s^T \neg \cZ_r d B_r  , \q s \in [t,T]  .
    \eeas
  Thus $\big( \cY , \cZ   \big) = \big(  Y^{t,\xi_0,\mu^0,\nu^0}  (\t_0, \eta_0),  Z^{t,\xi_0,\mu^0,\nu^0}  (\t_0, \eta_0)  \big)     $, proving \eqref{eq:xxa272}.

  \ss   Taking $\t_i =T$ and $\eta_i = g \big( X^{t,\xi_i,\mu^i,\nu^i}_T \big) \in \hL^p (\cF_T )$ for $i=0,\cds,n$,
  we see from \eqref{eq:xxa271}   that
 \beas
 \sum^n_{i=1}  \b1_{A_i} \eta_i  = \sum^n_{i=1}  \b1_{A_i} g(X^{t,\xi_i,\mu^i,\nu^i}_T)
 = \sum^n_{i=1}  \b1_{A_i} g \big( X^{t,\xi_0,\mu^0,\nu^0}_T \big) =  g \big( X^{t,\xi_0,\mu^0,\nu^0}_T \big)
 = \eta_0 ,  \q \pas
  \eeas
  Then   \eqref{eq:xxa272} shows that  \pas
  \beas
 \hspace{2.3cm}
 J( t,\xi_0,\mu^0,\nu^0) \neg = \neg  Y^{t,\xi_0,\mu^0,\nu^0}_t   \big(T, \eta_0  \big)
  \neg = \neg  \sum^n_{i=1} \b1_{A_i}  Y^{t,\xi_i,\mu^i,\nu^i}_t   \big(T,\eta_i \big)
  \neg = \neg   \sum^n_{i=1} \b1_{A_i} J( t,\xi_i, \mu^i, \nu^i )   .      \hspace{2.3cm}  \hb{\qed}
 \eeas

 In the next Lemma, we
 approach  $ I(t,x,\beta) \dfnn \underset{\mu \in \cU_t }{\esssup}
  \; Y^{t,x,\mu, \beta (  \mu )  }_t $  from above and
$w_1(t,x) = \underset{\beta \in \fB_t  }{\essinf} \; I(t,x,\beta)  $  from below:

\begin{lemm} \label{lem_mu_lmtu}
 Let $(t,x) \in [0,T] \times \hR^k$  and  $\e>0$. For any $\beta \in \fB_t$, there exist $\big\{(A_n, \mu^n ) \big\}_{n \in \hN} \subset \cF_t \times \cU_t$
     with $ \lmtu{n \to \infty} \b1_{A_n} = 1    $, \pas ~ such that for any $n \in \hN$
 \bea  \label{eq:xxa255}
    J \big( t,x,\mu^n, \beta (\mu^n) \big)  \ge  \big( I  (t,x,   \beta    )  - \e \big) \land \e^{-1}
        , \q    \pas     \hb{ on } A_n   .
 \eea
     Similarly, there exist $\big\{( \cA_n, \beta_n ) \big\}_{n \in \hN} \subset \cF_t \times \fB_t$
     with $ \lmtu{n \to \infty}   \b1_{\cA_n}    =1  $, \pas ~ such that for any $n \in \hN$
 \bea   \label{eq:xxa257}
  w_1(t,x)  \ge  I \big(t,x,   \beta_n   \big)-\e         , \q \pas \hb{ on }    \cA_n .
 \eea

   \end{lemm}

 \ss \no {\bf Proof:}     {\bf (i)}    Let $\beta \in \fB_t$.   Given $\mu^1, \mu^2 \in \cU_t$, we set
 $A \dfnn \{ J \big(t,x,\mu^1,  \beta(\mu^1)\big)  \ge  J \big(t,x,\mu^2,  \beta(\mu^2)\big) \} \in \cF_t$
 and define $\wh{\mu}_s  \dfnn \b1_A \mu^1_s + \b1_{A^c} \mu^2_s $, $s \in [t,T]$.  Clearly,
 $\wh{\mu} $ is an  $\bF-$progressively measurable process. For $i=1,2$, suppose that $ E \neg \int_t^T \neg  \big[ \mu^i_s \big]^{q_i}_{\overset{}{\hU}} \,ds < \infty $ for some $q_i>2$. It follows that
 $    E  \neg  \int_t^T   [ \wh{\mu}_s ]^{q_1 \land q_2}_{\overset{}{\hU}} \,ds
 \le E  \neg  \int_t^T \neg  \big[ \mu^1_s \big]^{q_1 \land q_2}_{\overset{}{\hU}} \,ds
  +E  \neg  \int_t^T \neg  \big[ \mu^2_s \big]^{q_1 \land q_2}_{\overset{}{\hU}} \,ds < \infty$.
    Thus, $\wh{\mu}  \in \cU_t $.    As $ \wh{\mu} = \mu^1 $ on $[t,T] \times A$,
 taking $(\t,A) = (t,A)$ in Definition \ref{def:strategy} yields that $ \beta(\wh{\mu} ) = \beta( \mu^1 )$,
  $ds \times dP-$a.s.  on $[t,T] \times A$. Similarly,  $ \beta(\wh{\mu} ) = \beta( \mu^2 )$,
   $ds \times dP-$a.s.  on $[t,T] \times A^c$.
 So $  \beta(\wh{\mu} )  = \b1_A  \beta( \mu^1 )  + \b1_{A^c}  \beta( \mu^2 ) $, $ds \times dP-$a.s.
   Then   \eqref{eq:xxa273} shows that
   \beas
 \q  J \big(t,x, \wh{\mu}, \beta (\wh{\mu})\big)
 = \b1_A J \big(t,x, \mu^1,  \beta (\mu^1)\big)
   + \b1_{A^c} J \big(t,x, \mu^2,  \beta  (\mu^2)\big)
   = J \big(t,x, \mu^1,  \beta (\mu^1)\big)
   \vee  J \big(t,x, \mu^2,  \beta  (\mu^2)\big) , \q \pas,
   \eeas
 which shows that the collection $ \big\{ J \big(t,x,  \mu, \beta ( \mu )\big) \big\}_{\mu \in \cU_t}$
 is directed upwards (see Theorem A.32 of \cite{Follmer_Schied_2004}).
 By Proposition VI-{\bf 1}-1 of \cite{Neveu_1975} or Theorem A.32 of \cite{Follmer_Schied_2004},
  there exists  a sequence
 $ \big\{ \wt{\mu}^i \big\}_{i \in \hN} \subset \cU_t$ such that
  \bea   \label{eq:xvx068}
   I(t,x,  \beta) = \esup{\mu  \in \cU_t}  J\big(t,x,\mu , \beta (\mu )\big)
   = \lmtu{i \to \infty} J \big(t,x,\wt{\mu}^i, \beta (\wt{\mu}^i)\big), \q \pas
  \eea
  So $ I(t,x,  \beta) $ is $  \cF_t-$measurable.

  For any $i \in \hN$, we set   $ \wt{A}_i \dfnn \big\{J \big(t,x,\wt{\mu}^i, \beta (\wt{\mu}^i)\big)
  \ge  \big( I(t,x,  \beta)   -   \e \big) \land \e^{-1}   \big\}  \in \cF_{t}$   and $\wh{A}_i \dfnn \wt{A}_i  \big\backslash  \underset{j < i }{\cup} \wt{A}_j \in \cF_t$.    
  Fix $n \in \hN$ and set  $A_n \dfnn \underset{i=1}{\overset{n}{\cup}}\wh{A}_i \in \cF_t $.
  Similar to $\wh{\mu}$, $\mu^n \dfnn \sum^n_{i=1} \b1_{\wh{A}_i} \wt{\mu}^i + \b1_{A^c_n} \wt{\mu}^1 $
   also defines a $ \cU_t -$process. 
  For $i=1,\cds, n$,   as $\mu^n = \wt{\mu}^i$ on $ [t,T] \times \wh{A}_i$, taking $(\t,A) = (t, \wh{A}_i)$ in Definition
   \ref{def:strategy} shows that $\beta( \mu^n )= \beta( \wt{\mu}^i )$, $ds \times dP-$a.s. on $ [t,T] \times \wh{A}_i$.
   Then   \eqref{eq:xxa273} implies that
   $
  \b1_{\wh{A}_i} J\big(t,x,\mu^n, \beta (\mu^n)\big)
  \neg = \neg \b1_{\wh{A}_i} J\big(t,x,\wt{\mu}^i, \beta (\wt{\mu}^i)\big)$, \pas ~
 Adding them up over $ i \neg \in \neg \{1, \cds \neg , n\}$ gives
   \beas
  \b1_{A_n} J\big(t,x,\mu^n, \beta (\mu^n)\big) 
  = \sum^n_{i=1} \b1_{\wh{A}_i} J\big(t,x,\wt{\mu}^i, \beta (\wt{\mu}^i)\big)
   \ge   \b1_{A_n} \big( \big( I(t,x,  \beta)  - \e \big) \land \e^{-1} \big) , \q \pas
  \eeas

  Let $\cN$ be the $P-$null set such that \eqref{eq:xvx068} holds on $\cN^c$. Clearly,
  $\{I(t,x,  \beta) < \infty \} \cap \cN^c \subset \underset{i \in \hN}{\cup}   \big\{J \big(t,x,\wt{\mu}^i, \beta (\wt{\mu}^i)\big)  \ge    I(t,x,  \beta)   -   \e       \big\}   $
  and $\{I(t,x,  \beta) = \infty \} \cap \cN^c \subset \underset{i \in \hN}{\cup}   \big\{J \big(t,x,\wt{\mu}^i, \beta (\wt{\mu}^i)\big)  \ge    \e^{-1}       \big\}   $.
  It follows that
  \beas
  \cN^c \subset  \underset{i \in \hN}{\cup}
   \Big( \big\{J \big(t,x,\wt{\mu}^i, \beta (\wt{\mu}^i)\big)
  \ge    I(t,x,  \beta)   -   \e       \big\} \cup \big\{J \big(t,x,\wt{\mu}^i, \beta (\wt{\mu}^i)\big)  \ge    \e^{-1}       \big\} \Big)
  =  \underset{i \in \hN}{\cup}   \wt{A}_i
  =  \underset{i \in \hN}{\cup}\wh{A}_i
  =  \underset{n \in \hN}{\cup} A_n .
  \eeas
  So  $ \lmtu{n \to \infty} \b1_{A_n} = 1    $, \pas

     \ss \no {\bf (ii)}     Let $\beta_1, \beta_2 \in \fB_t$. We just showed  that  $I  (t,x, \beta_1)$ and $I  (t,x, \beta_2)$
  are $\cF_t-$measurable, so $\cA_o \neg \dfnn \neg  \{ I  (t,x, \beta_1)  \neg \le \neg  I  (t,x, \beta_2) \} $
  belongs to $ \cF_t$. For any  $\mu  \neg  \in \neg \cU_t$, similar to $\wh{\mu}$ above,
     $   \beta_o (\mu)    \neg  \dfnn    \neg   \b1_{\cA_o}  \beta_1 (\mu)
     +      \b1_{\cA^c_o}  \beta_2 (\mu) $ defines a $  \cV_t -$process.
   For $i=1,2$, letting $C_i > 0$ be the constant associated to $\beta_i$ in Definition \ref{def:strategy} (i),  we see that
     \beas
    \big[ (\beta_o (\mu))_s \big]_{\overset{}{\hV}}    \neg  =    \neg   \b1_{\cA_o} \big[ ( \beta_1 (\mu) )_s \big]_{\overset{}{\hV}}  +      \b1_{\cA^c_o} \big[ ( \beta_2 (\mu))_s \big]_{\overset{}{\hV}}
    \le \k + ( C_1 \neg \vee \neg C_2 ) \, [\mu_s]_{\overset{}{\hU}} \, , \q ds \times dP-a.s .
     \eeas
        Let    $\mu^1, \mu^2  \neg \in \neg  \cU_t$ such that
 $\mu^1  \neg = \neg   \mu^2 $,  $ds \times d P -$a.s.  on $ \[t,\t\[ \, \cup \, \[\t, T\]_A   $
 for some $ \t \in   \cS_{t,T}$  and $ A  \in \cF_\t  $.
    By Definition \ref{def:strategy},  $\beta_1  (\mu^1)  \neg      =  \neg    \beta_1  (\mu^2)  $ and
   $  \beta_2 (\mu^1)   \neg  =   \neg    \beta_2 (\mu^2)$, $ds  \neg \times \neg  d P -$a.s.
   on $ \[t,\t\[ \, \cup \, \[\t, T\]_A   $. Then it follows that   for  $ds  \neg \times \neg  d P -$a.s.
      $(s,\o) \in \[t,\t\[ \, \cup \, \[\t, T\]_A $
   \bea   \label{eq:xxa303}
    \big(  \beta_o  (\mu^1) \big)_s(\o)
   \neg = \neg  \b1_{\cA_o}   \big( \beta_1  (\mu^1) \big)_s(\o)  \neg  + \neg  \b1_{\cA^c_o} \big( \beta_2  (\mu^1) \big)_s(\o)
     \neg = \neg  \b1_{\cA_o}   \big( \beta_1  (\mu^2) \big)_s(\o) \neg   + \neg \b1_{\cA^c_o} \big( \beta_2  (\mu^2) \big)_s(\o)    \neg = \neg    \big( \beta_o  (\mu^2) \big)_s(\o) . \q
   \eea
    Hence,  $  \beta_o    \in    \fB_t$.

    \ss    For any $\mu \in \cU_t$,     \eqref{eq:xxa273} shows that
   $   J \big(t,x, \mu,  \beta_o (\mu)\big)= \b1_{\cA_o} J \big(t,x, \mu,  \beta_1 (\mu)\big)
   + \b1_{\cA^c_o} J \big(t,x, \mu,  \beta_2 (\mu)\big) $,        \pas     ~
       Then taking essential supremum over $\mu \in \cU_t$ and using Lemma \ref{lem_ess} (2) yield that
   \beas
   I(t,x,   \beta_o )= \b1_{\cA_o} I(t,x,  \beta_1) + \b1_{\cA^c_o} I(t,x,  \beta_2)
   = I  (t,x, \beta_1) \land I  (t,x, \beta_2), \q \pas
   \eeas
   Thus   the collection $\{  I(t,x,  \beta) \}_{\beta \in \fB_t}$ is directed downwards (see Theorem A.32 of \cite{Follmer_Schied_2004}).
  By Proposition VI-{\bf 1}-1 of \cite{Neveu_1975} or Theorem A.32 of \cite{Follmer_Schied_2004},
  one can   find a sequence $ \big\{ \wt{\beta}_i \big\}_{i \in \hN} \subset \fB_{t} $ such that
   \bea  \label{eq:xvx068_b}
   w_1 (t,x) = \einf{\beta \in \fB_t} \, I  (t,x, \beta) = \lmtd{i \to \infty}  I  \big( t,x, \wt{\beta}_i \big),  \q \pas
   \eea

    For any $i \neg \in \neg  \hN$, we set   $ \wt{\cA}_i  \neg \dfnn \neg  \big\{ I  \big( t,x, \wt{\beta}_i \big)
  \neg \le \neg   w_1 (t,x)     +    \e   \big\}   \neg \in \neg  \cF_{t}$   and $\wh{\cA}_i  \neg \dfnn \neg  \wt{\cA}_i  \big\backslash  \underset{j < i }{\cup} \wt{\cA}_j  \neg \in \neg  \cF_t$.
   Fix     $n  \neg \in \neg  \hN$ and set  $\cA_n  \neg \dfnn \neg  \underset{i=1}{\overset{n}{\cup}}\wh{\cA}_i  \neg \in \neg  \cF_t $.
 For any $\mu   \neg \in \neg  \cU_t$,  similar to $\wh{\mu}$ above,
       $   \beta_n (\mu)
  \neg \dfnn \neg   \sum^n_{i = 1 } \b1_{\wh{\cA}_i}   \wt{\beta}_i (\mu)  \neg  + \neg  \b1_{\cA^c_n}  \wt{\beta}_1 (\mu)    $ defines a
  $\cV_t-$process. For $i \neg = \neg 1,\cds \neg , n $, let $\wt{C}_i \neg > \neg  0$ be the constant
  associated to $\wt{\beta}_i$ in Definition \ref{def:strategy} (i).
     Setting $C_n  \neg \dfnn \neg  \max\{\wt{C}_i \neg : i =1,\cds \neg , n \}$,
    we can deduce that
     \beas
    \big[ (\beta_n (\mu))_s \big]_{\overset{}{\hV}}    \neg  =    \neg   \sum^n_{i = 1 } \b1_{\wh{\cA}_i}  \big[ ( \wt{\beta}_i (\mu)  )_s \big]_{\overset{}{\hV}}  +      \b1_{\cA^c_n}   \big[ ( \wt{\beta}_1 (\mu) )_s \big]_{\overset{}{\hV}}
    \le \k + C_n   [\mu_s]_{\overset{}{\hU}} \, , \q ds \times dP-a.s .
     \eeas
      Let    $\mu^1, \mu^2  \neg \in \neg  \cU_t$ such that
 $\mu^1  \neg = \neg   \mu^2 $,  $ds \times d P -$a.s.  on $ \[t,\t\[ \, \cup \, \[\t, T\]_A   $ for some
 $ \t \in   \cS_{t,T}$  and $ A  \in \cF_\t  $.
   Similar to \eqref{eq:xxa303},  it holds  for  $ds \times d P -$a.s.    $(s,\o) \in \[t,\t\[ \, \cup \, \[\t, T\]_A    $
   that
   \beas
    ~ \;   \big(  \beta_n  (\mu^1) \big)_s(\o)
   \neg = \neg  \sum^n_{i = 1 } \b1_{\wh{\cA}_i}    \big( \wt{\beta}_i  (\mu^1) \big)_s(\o)
    \neg + \neg  \b1_{\cA^c_n}  \big( \wt{\beta}_1  (\mu^1) \big)_s(\o)
   \neg   = \neg  \sum^n_{i = 1 } \b1_{\wh{\cA}_i}    \big( \wt{\beta}_i  (\mu^2) \big)_s(\o)
    \neg + \neg  \b1_{\cA^c_n}  \big( \wt{\beta}_1  (\mu^2) \big)_s(\o)
     \neg = \neg    \big( \beta_n  (\mu^2) \big)_s(\o) .
    \eeas
        So  $ \beta_n    \in    \fB_t$.
    For any $\mu \in \cU_t$,  applying \eqref{eq:xxa273} again yields that
  $
   \b1_{\cA_n} J\big(t,x,\mu, \beta_n (\mu)\big) = \sum^n_{i = 1} \b1_{\wh{\cA}_i} J\big(t,x, \mu, \wt{\beta}_i ( \mu )\big)
   $,  \pas  ~
 Taking essential supremum over $\mu \in \cU_t$ and using     Lemma \ref{lem_ess} (2) again yield   that
   \beas
   \b1_{\cA_n}  I  (t,x, \beta_n  ) =    \sum^n_{i = 1} \b1_{\wh{\cA}_i} I \big(t,x,   \wt{\beta}_i  \big)
         \le   \b1_{\cA_n} \big(   w_1 (t,x )  + \e \big)    , \q \pas
  \eeas
  Let $\wt{\cN}$ be the $P-$null set such that \eqref{eq:xvx068_b} holds on $\wt{\cN}^c$.
  As $ |w_1 (t,x)|< \infty  $ by Proposition \ref{prop_value_bounds} and Proposition \ref{prop_value_constant},
  we see that $ \underset{n \in \hN}{\cup} \cA_n = \underset{i \in \hN}{\cup}\wh{\cA}_i
  =  \underset{i \in \hN}{\cup}\wt{\cA}_i  =   \wt{\cN}^c $. \qed

 In the proof of the weak dynamic programming principle below,   we
  first use Lemma \ref{lem_mu_lmtu} to  construct  approximately optimal controls/strategies
  by  pasting locally approximately optimal ones according to a finite partition of $\ol{O}_\d(t,x)$
  determined by the continuity of test functions $\phi$ and $\wt{\phi}$.
 After   a series of estimates on state processes and payoff processes,
 we  obtain  the weak dynamic programming principle
  by using the stochastic backward semigroup  property
  \eqref{eq:p677}, the continuous dependence of  payoff process  on the initial state
    (see Lemma \ref{lem_estimate_Y}) as well as
  the control-neutralizer assumption and the growth condition on strategies.

\ss \no {\bf Proof of Theorem \ref{thm_DPP}:}
 {\bf 1)}      For any $m \in \hN$ and $(s,\fx)  \neg \in \neg  [t,T]    \times    \hR^k$,
   the  continuity of  $ \f  $, $\wt{\f}$ shows that there exists
      a $  \d^m_{s,\fx} \in  ( 0,  1/m   )     $ such that
    \bea   \label{eq:xux391}
     \big| \f (s',\fx')  -  \f (s,\fx) \big| + \big| \wt{\f} (s',\fx')  -  \wt{\f} (s,\fx) \big|  \le    1/m   ,   \q
    \fa (s',\fx')  \in  \big[ (s- \d^m_{s,\fx}) \vee t  , (s + \d^m_{s,\fx}) \land T \big] \times \ol{O}_{  \d^m_{s,\fx} } (\fx) .
    \eea
        By  classical covering theory,   $ \big\{ \fD_m (s,\fx) \dfnn \big(  s- \d^m_{s,\fx} , s + \d^m_{s,\fx}  \big)
    \times O_{  \d^m_{s,\fx} } (\fx)   \big\}_{(s,\fx)  \in  [t,T]  \times  \hR^k}  $ has
    a finite subcollection $ \{ \fD_m (s_i,x_i) \}^{N_m}_{i=1}$ to cover $ \ol{O}_{  \d} (t, x) $.
      For  $  i   =1,\cds \neg, N_m  $, we set $ t_i \dfnn (s_i + \d^m_{s_i,x_i}) \land T $.

 \ss Fix $(\beta, \mu)  \in \fB_t \times \cU_t$  and  simply denote $\t_{\beta,\mu    }$ by $\t$.
 By Lemma \ref{lem_control_combine},   $  \wh{\mu}_s \dfnn  \b1_{ \{s < \t\} } \mu_s + \b1_{ \{s  \ge  \t\} } u_0  $,  $ s \in [t,T] $
 defines   a  $\cU_t-$control. We set   $\Th \dfnn \big(t,x,\mu, \beta ( \mu )\big)$ and
   $\wh{\Th} \dfnn \big(t,x, \wh{\mu}, \beta ( \wh{\mu} )\big)$.
     \if{0}
   Given $s \neg \in \neg  [t,T]$, since  $ \[t, \t\[     \in \dneg  \sP $,
      both $\cD_s   \neg \dfnn  \neg     \[t, \t\[
       \, \cap \,   ( [t,s]  \neg \times \neg  \O) $ and
       $ \cD^c  \neg =   \neg  ( [t,s]  \neg \times \neg  \O)
        \backslash \cD $  belong  to $  \sB\big([t,s]\big)  \neg \otimes \neg  \cF_s $.
       The $\bF-$progressive measurability of   $\mu$ then  implies that  for any $U  \neg \in \neg  \sB \big(\hU  \big)  $
     \bea
    \big\{ (r,  \o) \in \cD_s \neg  :      \wh{\mu}_r (\o)  \neg \in  \neg U \big\}
     &\tneg \dneg = &\tneg \dneg  \big\{ (r,  \o)  \neg \in \neg  \cD_s  \neg   :        \mu_r (\o)  \neg \in \neg  U \big\}
      \neg  =  \neg \cD_s  \neg \cap \neg  \big\{ (r,  \o)  \neg \in \neg  [t,s]  \neg \times \neg  \O  \neg   :    \mu_r (\o)
        \neg \in \neg  U \big\}    \neg  \in  \neg  \sB\big([t,s]\big)  \neg  \otimes \neg   \cF_s   ;   \qq \q
        \label{mu_progressive}  \\
       \big\{ (r,  \o)  \neg  \in \neg    \cD^c_s  \neg   :       \wh{\mu}_r (\o)  \neg  \in \neg   U \big\}
     &\tneg  \dneg = &\tneg \dneg  \big\{ (r,  \o)  \neg  \in \neg   \cD^c_s  \neg  :  u_0  \neg  \in  \neg  U \big\}
       \neg  = \neg   \cD^c_s  \neg  \cap \neg   \big\{ (r,  \o)  \neg  \in \neg   [t,s]  \neg  \times  \neg  \O
        \neg   :    u_0  \neg  \in \neg   U \big\}
        \neg  \in \neg   \sB\big([t,s]\big)  \neg  \otimes \neg   \cF_s   . \nonumber
    \eea
    Putting them together shows the $\bF-$progressive measurability of  $\wh{\mu}$. As
       $
   E \dneg \int_t^T \neg [    \wh{\mu}_s    ]^2_{\overset{}{\hU}} ds
     \neg \le   \neg E \dneg \int_t^T  \neg [     \mu_s    ]^2_{\overset{}{\hU}} ds
  \neg + \neg [    u_0    ]^2_{\overset{}{\hU}} (T \neg -\neg t) \neg     < \neg \infty$,
     we see that $\wh{\mu} \in \cU_t $.
   \fi

  \ss \no  {\bf 1a)} { \it Given   $ s  \in [t,T)  $,
   we first show that along $\wh{\mu}|_{[t,s]} $, the restriction of $\beta$ over $[s,T]$ is still an admissible strategy,
   which will be used in the next step to choose the locally approximately optimal controls, see \eqref{eq:xxd217a}. }

    Let  $ \wt{\mu} \in \cU_s$.
    The process  $ \big(\wh{\mu} \oplus_s \wt{\mu}\big)_r \dfnn \b1_{\{r<s\}} \wh{\mu}_r +  \b1_{\{r  \ge  s\}} \wt{\mu}_r  $, $r \in [t,T]$
  is clearly $\bF-$progressively measurable. Suppose that
   $E \neg \int_t^T  [ \mu_r ]^q_{\overset{}{\hU}} dr + E \neg \int_s^T  [ \wt{\mu}_r ]^{\wt{q}}_{\overset{}{\hU}} dr < \infty$ for some $q > 2$ and $\wt{q} > 2$.
   It follows that
   \beas
     E \neg \int_t^T \neg \big[ \big(   \wh{\mu} \oplus_s \wt{\mu} \big)_r   \big]^{q \land \wt{q}}_{\overset{}{\hU}} dr
     \le E \neg \int_t^T \neg \big[     \mu_r     \big]^{q \land \wt{q}}_{\overset{}{\hU}} dr
       + E \neg \int_s^T \neg \big[      \wt{\mu}_r    \big]^{q \land \wt{q}}_{\overset{}{\hU}} dr < \infty  .
    \eeas
        Thus,  $\wh{\mu} \oplus_s \wt{\mu} \in  \cU_t$. Then we  can  define
    \bea   \label{eq:xux511}
     \beta^s (\wt{\mu})    \dfnn  \big[ \beta (\wh{\mu} \oplus_s \wt{\mu} ) \big]^s \in \cV_s .
    \eea
      For $dr \times dP -$a.s. $(r,\o) \in [s,T] \times \O$,
   \beas
       \big[ (\beta^s (\wt{\mu}))_r (\o) \big]_{\overset{}{\hV}}
    =   \big[ ( \beta (\wh{\mu} \oplus_s \wt{\mu} ) )_r (\o) \big]_{\overset{}{\hV}}
    \le \k + C_\beta  \big[   (\wh{\mu} \oplus_s \wt{\mu} )_r (\o) \big]_{\overset{}{\hU}}
    = \k +  C_\beta  \big[     \wt{\mu}_r (\o) \big]_{\overset{}{\hU}} .
    \eeas
     Let    $\wt{\mu}^1, \wt{\mu}^2 \neg  \in \neg  \cU_s$ such that  $\wt{\mu}^1  \neg = \neg   \wt{\mu}^2 $,
$dr   \times   d P -$a.s. on $\[s,\z\[ \, \cup \, \[\z,T\]_A  $  for some  $  \z   \neg \in  \neg     \cS_{s,T}$ and $ A   \neg \in  \neg  \cF_\z   $.
Then $ \wh{\mu} \oplus_s \wt{\mu}^1 = \wh{\mu} \oplus_s \wt{\mu}^2 $, $dr   \times   d P -$a.s.
 on $\[t,\z\[ \, \cup \, \[\z,T\]_A  $. By Definition \ref{def:strategy},   $\beta (\wh{\mu} \oplus_s \wt{\mu}^1 ) = \beta (\wh{\mu} \oplus_s \wt{\mu}^2 )$, $dr   \times   d P -$a.s.
 on $\[t,\z\[ \, \cup \, \[\z,T\]_A  $.  It follows that for   $dr  \neg \times \neg  d P -$a.s.
  $ (r,\o) \in \[s,\z\[ \, \cup \, \[\z,T\]_A  $
  \beas
   \big( \beta^s (\wt{\mu}^1) \big)_r (\o)
  \neg =  \neg  \big( \beta (\wh{\mu} \oplus_s \wt{\mu}^1 ) \big)_r  (\o)
   \neg    =  \neg  \big( \beta (\wh{\mu} \oplus_s \wt{\mu}^2 ) \big)_r (\o)
    \neg  =   \neg   \big( \beta^s (\wt{\mu}^2 ) \big)_r (\o)    .
    \eeas
      Hence, $\beta^s   \neg \in \neg   \fB_s $.

       \ss \no  {\bf 1b)}
      { \it  Fix $m  \neg \in \neg  \hN$ with $ m \ge C^\f_{x,\d} \dfnn
        \sup \big\{|\f (s,\fx)| : (s,\fx) \in   \ol{O}_{\d + 3}(t,x) \cap ([t,T] \times \hR^k ) \big\}   $.
        According  to the   finite cover $ \{ \fD_m (s_i,x_i) \}^{N_m}_{i=1}$ of $ \ol{O}_{  \d} (t, x) $,
        we 
        use \eqref{eq:xxa255}
        to construct the $1/m-$optimal control $\mu^m$ for player I  under  strategy $\beta$
        by  pasting together local      $1/m-$optimal controls. }

        Given $  i   \neg  = \neg 1,\cds \neg, N_m  $,      \eqref{eq:xxa255} shows that
         there exists  $ \big\{ (A^{m,i}_n, \mu^{m,i}_n) \big\}_{n \in \hN} \neg \subset \neg  \cF_{t_i}  \neg \times \neg  \cU_{t_i}$
     with $ \lmtu{n \to \infty}  \b1_{A^{m,i}_n}  =1    $, \pas ~ such that for any $n \in \hN$
 \bea   \label{eq:xxd217a}
  J \big(t_i,x_i, \mu^{m,i}_n,  \beta^{t_i} (\mu^{m,i}_n)    \big)
  \ge \big( I \big(t_i,x_i,   \beta^{t_i}    \big) -  1/m \big) \land m   ,   \q \pas    \hb{  on }  A^{m,i}_n .
 \eea
 As $Y^{ \wh{\Th} }  \big(T, g \big(X^{ \wh{\Th} }_T\big)\big) \neg \in  \neg  \hC^p_\bF( [t,T] )$,
 the Monotone Convergence Theorem  shows that
 \beas
  \lmtd{n \to \infty} E \bigg[ \b1_{ (A^{m,i}_n)^c  } \Big( \underset{s \in [t,T]}{\sup}
  \Big| Y^{\wh{\Th}}_s \Big(T, g \big( X^{ \wh{\Th} }_T \big) \Big)\Big|^p + \big(C^\f_{x,\d}\big)^p \Big)   \bigg] =0   .
 \eeas
  So there exists an $n(m,i) \in \hN$ such that
  $E \bigg[ \b1_{ \big(A^{m,i}_{n(m,i)}\big)^c  } \Big( \underset{s \in [t,T]}{\sup} \Big| Y^{\wh{\Th}}_s \Big(T, g \big( X^{ \wh{\Th} }_T \big) \Big)\Big|^p + \big(C^\f_{x,\d}\big)^p \Big)   \bigg] \le m^{-(1+p)} N^{-1}_m$.
   Set  $(A^m_i,\mu^m_i) \dfnn \big( A^{m,i}_{n(m,i)}, \mu^{m,i}_{n(m,i)} \big)$
     and    $  \wt{A}^m_i \dfnn  \big\{ \big(\t , X^{\Th }_{\t } \big) \in \fD_m (s_i,x_i)
   \backslash \underset{j < i}{\cup} \fD_m (s_j,x_j) \big\} \in     \cF_{\t}    $.
  As $\wt{A}^m_i \neg \subset \neg  \big\{ \big(\t, X^{\Th}_{\t} \big)  \neg \in \neg  \fD_m (s_i,x_i)     \big\}
   \neg \subset  \neg \{ \t  \neg \le \neg  t_i    \} $, we see that  $ \wt{A}^m_i  \neg = \neg  \wt{A}^m_i
      \cap     \{ \t  \neg \le \neg  t_i    \}  \neg \in \neg  \cF_{t_i}$.

     By the continuity of process $ X^{\Th } $,   $(\t,X^{\Th }_\t) \in \pa O_\d (t,x) $, \pas ~ So
           $\{ \wt{A}^m_i \}^{N_m}_{i =1    }  $ forms a partition of $\cN^c$ for some $P-$null set $\cN$.
       Then  we can define
     an $\bF-$stopping time $\t_m \dfnn \sum^{N_m}_{i=1} \b1_{\wt{A}^m_i} \, t_i +  \b1_{\cN} T \ge  \t  $ as well as a process
     \beas
      \mu^m_s      & \tneg  \dfnn  & \tneg   \b1_{\{ s < \t_m   \}}    \wh{\mu}_s
     \neg + \neg  \b1_{\{ s  \ge  \t_m  \}}  \Big( \sum^{N_m}_{i=1}   \b1_{\wt{A}^m_i \cap A^m_i }
     \big(  \mu^m_i    \big)_s + \b1_{A_m} \wh{\mu}_s  \Big) \\
     & \tneg  =  &  \tneg  \b1_{A_m} \wh{\mu}_s +  \sum^{N_m}_{i=1}   \b1_{\wt{A}^m_i \cap A^m_i }
     \Big( \b1_{\{ s < t_i   \}}  \wh{\mu}_s
    \neg  +  \neg  \b1_{\{ s  \ge  t_i  \}}  \big(  \mu^m_i    \big)_s \Big) , \q   \fa  s  \in [t,T] ,  
    \eeas
  where $A_m \dfnn \Big( \underset{i=1}{\overset{N_m}{\cup}} ( \wt{A}^m_i \backslash A^m_i ) \Big) \cup \cN $.

     Let $s \neg \in \neg  [t,T]$ and $U  \neg \in \neg  \sB \big(\hU  \big)  $.
     As    $ \[t, \t\[ \, \in \neg \sP $,
     we see that   $\cD  \neg \dfnn   \neg    \[t, \t\[
      \, \cap \, ( [t,s] \neg \times \neg \O) \neg \in \neg \sB\big([t,s]\big)  \neg \otimes \neg  \cF_s $.
       The $\bF-$progressive measurability of   $\wh{\mu}$  then  implies that
     \bea   \label{mu_progressive}
    \big\{ (r,  \o) \neg \in \neg  \cD \neg :  \mu^m_r (\o)
      \neg \in \neg  U \big\}
       \neg  =  \neg  \big\{ (r,  \o)  \neg \in \neg  \cD  \neg   :      \wh{\mu}_r (\o)  \neg \in \neg  U \big\}
       \neg = \neg  \cD  \cap  \big\{ (r,  \o)  \neg \in \neg  [t,s]  \neg \times  \neg \O  \neg  :
         \wh{\mu}_r (\o)  \neg \in \neg  U \big\}
        \neg \in \neg  \sB\big([t,s]\big)  \neg \otimes \neg  \cF_s   .   \q
    \eea
Given $i  \neg = \neg  1, \cds \neg , N_m$, we set $\ol{A}^m_i \dfnn (\wt{A}^m_i \backslash  A^m_i) \cup \cN \in \cF_{t_i} $.
 If $s \neg   < \neg   t_i   $, both
$\cD^m_i   \neg   \dfnn  \neg   \[\t,T\]_{\wt{A}^m_i \cap A^m_i}    \cap    ( [t,s]    \times    \O)
      \neg    =   \neg    ( [t_i,T]  \cap  [t,s] )     \times    ( \wt{A}^m_i \cap A^m_i )
$ and $\wh{\cD}^m_i    \neg    \dfnn   \neg    \[\t,T\]_{\ol{A}^m_i  }    \cap    ( [t,s]    \times    \O)
      \neg    =   \neg    ( [t_i,T]  \cap  [t,s] )     \times   \ol{A}^m_i  $  are empty.
 Otherwise, if $ s  \neg  \ge  \neg  t_i  $,   both
         $ \cD^m_i  \neg = \neg    [t_i, s]       \times    ( \wt{A}^m_i \cap A^m_i )    $
        and   $ \wh{\cD}^m_i  \neg = \neg    [t_i, s]    \neg  \times \neg  \ol{A}^m_i
          $ belong to $\sB\big([t_i,s]\big)  \neg \otimes \neg  \cF_s$.
       Using a  similar  argument to \eqref{mu_progressive} on
      the $\bF-$progressive measurability of process  $\mu^m_i$   yields that
               \beas
       \big\{ (r,  \o) \in \cD^m_i \neg :  \mu^m_r (\o) \neg \in \neg U \big\}
      & \tneg \dneg = & \tneg  \dneg  \big\{ (r,  \o) \in \cD^m_i  \neg  :  \big(  \mu^m_i \big)_r (\o)  \neg \in \neg U \big\}
       \neg \in  \neg \sB\big([t_i,s]\big)  \neg \otimes \neg  \cF_s    \neg  \subset   \neg \sB\big([t,s]\big)  \neg \otimes \neg  \cF_s      \\
    \hb{and } \;   \big\{ (r,  \o)   \in    \wh{\cD}^m_i \neg :  \mu^m_r (\o)
      \neg \in \neg  U \big\}      &  \tneg  \dneg  = &  \tneg   \dneg  \big\{ (r,  \o)  \neg \in \neg  \wh{\cD}^m_i  \neg   :      \wh{\mu}_r (\o)  \neg \in \neg  U \big\}
            \neg \in \neg    \sB\big([t,s]\big)  \neg \otimes \neg  \cF_s ,
     \eeas
    both of   which together with \eqref{mu_progressive} shows the $\bF-$progressive measurability of $\mu^m $.
  For $i =1,\cds \neg , N_m$, suppose that  $E  \neg  \int_{t_i}^T  \neg \big[ \big( \mu^m_i \big)_r \big]^{q_i}_{\overset{}{\hU}} \, dr
    \neg < \neg  \infty $ for some $ q_i > 2$. Setting $q_* \dfnn q \land \min\{ q_i: i =1,\cds \neg , N_m \}$, we can deduce that
    \beas   
  E \neg  \int_t^T  \neg \big[ \mu^m_r  \big]^{q_*}_{\overset{}{\hU}} \, dr
   \neg \le \neg   E  \neg \int_t^T     [  \mu_r   ]^{q_*}_{\overset{}{\hU}} \, dr
     +    \sum^{N_m}_{i =1} E \neg \int_{t_i}^T  \neg \big[ \big( \mu^m_i \big)_r \big]^{q_*}_{\overset{}{\hU}} \, dr
    \neg < \neg  \infty .
    \eeas
     Hence,     $\mu^m  \neg \in \neg  \cU_t$.

      \ss \no     {\bf 1c)}
   {\it  Next,    set     $   \Th_m \neg \dfnn \neg   \big( t , x , \mu^m ,   \beta (\mu^m) \big)$.
    We shall use  a series of estimates on state processes $X^{t,\xi,\mu,\nu}$/payoff processes $Y^{t,\xi,\mu,\nu}$,
    a stochastic backward semigroup  property  \eqref{eq:p677}
    as well as   the continuous dependence of $Y^{t,\xi,\mu,\nu}$  on $\xi$
        to demonstrate how  $ J \big(t , x , \mu^m ,   \beta (\mu^m)  \big)$ deviates from
      $ Y^{\Th}_{t}  \big(\t,    \f \big(\t , X^{\Th}_{\t }  \big)   \big) $, which will eventually lead to
         \bea   \label{eq:xqxqx217}
     w_1(t,x)  \ge     \underset{\beta \in \fB_t }{\essinf} \; \underset{\mu \in \cU_t}{\esssup}    \;
    Y^{t,x,\mu, \beta (  \mu ) }_t \Big(\t_{\beta,\mu},
    \f \big(\t_{\beta,\mu},  X^{t,x,\mu, \beta (  \mu ) }_{\t_{\beta,\mu}} \big) \Big), \q \pas
    \eea
    }
  \q  As  $ \mu^m  \neg = \neg \wh{\mu}  \neg = \neg   \mu         $ on $\[t,\t \[$, taking $(\t, A) = (\t, \es)$ in Definition \ref{def:strategy} shows that $  \beta    ( \mu^m)   \neg = \neg    \beta   ( \mu)        $, $ds \times dP -$a.s. on $\[t,\t \[$, and then applying
   \eqref{eq:p611} with $(\t, A) = (\t, \es)$ yields that \pas
   \bea    \label{eq:q101}
    X^{\Th_m  }_s = X^{\Th}_s  \in \ol{O}_\d(x)  , \q \fa s \in [t,\t] .
   \eea
     Thus,  for any $ \eta\in \hL^p (\cF_{\t})  $,      the BSDE$ \big(t, \eta ,  f^{\Th_m}_{\t}    \big)$
      and the BSDE$ \big(t, \eta , f^{\Th}_{\t}    \big)$      are essentially the same. To wit,
   \bea  \label{eq:q104b}
      \big(Y^{\Th_m}  (\t, \eta),Z^{\Th_m}   (\t, \eta)  \big)
      \neg  =  \neg  \big(Y^\Th  \neg   (\t, \eta),Z^\Th  \neg   (\t, \eta)  \big)    . \q
    \eea

  Given $A \in \cF_t$, we see from  \eqref{eq:q101}  that
  \beas
\b1_A  X^{\Th_m }_{  \t_m \land s} & \tneg \neg =&  \neg \tneg \b1_A   X^{\Th_m }_{\t \land s} \neg +
  \neg \b1_A \neg \int_{\t \land s}^{ \t_m \land s}    b \Big(r, X^{\Th_m }_r,\mu^m_r, \big(\beta (\mu^m)\big)_r \Big)  dr
  \neg + \neg \b1_A  \neg \int_{\t \land s}^{ \t_m \land s}
    \si \Big(r, X^{\Th_m }_r,\mu^m_r, \big(\beta ( \mu^m )\big)_r  \Big)  dB_r,    \\
 & \tneg \neg  =& \neg  \tneg \b1_A  X^{\Th }_{\t \land s}  \neg + \neg  \int_{\t \land s}^{ \t_m \land s}
     \b1_A  b \big(r, X^{\Th_m }_{\t_m \land r}, u_0, \big( \beta( \mu^m ) \big)_r \big)     dr
  \neg + \neg \int_{\t \land s}^{ \t_m \land s} \b1_A  \si \big(r, X^{\Th_m }_{\t_m \land r} , u_0,  \big( \beta( \mu^m ) \big)_r \big)  dB_r   , \q  s  \neg \in \neg  [t,T].
 \eeas
 It  follows that
 \bea
 \b1_A \underset{r \in [t,s]}{\sup} \big| X^{\Th_m }_{  \t_m \land r} \neg  - \neg  X^{\Th }_{\t \land r} \big|
 & \tneg \le & \tneg  \int_{\t \land s}^{ \t_m \land s} \neg
  \b1_A  \big|   b \big(r, X^{\Th_m }_{\t_m \land r}, u_0, \big( \beta( \mu^m ) \big)_r \big) \big|    dr  \nonumber \\
 &\tneg  & \tneg   + \neg  \underset{r \in [t,s]}{\sup} \left| \int_{\t \land r}^{ \t_m \land r}  \neg \b1_A  \si \big(r', X^{\Th_m }_{\t_m \land r'},u_0,  \big( \beta( \mu^m ) \big)_{r'} \big)     dB_{r'}   \right|  , \q   s  \neg \in \neg  [t,T].
 \label{eq:xvx103}
 \eea
 Let  $C(\k,x,\d)$   denote a generic constant,  depending on $\k \neg + \neg |x| \neg + \neg \d$, $C^\f_{x,\d}  $,
   $T$, $\g $,  $p$ and $|g(0)|$,        whose form  may vary from line to line.
   Squaring both sides of \eqref{eq:xvx103} and taking expectation, we can deduce from H\"older's inequality, Doob's martingale inequality,
   \eqref{b_si_linear_growth}, \eqref{b_si_Lip}, \eqref{eq:q101} and Fubini's Theorem that
 \bea
\q && \hspace{-1.2cm} E \bigg[ \b1_A  \underset{r \in [t,s]}{\sup} \big| X^{\Th_m }_{  \t_m \land r} \neg  - \neg X^{\Th }_{\t \land r} \big|^2 \bigg] \nonumber \\
  &&   \le  \neg    4 E  \neg  \int_{\t \land s}^{ \t_m \land s}  \neg
  \b1_A  \big|   b \big(r, X^{\Th_m }_{\t_m \land r}, u_0, \big( \beta( \mu^m ) \big)_r \big) \big|^2    dr
  \neg    +  \neg     8 E  \neg   \int_{\t \land s}^{ \t_m \land s}  \neg  \b1_A
  \big|  \si \big(r, X^{\Th_m }_{\t_m \land r}, u_0,  \big( \beta( \mu^m ) \big)_r \big) \big|^2     dr     \nonumber     \\
  &&    \le \neg   12 \g^2    E \int_{\t \land s}^{ \t_m \land s} \b1_A \Big(  \big|   X^{\Th_m }_{\t_m \land r}  \neg - \neg
     X^{\Th }_{\t \land r}  \big|  \neg + \neg   \big|X^{\Th }_{\t \land r}  \big|  \neg + \neg  1
     \neg + \neg  \big[ \big( \beta( \mu^m ) \big)_r \big]_{\overset{}{\hV}} \Big)^2 dr  \nonumber  \\
 && \le  \neg   24 \g^2     \int_t^s    E \Big[ \b1_A \underset{r' \in [t,r]}{\sup} \big|   X^{\Th_m }_{\t_m \land r'}   \neg -   \neg    X^{\Th }_{\t \land r'}  \big|^2 \Big] dr
 \neg + \neg  \frac{C(\k,x,\d)}{m }   P  ( A ) , \q \fa  s  \neg \in \neg  [t,T]  ,   \qq \label{eq:xvx113}
 \eea
    where we used the facts that
 \bea   \label{eq:xvx063}
 \t_m - \t \le \sum^{N_m}_{i=1} \b1_{\wt{A}^m_i} \, 2 \d^m_{s_i,x_i} < \frac{2}{m} ,~ \pas \q \hb{ and } \q
  \big[ \big( \beta( \mu^m ) \big)_r \big]_{\overset{}{\hV}}  \le \k   ,~ dr \times dP-a.s.  \hb{ on } \[\t,\t_m\[ .
 \eea
   Then an application of    Gronwall's inequality yields that
  \beas
  E \bigg[ \b1_A \underset{r \in [t,s]}{\sup} \big| X^{\Th_m }_{  \t_m \land r} -X^{\Th }_{\t \land r} \big|^2 \bigg]
  \le  \frac{C(\k,x,\d)}{m } P  ( A )
     e^{24 \g^2  (s-t)}  , \q \fa  s  \neg \in \neg  [t,T].
  \eeas
  In particular,   $    E \bigg[ \b1_A \underset{r \in [t,T]}{\sup} \big| X^{\Th_m }_{  \t_m \land r}  \neg - \neg X^{\Th }_{\t \land r} \big|^2 \bigg]
   \neg \le \neg  \frac{ C(\k,x,\d)}{m} P(A)$.     Letting $A$ vary in $\cF_t$ yields that
  \bea   \label{eq:xvx043}
  E \bigg[  \underset{r \in [t,T]}{\sup} \big| X^{\Th_m }_{  \t_m \land r}  \neg - \neg X^{\Th }_{\t \land r} \big|^2  \Big| \cF_t \bigg]
  \le  \frac{ C(\k,x,\d)}{m}, \q \pas
  \eea

 Let $i  \neg = \neg  1, \cds, N_m $ and set  $\Th^{t_i}_m  \neg \dfnn \neg \big(t_i, X^{\Th_m}_{t_i}, [\mu^m]^{t_i},
 [\beta(\mu^m)]^{t_i}  \big) $. We see from \eqref{eq:xxa603}  that
 $X^{\Th_m}_T = X^{\Th^{t_i}_m}_T$, \pas ~  It then follows from   \eqref{eq:xxa605} that
       \bea \label{eq:xxd213}
        Y^{\Th_m}_{t_i} \big(T, g \big(X^{\Th_m}_T\big)\big)
=   Y^{\Th^{t_i}_m}_{t_i} \big(T, g \big(X^{\Th_m}_T\big)\big)
=   Y^{\Th^{t_i}_m}_{t_i} \big(T, g \big(X^{\Th^{t_i}_m}_T\big)\big)
=   J( \Th^{t_i}_m ), \q \pas
    \eea
  Similar to $\mu^m$,
 \beas
  \big(\wh{\mu}^m_i\big)_s
    &\tneg \dfnn & \tneg  \b1_{\{ s < \t_m   \}}    \wh{\mu}_s
     \neg + \neg  \b1_{\{ s  \ge  \t_m  \}}  \Big(\b1_{\wt{A}^m_i \cap A^m_i}
     \big(  \mu^m_i    \big)_s + \b1_{(\wt{A}^m_i \cap A^m_i)^c}
        \wh{\mu}_s \Big)     \\
        &\tneg = & \tneg  \b1_{ \wt{A}^m_i \cap A^m_i } \Big( \b1_{\{ s < t_i  \}} \wh{\mu}_s
  \neg + \neg  \b1_{\{ s  \ge  t_i  \}}
    \big( \mu^m_i \big)_s    \Big) \neg + \neg  \b1_{( \wt{A}^m_i \cap A^m_i )^c}  \wh{\mu}_s     ,  \q  s   \neg \in \neg  [t,T]
    \eeas
      also defines a   $\cU_t-$process.
As $\mu^m  \neg = \neg  \wh{\mu}^m_i $ on $ \[t,\t_m\[ \, \cup \, \[\t_m,T\]_{\wt{A}^m_i  \cap A^m_i }$
and $ \wh{\mu}^m_i  \neg = \neg  \wh{\mu}    \oplus_{t_i}    \mu^m_i  $ on $ \big( [t, t_i  )  \neg \times \neg  \O \big)  \cup  \big( [t_i, T]  \neg \times \neg  (\wt{A}^m_i \cap A^m_i ) \big) $,
 Definition \ref{def:strategy} shows that  $\beta( \mu^m )  \neg  = \neg \beta \big( \wh{\mu}^m_i \big)$, $ds \times dP-$a.s. on $ \[t,\t_m\[  \cup  \[\t_m,T\]_{\wt{A}^m_i \cap A^m_i  }$ and $ \beta \big( \wh{\mu}^m_i \big)  \neg = \neg  \beta \big(  \wh{\mu} \oplus_{t_i} \mu^m_i \big)  $, $ds \times dP-$a.s.  on $ \big( [t, t_i  )  \neg \times \neg  \O \big)  \cup  \big( [t_i, T] \neg \times \neg ( \wt{A}^m_i \cap A^m_i ) \big) $. Thus $\big( \mu^m, \beta( \mu^m ) \big)  \neg = \neg   \Big( \wh{\mu} \oplus_{t_i} \mu^m_i,   \beta \big(  \wh{\mu} \oplus_{t_i} \mu^m_i \big) \Big) $, $ds \times dP-$a.s.  on $   \[\t_m,T\]_{\wt{A}^m_i \cap A^m_i}  \neg = \neg    [t_i, T]  \neg \times \neg  (\wt{A}^m_i \cap A^m_i)  $.
From \eqref{eq:xux511}, one has   $\big( [\mu^m]^{t_i}, [\beta( \mu^m )]^{t_i} \big)
   \neg = \neg   \big(    \mu^m_i ,     \beta^{t_i}  (    \mu^m_i  )  \big)
     $, $ds \times dP-$a.s.  on $       [t_i, T]  \neg \times \neg ( \wt{A}^m_i \cap A^m_i)  $.
      Then  by \eqref{eq:xxd213},      \eqref{eq:xxa273} and      \eqref{eq:s025},
      it holds \pas ~ on $   \wt{A}^m_i \cap A^m_i \in \cF_{t_i}  $ that
   $$
   Y^{ \Th_m}_{\t_m} \big(T, g \big(X^{ \Th_m}_T\big)\big)
    \neg   =  \neg    Y^{ \Th_m}_{t_i} \big(T, g \big(X^{ \Th_m}_T\big)\big)
   \neg  =  \neg   J \Big(  t_i, X^{ \Th_m}_{\t_m},  \mu^m_i, \beta^{t_i}  \big( \mu^m_i \big)  \Big)
      \neg    \ge    \neg    J \Big(  t_i, X^{\Th}_{\t}, \mu^m_i, \beta^{t_i}  \big( \mu^m_i \big)  \Big)
     \neg  -  \neg  c_0 | X^{ \Th_m}_{\t_m}  \neg - \neg  X^{\Th}_{\t} |^{2/p}      .
   $$
  Since $  \fD_m (s_i,x_i)  \cap \ol{O}_\d (t, x) \ne \es $, it is easy to see that
 \beas
      \ol{\fD}_m (s_i,x_i) = \big[  s_i - \d^m_{s_i,x_i} , s_i + \d^m_{s_i,x_i}  \big]
     \times \ol{O}_{  \d^m_{s_i,x_i} } (x_i) \subset \ol{O}_{\d+ 2 \sqrt{2} \d^m_{s_i,x_i} } (t,x)
     \subset \ol{O}_{\d+\frac{2 \sqrt{2}}{m} } (t,x)
     \subset \ol{O}_{\d+3 } (t,x)  .
 \eeas
  So $\f(t_i,x_i) \le C^\f_{x,\d}  <  m + 1/m  $. On the other hand, one has
 $\f(t_i,x_i) \le w_1(t_i,x_i) \le  I \big(t_i,x_i,   \beta^{t_i}    \big)$, \pas  ~  Then
   it follows from \eqref{eq:xxd217a} that
  \beas  
 \f(t_i,x_i)   \le  I \big(t_i,x_i,   \beta^{t_i}    \big) \land  ( m + 1/m    )
 \le J \big(t_i,x_i, \mu^m_i,  \beta^{t_i} (\mu^m_i)    \big) + 1/m ,   \q \pas    \hb{  on }  A^m_i .
 \eeas
     As $\dis \big| X^{\Th}_{\t} - x_i \big|^{2/p}
 <  ( \d^m_{s_i,x_i})^{2/p} < m^{-2/p} \le 1/m $ on $\wt{A}^m_i$,
    we can also deduce from
 \eqref{eq:s025},     \eqref{eq:xux391} and the continuity of $\f$
   that   it holds \pas ~ on $   \wt{A}^m_i \cap A^m_i    $  that
   $$
       J \Big(  t_i, X^{\Th}_{\t}, \mu^m_i, \beta^{t_i}  \big( \mu^m_i \big)  \Big)
               \neg    \ge    \neg      J \big( t_i,x_i , \mu^m_i, \beta^{t_i}  \big( \mu^m_i \big)  \big)
 - \frac{c_0}{m}
     \neg   \ge    \neg         \f \big(t_i, x_i  \big)    -    \frac{c_0}{m}
   \neg    \ge    \neg     \f \big(s_i, x_i  \big)   -     \frac{c_0}{m}
        \ge      \f \big(\t, X^{ \Th}_{\t}  \big)   -     \frac{c_0}{m} \neg \dfnn \neg  \eta_m
        \neg  \in   \neg  \hL^\infty ( \cF_\t )  .
   $$
 Thus it holds \pas ~ on $   \cup^{N_m}_{i=1} (\wt{A}^m_i \cap A^m_i)    $  that
  \bea
  Y^{ \Th_m}_{\t_m} \big(T, g \big(X^{ \Th_m}_T\big)\big)  \ge
  \eta_m  \neg  -  \neg  c_0 | X^{ \Th_m}_{\t_m}  \neg - \neg  X^{\Th}_{\t} |^{2/p}
  \dfnn   \wt{\eta}_m \in \hL^p ( \cF_{\t_m} ) .     \label{eq:xux671}
  \eea
         By \eqref{eq:xvx071},     it holds  \pas ~ that
   \bea   \label{eq:xux683}
    \big|  Y^{\Th}_{t}  (\t,  \eta_m  )
   \neg  - \neg  Y^{\Th}_{t}   (\t,    \f \big(\t , X^{\Th}_{\t }   )   \big)  \big|^p
       \neg \le  \neg c_0   E \Big[ \big|   \eta_m      \neg  -  \neg
        \f \big(\t , X^{\Th}_{\t }  \big)     \big|^p   \Big|\cF_t\Big]    \le
          \frac{c_0}{m^p}           .
   \eea

  Let $(Y^m, Z^m) \in \hG^p_\bF([t,T])$ be the unique solution of the following BSDE with zero generator:
   \beas
   Y^m_s = Y^{\Th_m}_{\t}  (\t_m, \eta_m )  - \int_s^T Z^m_r d B_r      , \q    s \in [t,T] .
   \eeas
  For any  $s \in [t,T]$, one can deduce that
   \beas
      Y^m_{\t \land s} = E[Y^m_{\t \land s} | \cF_\t] = E \bigg[ Y^{\Th_m}_{\t}  (\t_m, \eta_m )  - \int_{\t \land s}^T Z^m_r d B_r \Big| \cF_\t \bigg]
   = Y^{\Th_m}_{\t}  (\t_m, \eta_m )  - \int_{\t \land s}^\t Z^m_r d B_r , \q  \pas ~
   \eeas
       By the continuity of process $Y^m$,  it holds \pas ~ that
    \bea   \label{eq:xvx031}
   Y^m_{\t \land s}  = Y^{\Th_m}_{\t}  (\t_m, \eta_m )  - \int_{\t \land s}^\t Z^m_r d B_r  = Y^{\Th_m}_{\t}  (\t_m, \eta_m )  - \int_s^T \b1_{\{r < \t\}} Z^m_r d B_r     , \q    s \in [t,T] .
   \eea
   Thus, we see that  $(Y^m_s, Z^m_s) \neg = \neg  \big(  Y^m_{\t \land s}, \b1_{\{s < \t\}} Z^m_s \big),
   s  \neg \in \neg  [t,T]   $. Also,
   taking $[\cd|\cF_{\t \land s}]$ in \eqref{eq:xvx031} shows that \pas
   \beas
    Y^m_s  \neg = \neg   Y^m_{\t \land s}  \neg = \neg  E\big[Y^{\Th_m}_{\t}  (\t_m,  \eta_m ) \big|\cF_{\t \land s}\big] ,
    \q   \fa  s  \neg \in \neg  [t,T] .
   \eeas
    On the other hand,
   let $(\wt{Y}^m, \wt{Z}^m ) \in \hG^p_\bF([t,T])$ be the unique solution of the following BSDE with zero generator:
   \bea  \label{eq:xvx036}
   \wt{Y}^m_s =    \eta_m     - \int_s^T \wt{Z}^m_r d B_r      , \q    s \in [t,T] .
   \eea
   Similar to $(Y^m, Z^m) $, it holds \pas ~ that
     \bea   \label{eq:xvx061}
     \big(\wt{Y}^m_s, \wt{Z}^m_s\big) \neg = \neg \big(  \wt{Y}^m_{\t \land s}, \b1_{\{s < \t\}} \wt{Z}^m_s \big)
       \q \hb{ and } \q      \wt{Y}^m_s   \neg  = \neg  E  [   \eta_m  |\cF_{\t \land s}  ] ,
       \q        \fa s \neg \in \neg  [t,T]          .
 \eea

 We   can  deduce that  $( \cY^m,  \cZ^m )
   \dfnn  \big\{ \big( \b1_{\{s < \t\}} Y^m_s \neg + \neg  \b1_{\{s  \ge  \t\}} Y^{\Th_m}_s  (\t_m, \eta_m ),
    \b1_{\{s < \t\}} Z^m_s  \neg + \neg  \b1_{\{s  \ge  \t\}} Z^{\Th_m}_s  (\t_m, \eta_m )     \big) \big\}_{s \in [t,T]} \in \hG^p_\bF([t,T]) $
   solves the following BSDE
   \bea
   \cY^m_s & \tneg =& \tneg  \b1_{\{s  \ge  \t\}} Y^{\Th_m}_s  (\t_m, \eta_m )
   \neg + \neg   \b1_{\{s < \t\}}Y^{\Th_m}_{\t}  (\t_m, \eta_m )   \neg - \neg  \b1_{\{s < \t\}} \int_s^T Z^m_r d B_r
   \neg  = \neg  Y^{\Th_m}_{\t \vee s}  (\t_m, \eta_m )  \neg - \neg \b1_{\{s < \t\}} \int_s^T \b1_{\{r < \t\}} Z^m_r    d B_r  \nonumber \\
   & \tneg =& \tneg  \eta_m  \neg + \neg  \int_{\t \vee s}^T f^{\Th_m}_{\t_m} \big(r,Y^{\Th_m}_r  (\t_m, \eta_m ), Z^{\Th_m}_r  (\t_m, \eta_m ) \big) dr
     \neg - \neg  \int_{\t \vee s}^T Z^{\Th_m}_r  (\t_m, \eta_m ) dB_r  \neg - \neg  \int_s^T \b1_{\{r < \t\}} Z^m_r    d B_r \nonumber   \\
   & \tneg =& \tneg  \eta_m  \neg + \neg  \int_s^T \b1_{\{r  \ge  \t \}} f^{\Th_m}_{\t_m} \big(r,\cY^m_r , \cZ^m_r \big) dr
     \neg - \neg  \int_s^T   \cZ^m_r   dB_r     , \q s \in [t,T] .   \label{eq:xvx034}
   \eea
   Since \eqref{f_Lip}, H\"older's inequality and \eqref{eq:s031} imply that
\beas
E \bigg[       \int_t^T \neg \b1_{\{s  \ge  \t \}} \big| f^{\Th_m}_{\t_m} \big(s,  \wt{Y}^m_s ,
    \wt{Z}^m_s   \big)  \big|^p ds       \bigg]
   \neg \le \neg  c_p  E \bigg[       \int_t^T \neg  \big| f^{\Th_m}_{\t_m} \big(s,  0 ,
    0   \big)  \big|^p ds
    \neg  + \neg     \underset{s \in [t,T]}{\sup}  \big|  \wt{Y}^m_s \big|^p
    \neg  + \neg     \Big(     \int_t^T \neg  \big|   \wt{Z}^m_s \big|^2       ds \Big)^{p/2}      \bigg] < \infty ,
\eeas
  applying  \eqref{eq:n211} to $  \cY^m - \wt{Y}^m $ and using \eqref{eq:xvx061} yield that
   \bea
     E \Big[      \big|   Y^{\Th_m}_\t  (\t_m, \eta_m ) \neg   - \neg   \eta_m     \big|^p  \Big| \cF_t  \Big]
    & \tneg  \dneg  = & \tneg  \dneg     E \Big[      \big|     \cY^m_\t  \neg  - \neg   \wt{Y}^m_\t     \big|^p   \Big|\cF_t \Big]
       \neg    \le  \neg    E \Big[ \, \underset{s \in [t,T]}{\sup}  \big|    \cY^m_s  \neg  - \neg   \wt{Y}^m_s     \big|^p  \Big|\cF_t  \Big]
     \neg   \le  \neg    c_0 E \bigg[       \int_\t^T \neg  \big| f^{\Th_m}_{\t_m} \big(s,  \wt{Y}^m_s ,
    \wt{Z}^m_s   \big)  \big|^p ds   \Big|\cF_t   \bigg]   \nonumber \\
    & \tneg  \dneg   = & \tneg  \dneg   c_0 E \bigg[      \int_\t^{\t_m} \neg  \big| f  \big(s, X^{\Th_m}_{\t_m \land s}, \eta_m ,
    0 ,u_0, (\beta(\mu^m))_s  \big)  \big|^p ds  \Big| \cF_t  \bigg] , \q \pas      \label{eq:xvx114}
     \eea
     Then one can deduce from  \eqref{eq:xvx071},  \eqref{f_linear_growth},
     \eqref{f_Lip}, \eqref{eq:q101}, \eqref{eq:xvx063}   and \eqref{eq:xvx043}    that
     \bea
    && \hspace{-1cm} \big| Y^{\Th_m}_{t} \big(\t ,        Y^{\Th_m}_{\t} \big(\t_m,         \eta_m     \big)     \big)
    -  Y^{\Th_m}_{t} \big(\t,         \eta_m     \big)   \big|^p
    \le c_0 E \Big[      \big|   Y^{\Th_m}_\t  (\t_m, \eta_m )   -\eta_m     \big|^p  \Big| \cF_t  \Big]  \nonumber  \\
   &&     \le c_0 E \bigg[      \int_\t^{\t_m} \neg \Big(1+ \big|X^{\Th_m}_{\t_m \land s}
    - X^\Th_{\t \land s} \big|^2 + | X^\Th_{\t \land s}  |^2
    +  | \eta_m  |^p +  \big[ (\beta(\mu^m))_s \big]^2_{\overset{}{\hV}}  \Big) ds  \Big| \cF_t  \bigg] \nonumber \\
  &&  \le  c_0  E \bigg[        (\t_m - \t) \cd  \underset{s \in [t,T]}{\sup} \big|X^{\Th_m}_{\t_m \land s} - X^\Th_{\t \land s} \big|^2     \Big| \cF_t  \bigg]
    + \frac{c_0}{m}   \Big\{ 1   + (|x|+\d)^2
    +  \Big(C^\f_{x,\d} + \frac{c_0}{m}   \Big)^p +  \k^2  \Big\} \nonumber \\
  &&  \le \frac{C(\k,x,\d)}{m^2}+   \frac{C(\k,x,\d)}{m}
        +   \frac{c_0}{m^{p+1}}   \le   \frac{C(\k,x,\d)}{m}   , \q \pas   \label{eq:xvx117}
     \eea

    Applying  \eqref{eq:p677} with $(\z,\t, \eta) = (\t,\t_m, \eta_m)$, applying
     \eqref{eq:q104b} with $\eta =\eta_m$ and using \eqref{eq:xux683} yield that \pas
    \bea
 Y^{\Th_m}_{t} \big(\t_m,         \eta_m     \big)
 & \tneg =& \tneg  Y^{\Th_m}_{t} \big(\t ,        Y^{\Th_m}_{\t} \big(\t_m,         \eta_m     \big)     \big)
  \ge   Y^{\Th_m}_{t} \big(\t,         \eta_m     \big)  -    \frac{C(\k,x,\d)}{m^{1/p}} \nonumber \\
 & \tneg =& \tneg   Y^{\Th}_{t} \big(\t,         \eta_m     \big)-    \frac{C(\k,x,\d)}{m^{1/p}}
  \ge  Y^{\Th}_{t}  \big(\t,    \f \big(\t , X^{\Th}_{\t }  \big)   \big)
   -    \frac{C(\k,x,\d)}{m^{1/p}}.   \label{eq:h041}
    \eea
        As  $ \mu^m  \neg = \neg \wh{\mu}            $ on $\[t,\t_m \[$,
    taking $(\t, A) = (\t_m, \es)$ in Definition \ref{def:strategy} shows that
    $  \beta    ( \mu^m)   \neg = \neg    \beta   ( \wh{\mu} )        $, $ds \times dP -$a.s.
    on $\[t,\t_m \[$, and then applying
   \eqref{eq:p611} with $(\t, A) = (\t_m, \es)$ yields that \pas
   \bea    \label{eq:q101b}
    X^{\Th_m  }_s = X^{\wh{\Th}}_s     , \q \fa s \in [t,\t_m] .
   \eea
     Given $i = 1, \cds, N_m$, \eqref{eq:q101b} shows that $ X^{\Th_m  }_{t_i} = X^{\wh{\Th}}_{t_i} $, \pas ~ on $\wt{A}^m_i \backslash A^m_i$.
   As $\mu^m  \neg = \neg  \wh{\mu}  $ on $ \[t,\t_m\[ \, \cup \, \[\t_m,T\]_{\wt{A}^m_i  \backslash A^m_i }$,
 Definition \ref{def:strategy} shows that  $\beta( \mu^m )  \neg  = \neg \beta \big( \wh{\mu}  \big)$, $ds \times dP-$a.s. on $ \[t,\t_m\[ \, \cup \, \[\t_m,T\]_{\wt{A}^m_i \backslash A^m_i  }$.
  So     $\big( [\mu^m]^{t_i}, [\beta( \mu^m )]^{t_i} \big)
   \neg = \neg   \big(    [\wh{\mu}]^{t_i}, [\beta( \wh{\mu} )]^{t_i}  \big)
     $ holds  $ds \times dP-$a.s.  on $       \[\t_m,T\]_{\wt{A}^m_i \backslash A^m_i}  \neg = \neg    [t_i, T]  \neg \times \neg  (\wt{A}^m_i \backslash A^m_i)  $. Then by   \eqref{eq:xxa273} and a similar argument to \eqref{eq:xxd213},
 it holds \pas ~ on $\wt{A}^m_i \backslash A^m_i$ that
   \bea   \label{eq:xvx073}
   Y^{\Th_m}_{{\t_m}} \big(T, g \big( X^{ \Th_m }_T \big) \big)
   = Y^{\Th_m}_{t_i} \big(T, g \big( X^{ \Th_m }_T \big) \big)
   = J \big(  \Th^{t_i}_m \big)
    =       J \big( \wh{\Th}^{t_i} \big)
     =       Y^{\wh{\Th}}_{t_i} \big(T, g \big( X^{ \wh{\Th} }_T \big) \big)
   =       Y^{\wh{\Th}}_{{\t_m}} \big(T, g \big( X^{ \wh{\Th} }_T \big) \big) ,
   \eea
    where $ \wh{\Th}^{t_i}    \neg \dfnn \neg \big(t_i, X^{\wh{\Th}}_{t_i}, [\wh{\mu}]^{t_i},
 [\beta(\wh{\mu})]^{t_i}  \big) $.

  \ss  Let $ \wh{\eta}_m \dfnn Y^{\Th_m}_{{\t_m}} \big(T, g \big( X^{ \Th_m }_T \big) \big) \land \wt{\eta}_m
    \in \hL^p \big( \cF_{\t_m} \big)  $ and set $ \wt{A}_m \dfnn \big\{ Y^{ \Th_m}_{\t_m} \big(T, g \big(X^{ \Th_m}_T\big)\big)  <
      \wt{\eta}_m \big\} \in \cF_{\t_m} $. Clearly,  $\b1_{\wt{A}_m} \le    \b1_{A_m} $, \pas ~
   Applying \eqref{eq:xvx071} again, we can deduce from \eqref{eq:xvx043} and \eqref{eq:xvx073} that  \pas
    \bea
             \big|  Y^{\Th_m}_{t}  ({\t_m},  \wh{\eta}_m  )
   \neg  - \neg  Y^{\Th_m}_{t}  \big({\t_m},      \eta_m   \big)  \big|^p
      & \tneg \dneg \le  & \tneg  \dneg     c_0    E \Big[ \big|   \wh{\eta}_m      \neg  -  \neg
                \eta_m      \big|^p   \Big|\cF_t\Big]
                    =     c_0    E \Big[ \b1_{ \wt{A}^c_m   } \big|  \wt{\eta}_m       \neg  -  \neg
                \eta_m     \big|^p + \b1_{ \wt{A}_m  } \big|   Y^{\Th_m}_{{\t_m}} \big(T, g \big( X^{ \Th_m }_T \big) \big)      \neg  -  \neg
                \eta_m     \big|^p   \Big|\cF_t\Big]
        \nonumber    \\
         & \tneg  \dneg  \le  & \tneg  \dneg    c_0    E \Big[   \big|  X^{\Th_m}_{\t_m}       \dneg  -  \dneg
                 X^{\Th}_{\t}      \big|^2  \neg + \neg  \b1_{ A_m  } \big| Y^{\Th_m}_{{\t_m}} \big(T, g \big( X^{ \Th_m }_T \big) \big)         \neg  -  \neg   \eta_m     \big|^p   \Big|\cF_t\Big]   \nonumber    \\
         & \tneg  \dneg     \le  & \tneg  \neg   \neg  \frac{C(\k,x,\d)}{m} \neg + \neg   c_0    E \Big[      \b1_{ A_m  } \big|   Y^{\wh{\Th}}_{{\t_m}} \big(T, g \big( X^{ \wh{\Th} }_T \big) \big)    \neg  -  \neg   \f(\t, X^\Th_\t)        \big|^p   \Big|\cF_t\Big]  \neg + \neg \frac{c_0}{m^p}         \nonumber    \\
    & \tneg    \dneg   \le   & \tneg  \dneg  \frac{C(\k,x,\d)}{m}
  + c_0  E \Big[ \b1_{ A_m  } \Big( \underset{s \in [t,T]}{\sup} \Big| Y^{\wh{\Th}}_s \big(T, g \big( X^{ \wh{\Th} }_T \big) \big)\Big|^p + \big(C^\f_{x,\d}\big)^p \Big)  \Big|\cF_t\Big]     .
    \label{eq:xux711}
        \eea
   Applying  \eqref{eq:p677} with $(\z,\t,\eta) = \big(\t_m, T,g \big( X^{ \Th_m }_T \big) \big)$,
     we see from Proposition \ref{prop_BSDE_estimate_comparison} (2), \eqref{eq:xux711} and \eqref{eq:h041}   that \pas
    \bea
  && \hspace{-1cm}    Y^{\Th_m}_{t} \Big(T,g \big( X^{ \Th_m }_T \big) \Big)
    =     Y^{\Th_m}_{t} \Big(\t_m,
    Y^{\Th_m}_{\t_m} \big(T, g \big( X^{ \Th_m }_T \big) \big)   \Big)
    \neg  \ge  \neg    Y^{\Th_m}_{t}  \big( \t_m,         \wh{\eta}_m      \big)
     \nonumber \\
  &&    \ge     Y^{\Th}_{t}  \big(\t,    \f \big(\t , X^{\Th}_{\t }  \big)   \big)
   -    \frac{ C(\k,x,\d) }{m^{1/p}} - c_0  \bigg\{  E \Big[ \b1_{ A_m  } \Big( \underset{s \in [t,T]}{\sup} \Big| Y^{\wh{\Th}}_s \big(T, g \big( X^{ \wh{\Th} }_T \big) \big)\Big|^p + \big(C^\f_{x,\d}\big)^p \Big)  \Big|\cF_t\Big]  \bigg\}^{\frac{1}{p}}     .
              \label{eq:xux623}
  \eea

 Letting $      \wh{A}_m \dfnn \Big\{ E \Big[ \b1_{ A_m  } \Big( \underset{s \in [t,T]}{\sup} \Big| Y^{\wh{\Th}}_s \big(T, g \big( X^{ \wh{\Th} }_T \big) \big)\Big|^p + \big(C^\f_{x,\d}\big)^p \Big)  \Big|\cF_t\Big]  > 1/ m   \Big\} $,
       one can deduce  that
        \beas
 P( \wh{A}_m  )
 & \tneg \le & \tneg  m E \bigg[  E \Big[ \b1_{ A_m  } \Big( \underset{s \in [t,T]}{\sup} \Big| Y^{\wh{\Th}}_s \big(T, g \big( X^{ \wh{\Th} }_T \big) \big)\Big|^p + \big(C^\f_{x,\d}\big)^p \Big)  \Big|\cF_t\Big] \bigg] \\
 & \tneg  \le & \tneg  \sum^{N_m}_{i=1} m E \Big[ \b1_{  ( A^m_i )^c } \Big( \underset{s \in [t,T]}{\sup} \Big| Y^{\wh{\Th}}_s \big(T, g \big( X^{ \wh{\Th} }_T \big) \big)\Big|^p + \big(C^\f_{x,\d}\big)^p \Big)  \Big] \le m^{-p} .
 \eeas
     Multiplying $ \b1_{\wh{A}_m^c} $ to both sides of \eqref{eq:xux623}      yields that
    \bea   \label{eq:xxd231}
      \b1_{\wh{A}_m^c}  I \big(t,x,    \beta  \big)    \ge   \b1_{\wh{A}_m^c} J \big(t , x , \mu^m ,   \beta (\mu^m)  \big)
      \neg    \ge   \neg \b1_{\wh{A}_m^c} Y^{\Th}_{t}  \big(\t,    \f \big(\t , X^{\Th}_{\t }  \big)   \big)
        \neg - \neg  \frac{C(\k,x,\d)}{m^{1/p}}, \q \pas
\eea
    As  $ \dis   \sum_{m \in \hN} P \big( \wh{A}_m \big) \le \sum_{m \in \hN} m^{-p} < \infty $, Borel-Cantelli theorem shows that $P\big( \lsup{m \to \infty} \b1_{\wh{A}_m} = 1 \big) = 0$. It follows that
    $P\big( \lsup{m \to \infty} \b1_{\wh{A}_m} = 0 \big) = 1 $ and thus
    \bea \label{eq:xux523}
      \lmt{m \to \infty} \b1_{\wh{A}_m} = 0 ,  \q  \pas
      \eea
    So letting $m \to \infty$ in \eqref{eq:xxd231} yields that
   $  I \big(t,x,    \beta  \big)  \ge     Y^{t,x,\mu,\beta(\mu)}_{t}  \big(\t_{\beta,\mu},    \f \big(\t_{\beta,\mu} , X^{t,x,\mu,\beta(\mu)}_{\t_{\beta,\mu} }  \big)   \big)    $,   \pas ~
       Taking essential supremum over $ \mu \in \cU_t $  and then  taking essential infimum over $\beta \in \fB_{t}$, we obtain \eqref{eq:xqxqx217}.

 \no {\bf 1d)} {\it Now let us  show the other inequality of Theorem \ref{thm_DPP} \(1\).
  Similar to $\mu^m$,         we shall first   use \eqref{eq:xxa257}
        to construct the $1/m-$optimal strategy $\beta_m$ by pasting
        together local      $1/m-$optimal strategies with respect to
        the   finite cover $ \{ \fD_m (s_i,x_i) \}^{N_m}_{i=1}$ of $ \ol{O}_{  \d} (t, x) $.  }

       Fix $m \in \hN$.      For  $  i   =1,\cds \neg, N_m  $,
              \eqref{eq:xxa257} shows that    there exists  $(\cA^m_i, \beta^m_i) \in \cF_{t_i} \times \fB_{t_i}$
     with $  P \big(  \cA^m_i  \big)  \ge   1- m^{\frac{1+p^2}{1-p}} N^{-1}_m   $ such that
 \bea   \label{eq:xxd083}
 \wt{\f}(t_i,x_i)  \ge  w_1(t_i,x_i)  \ge   I \big(t_i,x_i,   \beta^m_i   \big) - 1/m ,   \q \pas    \hb{  on }  \cA^m_i .
 \eea

      Let $\beta_\psi$ be the $\fB_t -$strategy  considered in \eqref{def_beta_psi} and fix $\beta \neg \in \neg \fB_t$.
      For any $\mu \neg \in \neg \cU_t$,  we  simply denote $\t_{\beta,\mu}$ by $\t_\mu$ and  define
     \beas
     \big(  \wh{\beta} (\mu) \big)_s        \dfnn     \b1_{\{ s < \t_\mu  \}} \big( \beta (\mu) \big)_s
     \neg + \neg  \b1_{\{ s \ge \t_\mu    \}} \big( \beta_\psi (\mu) \big)_s
       , \q \fa  s  \in [t,T] ,
    \eeas
    which is a $\cV_t-$control by Lemma \ref{lem_control_combine}.
      By   (A-u),  it holds        $ds \times dP-$a.s.  that
         \bea   \label{eq:xux717a}
        \big[ ( \wh{\beta} (\mu))_s \big]_{\overset{}{\hV}}
        = \b1_{\{ s < \t_\mu  \}} \big[\big( \beta (\mu) \big)_s\big]_{\overset{}{\hV}}
     \neg + \neg  \b1_{\{   s  \ge \t_\mu  \}} \big[ \big( \beta_\psi (\mu) \big)_s \big]_{\overset{}{\hV}}
     \le \k + ( C_\beta \vee \k)  [\mu_s ]_{\overset{}{\hU}} .     \q
         \eea

   To see $\wh{\beta} \in \fB_t$, we let    $\mu^1, \mu^2 \neg  \in \neg  \cU_t $ such that  $\mu^1  \neg = \neg   \mu^2 $,
$ds  \neg \times \neg  d P -$a.s. on $ \[t,\t\[ \, \cup \, \[\t, T\]_A   $  for some  $ \t \neg \in  \neg   \cS_{t,T}$
  and $ A  \neg  \in  \neg  \cF_\t  $.
 Since $ \beta(\mu^1) \neg  = \neg  \beta(\mu^2) $, $ds  \neg \times \neg  dP-$a.s. on $ \[t,\t\[ \, \cup \, \[\t, T\]_A   $
 by Definition \ref{def:strategy},
  it holds    $ds  \neg \times \neg  dP-$a.s. on $       \big( \[t,\t\[ \, \cup \, \[t, T\]_A \big) \cap
  \[ t, \t_{\mu^1} \land \t_{\mu^2} \[ $ that
 \bea \label{eq:xxd075a}
 \big(\wh{\beta}(\mu^1) \big)_s  = \big( \beta (\mu^1) \big)_s  = \big( \beta (\mu^2) \big)_s  = \big(\wh{\beta}(\mu^2) \big)_s \,.
 \eea
   And  \eqref{eq:p611} shows that   except on a $P-$null set $ \cN $
 \bea  \label{eq:xxd071a}
  \b1_A X^{\Th_{\mu^1}}_s   + \b1_{A^c} X^{\Th_{\mu^1}}_{\t \land s}
  = \b1_A  X^{\Th_{\mu^2}}_s   + \b1_{A^c} X^{\Th_{\mu^2}}_{\t \land s}  , \q   \fa   s \in [t, T ] .
 \eea
   Then it holds
 for any $ \o \in A  \cap \cN^c  $ that
  \beas
   \t_{\mu^1} (\o)
 =   \inf \Big\{s  \neg \in \neg  ( t , T  ] \neg  : \big( s,X^{\Th_{\mu^1}}_s (\o) \big)
  \neg  \notin   \neg   O_\d (t, x) \Big\}
   =  \inf \Big\{s  \neg \in \neg  ( t , T  ] \neg  : \big( s,X^{\Th_{\mu^2}}_s (\o) \big)
  \neg  \notin   \neg   O_\d (t, x) \Big\} = \t_{\mu^2}  (\o)      .
 \eeas
 Let $A_o \dfnn    \{ \t   \ge    {\t_{\mu^1} \land \t_{\mu^2}}  \}   $. We can deduce from  \eqref{eq:xxd071a} that
 for any $ \o \in A_o \cap \{\t_{\mu^1} \le \t_{\mu^2} \} \cap \cN^c  $
    \beas
   \t_{\mu^1} (\o) & \tneg =& \tneg   \inf \Big\{s  \neg \in \neg  ( t , T  ] \neg  : \big( s,X^{\Th_{\mu^1}}_s (\o) \big)
 \neg  \notin   \neg   O_\d (t, x) \Big\}
 =   \inf \Big\{s  \neg \in \neg  ( t , \t(\o)  ] \neg  : \big( s,X^{\Th_{\mu^1}}_s (\o) \big)
  \neg  \notin   \neg   O_\d (t, x) \Big\}  \nonumber  \\
& \tneg =&  \tneg  \inf \Big\{s  \neg \in \neg  ( t , \t(\o)  ] \neg  : \big( s,X^{\Th_{\mu^2}}_s (\o) \big)
  \neg  \notin   \neg   O_\d (t, x) \Big\}  \ge   \inf \Big\{s  \neg \in \neg  ( t , T  ] \neg  : \big( s,X^{\Th_{\mu^2}}_s (\o) \big)
  \neg  \notin   \neg   O_\d (t, x) \Big\} =
 \t_{\mu^2}  (\o)  \ge  \t_{\mu^1}  (\o)    .
 \eeas
Similarly, it holds    on  $  A_o \cap \{\t_{\mu^2} \neg \le \t_{\mu^1} \} \cap \cN^c $ that   $\t_{\mu^1} = \t_{\mu^2}$.
So
 \bea \label{eq:xxd031a}
 \t_{\mu^1} = \t_{\mu^2} \; \hb{ on $\wt{A} \dfnn ( A \cup   A_o) \cap \cN^c  $. }
 \eea
  Since $      \[t,\t\[  \,  \cap \,  \[ \t_{\mu^1}  \neg \land \neg  \t_{\mu^2}, T \]
   \neg  =  \neg  \[ \t_{\mu^1}  \neg \land \neg  \t_{\mu^2},   \t \[_{A_o} $
 and   $\[t,T\]_A    \cap   \[ \t_{\mu^1}  \neg \land \neg  \t_{\mu^2}, T \]
   \neg  =  \neg  \[ \t_{\mu^1}  \neg \land \neg  \t_{\mu^2}, T \]_A $,
   \eqref{eq:xxd031a} leads to  that
    \beas
      \big( \[t,\t\[ \, \cup \, \[t, T\]_A \big) \cap
  \[ \t_{\mu^1} \land \t_{\mu^2}, T \]_{\cN^c} \subset \[ \t_{\mu^1} \land \t_{\mu^2}, T \]_{\wt{A}  }
  =  \[ \t_{\mu^1}  , T   \]_{\wt{A}  } \, \cap \, \[ \t_{\mu^2}  , T   \]_{\wt{A}  }   .
 \eeas
  Thus it holds  $ds  \neg \times \neg  dP-$a.s. on $ \big( \[t,\t\[ \, \cup \, \[t, T\]_A \big) \cap
  \[ \t_{\mu^1} \land \t_{\mu^2}, T \]  $ that
  $
  \big(\wh{\beta}(\mu^1) \big)_s \neg = \neg  \psi(s,\mu^1_s )  \neg = \neg  \psi(s,\mu^2_s )
   \neg = \neg  \big(\wh{\beta}(\mu^2) \big)_s   $,
 which together with   \eqref{eq:xxd075a}   shows that $\wh{\beta} \in \fB_t$.

  Given $\mu \in \cU_t$, we set $\Th_\mu    \dfnn    \big(t,x,\mu, \beta ( \mu )\big)$  and  $\wh{\Th}_\mu   \neg   \dfnn   \neg   \big(t,x,\mu, \wh{\beta} ( \mu )\big)$.
  For  $  i   =1,\cds \neg, N_m  $,
 analogous to  $\wt{A}^m_i$ of part (1b), $   \cA^{\mu,m}_i \dfnn  \big\{ \big(\t_\mu, X^{\Th_\mu}_{\t_\mu} \big) \in \fD_m (s_i,x_i)   \backslash \underset{j < i}{\cup} \fD_m (s_j,x_j) \big\}       $ belongs to $  \cF_{\t_\mu} \cap \cF_{t_i}$.
    By the continuity of process $ X^{\Th_\mu } $,   $(\t_\mu,X^{\Th_\mu }_{\t_\mu}) \in \pa O_\d (t,x) $, \pas ~ So
           $\{ \cA^{\mu,m}_i \}^{N_m}_{i =1    }  $ forms a partition of $\cN^c_\mu$ for some $P-$null set $\cN_\mu$.
       Then  we can define
           an $\bF-$stopping time $\t^m_\mu \dfnn \sum^{N_m}_{i=1} \b1_{\cA^{\mu,m}_i} \, t_i +  \b1_{\cN_\mu} T  \ge  \t_\mu $ as well as a process
          \bea
     \big(  \beta_m (\mu) \big)_s     & \dneg  \dfnn &  \dneg   \b1_{\{ s < \t^m_\mu  \}} \big( \wh{\beta} (\mu) \big)_s
     \neg + \neg  \b1_{\{ s  \ge  \t^m_\mu  \}} \Big(  \sum^{N_m}_{i=1}   \b1_{\cA^{\mu,m}_i \cap \cA^m_i}
     \big(  \beta^m_i  ( [\mu]^{t_i}   ) \big)_s + \b1_{\cA^m_\mu} \big( \wh{\beta} (\mu) \big)_s \Big) \nonumber    \\
     & \dneg  =  &  \dneg  \b1_{\cA^m_\mu} \big( \wh{\beta} (\mu) \big)_s + \sum^{N_m}_{i=1}
      \b1_{\cA^{\mu,m}_i \cap \cA^m_i}     \Big( \b1_{\{ s < t_i  \}} \big( \wh{\beta} (\mu) \big)_s
        \neg + \neg \b1_{\{ s  \ge  t_i  \}}  \big(  \beta^m_i  ( [\mu]^{t_i}   ) \big)_s \Big) , \q \fa  s  \in [t,T] ,    \label{eq:q211}     \qq
    \eea
   where $ \cA^m_\mu =  \Big( \underset{i=1}{\overset{N_m}{\cup}}
   \big( \cA^{\mu,m}_i \backslash \cA^m_i \big)  \Big) \cup \cN_\mu $.

 \ss  We claim that $\beta_m$ is a $\fB_t-$strategy.
 \if{0}
  To see this, let us first discuss the $\bF-$progressive measurability of   process  $\beta_m (\mu) $.
     Let $s \in [t,T]$ and $V \in \sB \big(\hV  \big)  $. As $\cD^m_\mu \dfnn \[t, \t^m_\mu \[
      \, \cap \, ( [t,s] \neg \times  \neg  \O) \in \sP$,
     similar to \eqref{mu_progressive}, one can deduce that
             \bea
         \big\{ (r,  \o) \neg \in \neg  \cD^m_\mu  \neg  :   \big( \beta_m (\mu) \big)_r (\o) \neg \in \neg V \big\}
        \neg   =  \neg     \big\{ (r,  \o)  \neg \in \neg   \cD^m_\mu  \neg  :
             \big( \wh{\beta} ( \mu )   \big)_r (\o) \neg   \in \neg V \big\}
               \neg   \in \neg  \sB ([t,s] )  \neg \otimes \neg  \cF_s  . \q
               \label{eq:xux517c}
       \eea
     Given $i  = 1, \cds \neg , N_m$, we set $\ol{\cA}^{\mu,m}_i \dfnn \big( \cA^{\mu,m}_i \backslash \cA^m_i \big)  \cup \cN_\mu \in \cF_{t_i}$. If $s < t_i   $, both
     $     \cD^{\mu,m}_i \dfnn  \[\t^m_\mu, T\]_{\cA^{\mu,m}_i \cap \cA^m_i} \cap  ( [t,s] \times \O )
     =  ( [t_i,T] \, \cap \, [t,s] )  \times  \big(\cA^{\mu,m}_i \cap \cA^m_i \big)
        $ and  $     \wh{\cD}^{\mu,m}_i \dfnn  \[\t^m_\mu, T\]_{\ol{\cA}^{\mu,m}_i} \cap  ( [t,s] \times \O )
     =  ( [t_i,T] \, \cap \, [t,s] )  \times  \ol{\cA}^{\mu,m}_i  $ are empty.
        Otherwise, if $s \ge t_i   $,
      Otherwise, both $ \cD^{\mu,m}_i =   [t_i, s]    \times \big(\cA^{\mu,m}_i \cap \cA^m_i \big)$ and
     $     \wh{\cD}^{\mu,m}_i       =  [t_i, s]   \times  \ol{\cA}^{\mu,m}_i  $ belong to
     $     \sB\big([t_i,s]\big)  \neg \otimes \neg  \cF_s  $.
      Using a similar argument to \eqref{mu_progressive}
      on the $\bF-$progressive measurability of process $ \beta^m_i \big( [\mu]^{t_i} \big) $ yields that
             \beas
    \big\{ (r,  \o) \neg \in  \neg  \cD^{\mu,m}_i \neg :  \big( \beta_m (\mu) \big)_r (\o) \neg \in \neg V \big\}
   & \tneg \dneg = & \tneg  \dneg  \big\{ (r,  \o)  \neg \in  \neg  \cD^{\mu,m}_i \neg : \big( \beta^m_i \big( [\mu]^{t_i} \big) \big)_r (\o)   \neg \in \neg  V \big\}
       \neg \in \neg  \sB\big([t_i,s]\big)  \neg \otimes \neg  \cF_s
        \neg   \subset \neg  \sB\big([t,s]\big)  \neg \otimes \neg  \cF_s \\
 \hb{and } \;  \big\{ (r,  \o) \neg \in  \neg  \cD^{\mu,m}_i \neg :  \big( \beta_m (\mu) \big)_r (\o) \neg \in \neg V \big\}
   & \tneg \dneg = & \tneg  \dneg  \big\{ (r,  \o)  \neg \in  \neg  \wh{\cD}^{\mu,m}_i \neg :  \big( \wh{\beta} (  \mu ) \big)_r (\o)   \neg \in \neg  V \big\}
       \neg \in \neg       \sB\big([t,s]\big)  \neg \otimes \neg  \cF_s     ,
     \eeas
    both of  which together with     \eqref{eq:xux517c} shows  the $\bF-$progressive measurability of
      $\beta_m (\mu) $.
 \fi
       Using a similar argument to that in part (1b) for the measurability of the pasted control $\mu^m$, one can
       deduce that the process $\beta_m (\mu) $ is  $\bF-$progressively measurable.
        For $i =1, \cds \neg , N_m $, let $C^m_i  >0 $ be the constant associated to $\beta^m_i$ in Definition \ref{def:strategy} (i).
        Setting $ C_m =C_\beta \vee \k \vee \max\{ C^m_i: i =1,\cds \neg , N_m \} $, we can deduce from
         \eqref{eq:xux717a} and (A-u)  that
        $ds \times dP-$a.s.
         \bea
    && \hspace{-0.7cm}   \big[ ( \beta_m (\mu))_s \big]_{\overset{}{\hV}}
          \neg =  \neg   \b1_{\{ s < \t^m_\mu  \}} \big[\big( \wh{\beta} (\mu) \big)_s\big]_{\overset{}{\hV}}
     \neg + \neg  \b1_{\{ s  \ge  \t^m_\mu  \}}  \Big( \sum^{N_m}_{i=1}   \b1_{\cA^{\mu,m}_i \cap \cA^m_i}
     \big[ \big(  \beta^m_i  ( [\mu]^{t_i}   ) \big)_s\big]_{\overset{}{\hV}} + \b1_{\cA^m_\mu} \big[\big( \wh{\beta} (\mu) \big)_s\big]_{\overset{}{\hV}} \Big) \qq \qq \nonumber \\
     &&    \le      \big( \b1_{\{ s < \t^m_\mu  \}}
     \neg + \neg  \b1_{\{ s   \ge  \t^m_\mu  \}} \b1_{\cA^m_\mu} \big)
     \big( \k  \neg + \neg  ( C_\beta  \neg \vee \neg  \k)  [\mu_s ]_{\overset{}{\hU}} \big)
     \neg + \neg  \b1_{\{ s  \ge  \t^m_\mu  \}}    \sum^{N_m}_{i=1}   \b1_{\cA^{\mu,m}_i \cap \cA^m_i}
     \big( \k  \neg + \neg  C^m_i \big[   [\mu]^{t_i}_s\big]_{\overset{}{\hU}} \big)
      \neg \le \neg  \k  \neg + \neg   C_m  [\mu_s ]_{\overset{}{\hU}} .   \qq \q  \label{eq:xux717}
         \eea
      Let $E \int_t^T \neg \big[   \mu_s   \big]^q_{\overset{}{\hU}} ds   < \infty$ for some  $q >2$. It follows from
      \eqref{eq:xux717} that
      \beas
      E \neg  \int_t^T  \neg \big[ ( \beta_m (\mu))_s \big]^q_{\overset{}{\hV}} \, ds
      \le 2^{q-1} \k^q T +  2^{q-1}   C^q_m E  \neg \int_t^T  \neg \big[   \mu_s   \big]^q_{\overset{}{\hU}} \, ds < \infty.
      \eeas
       Hence $ \beta_m (\mu) \in \cV_t$.

\ss Let    $\mu^1, \mu^2 \neg  \in \neg  \cU_t $ such that  $\mu^1  \neg = \neg   \mu^2 $,
$ds  \neg \times \neg  d P -$a.s. on $ \[t,\t\[ \,  \cup \,  \[\t, T\]_A   $  for some  $ \t \neg \in  \neg   \cS_{t,T}$
  and $ A  \neg  \in  \neg  \cF_\t  $.
 As $ \wh{\beta} (\mu^1)  \neg  = \neg   \wh{\beta} (\mu^2) $, $ds  \neg \times \neg  dP-$a.s.
 on $ \[t,\t\[  \, \cup \,  \[\t, T\]_A   $
 by Definition \ref{def:strategy},
  it holds    $ds  \neg \times \neg  dP-$a.s. on $       \big( \[t,\t\[ \, \cup \, \[t, T\]_A \big)  \neg \cap \neg
  \[ t, \t^m_{\mu^1} \neg \land \neg  \t^m_{\mu^2} \[ $ that
 \bea \label{eq:xxd075}
 \big(\beta_m(\mu^1) \big)_s  = \big( \wh{\beta}  (\mu^1) \big)_s  = \big( \wh{\beta}  (\mu^2) \big)_s  = \big(\beta_m(\mu^2) \big)_s \,.
 \eea
 Definition \ref{def:strategy} also shows that
  $ \big( \mu^1, \beta(\mu^1)\big) \neg  = \neg \big( \mu^2, \beta(\mu^2) \big) $, $ds  \neg \times \neg  dP-$a.s. on $ \[t,\t\[ \, \cup \, \[\t, T\]_A    $.  So we again have
   \eqref{eq:xxd071a}   except on a $P-$null set $ \cN $, and
   \eqref{eq:xxd031a} still holds  on $\wt{A} \dfnn ( A \cup   A_o) \cap \cN^c  $
   with $A_o \neg =  \neg   \{ \t   \neg  \ge  \neg    {\t_{\mu^1}  \neg \land \neg  \t_{\mu^2}}  \}   $.
    Plugging  \eqref{eq:xxd031a}    into \eqref{eq:xxd071a} yields that
      \bea  \label{eq:xxd031d}
     X^{\Th_{\mu^1}}_{ \t_{\mu^1} }     = X^{\Th_{\mu^2}}_{\t_{\mu^2}} \hb{ holds   on }
      \wt{A}    .
    \eea
    Given $i  \neg  = \neg  1, \cds \neg , N_m  $.
  since it holds $ds \neg \times \neg  dP-$a.s. on $\big( \[t,\t\[ \, \cup \, \[\t, T\]_A \big)  \neg \cap \neg
  \big([t_i,T]  \neg \times \neg  \O \big)   \neg  = \neg  \[ t_i, \t  \neg \vee \neg  t_i \[ \,\cup \, \[\t
   \neg \vee \neg  t_i, T \]_A  $ that
 $  \big( [\mu^1]^{t_i} \big)_s   \neg  = \neg  \mu^1_s   \neg = \neg   \mu^2_s  \neg = \neg \big( [\mu^2]^{t_i} \big)_s $,
  taking $ ( \t, A)  \neg  =  \neg  (   \t \vee t_i, A )   $
    in Definition \ref{def:strategy} with respect to $\beta^m_{i} $ yields that
 for $ds  \neg \times \neg  dP-$a.s. $(s,\o)  \neg \in  \neg     \[ t_i, \t \vee t_i \[ \,\cup \, \[\t \vee t_i, T \]_A
 \neg  = \neg  \big( \[t,\t\[ \, \cup \, \[t, T\]_A \big)  \neg \cap \neg
  \big([t_i,T]  \neg \times \neg  \O \big) $
 \bea   \label{eq:xxd073}
   \big( \beta^m_{i } ( [ \mu^1 ]^{t_i}  ) \big)_s (\o)  = \big( \beta^m_{i } ( [ \mu^2 ]^{t_i} ) \big)_s (\o) .
 \eea
   Given $\o \in  \cA_{i} \dfnn  \wt{A}  \cap \cA^{\mu^1,m}_i    $,   \eqref{eq:xxd031a}  and \eqref{eq:xxd031d}
   imply that
  \beas
  \Big(\t_{\mu^2}(\o), X^{\Th_{\mu^2}}_{\t_{\mu^2}(\o)} (\o)\Big)
  = \Big(\t_{\mu^1}(\o), X^{\Th_{\mu^1}}_{\t_{\mu^1}(\o)} (\o)\Big)
  \in      \fD_m (s_i,x_i)   \backslash \underset{j < i}{\cup} \fD_m (s_j,x_j) , \q \hb{i.e., $ \o \in  \cA^{\mu^2,m }_{i }$.}
 \eeas
       So $ \cA_{i} \subset \cA^{\mu^1,m}_i  \cap \cA^{\mu^2,m }_{i }$,
  and it follows that  $  \b1_{\cA_i} \t^m_{\mu^1} = \b1_{\cA_i} t_i = \b1_{\cA_i} \t^m_{\mu^2}$.  Then one can deduce that
   \bea   \label{eq:xvx091}
   \big( \[t,\t\[ \,  \cup \,  \[t, T\]_A \big)    \cap      \[ \t^m_{\mu^1} \neg \land \neg  \t^m_{\mu^2},T \]_{\cA_{i} \cap \cA^m_i}
  & \tneg   = & \tneg   \big( \[t,\t\[ \, \cup \, \[t, T\]_A \big)  \neg  \cap  \neg  \big( [t_i,T]  \neg \times \neg  ( \cA_{i} \neg   \cap \cA^m_i ) \big)  \qq   \qq   \nonumber   \\
   & \tneg \subset & \tneg    [t_i, T]  \neg \times \neg
    \big( \cA^{\mu^1,m}_{i}   \neg  \cap    \cA^{\mu^2,m}_{i}  \neg  \cap \cA^m_i  \big) \, ,    \q
  \eea
   which together  with      \eqref{eq:xxd073}   shows   that
 for $ds \times dP-$a.s. $(s,\o) \in    \big( \[t,\t\[ \, \cup \, \[t, T\]_A \big) \cap   \[ \t^m_{\mu^1} \land \t^m_{\mu^2},T \]_{\cA_{i} \cap \cA^m_i}   $
  \bea    \label{eq:xvx093}
 \big(\beta_m(\mu^1) \big)_s (\o) = \big( \beta^m_{i} ( [\mu^1]^{t_i}  ) \big)_s (\o)
 = \big( \beta^m_{i} ( [\mu^2]^{t_i}  ) \big)_s (\o)
 =   \big(\beta_m(\mu^2) \big)_s (\o) .
 \eea
 Analogous to \eqref{eq:xvx091},
  $  \big( \[t,\t\[ \, \cup \, \[t, T\]_A \big) \cap   \[ \t^m_{\mu^1}
  \neg \land \neg \t^m_{\mu^2},T \]_{\cA_{i} \backslash \cA^m_i}
    \neg \subset \neg    [t_i, T]  \neg \times \neg
    \big( ( \cA^{\mu^1,m}_{i}  \backslash \cA^m_i )  \cap (  \cA^{\mu^2,m}_{i}  \backslash  \cA^m_i )  \big) $. So
  \eqref{eq:xxd075} also  holds   $ds \times dP-$a.s. on  $     \big( \[t,\t\[ \, \cup \, \[t, T\]_A \big)
    \cap   \[ \t^m_{\mu^1} \land \t^m_{\mu^2},T \]_{\cA_{i} \backslash \cA^m_i}   $.
  Combining this with \eqref{eq:xvx093} and then
  letting $i$ run over $\{ 1, \cds \neg , N_m \}$ yield that
  \bea \label{eq:xvx095}
  \big(\beta_m(\mu^1) \big)_s   =   \big(\beta_m(\mu^2) \big)_s   , \q  ds \times dP-a.s. \hb{ on }
  \big( \[t,\t\[ \, \cup \, \[t, T\]_A \big) \cap   \[ \t^m_{\mu^1} \land \t^m_{\mu^2},T \]_{A \cup A_o } .
  \eea
  As $ \[ \t^m_{\mu^1} \land \t^m_{\mu^2},T \]_{A^c \cap A^c_o } \subset
  \[ \t_{\mu^1} \land \t_{\mu^2},T \]_{A^c \cap A^c_o } \subset \[\t,T\]_{A^c \cap A^c_o }
  \subset \[\t,T\]_{A^c   } $, one can deduce  that
  $ \big( \[t,\t\[ \, \cup \, \[t, T\]_A \big) \cap   \[ \t^m_{\mu^1} \land \t^m_{\mu^2},T \]_{A \cup A_o }
  = \big( \[t,\t\[ \, \cup \, \[t, T\]_A \big) \cap   \[ \t^m_{\mu^1} \land \t^m_{\mu^2},T \] $.
  Therefore, \eqref{eq:xvx095}  together with \eqref{eq:xxd075}  implies that  $\beta_m \in \fB_t$.

    \ss \no {\bf 1e)} {\it
       Next,   let $\mu \neg  \in  \neg    \cU_t $
     and    $ \Th^m_\mu   \neg  \dfnn   \neg    \big( t , x , \mu ,  \beta_m (\mu) \big)$.
   We shall do similar estimates to those in   part \(c\) to conclude
   \bea \label{eq:xqxqx223}
             w_1(t,x)  \le    \underset{\beta \in \fB_t }{\essinf} \; \esup{\mu \in \cU_t} \;  Y^{t,x,\mu, \beta (  \mu ) }_t \Big(\t_{\beta,\mu},
    \wt{\f}\big(\t_{\beta,\mu},  X^{t,x,\mu, \beta (  \mu ) }_{\t_{\beta,\mu}} \big) \Big)  , \q   \pas
   \eea
    }

                As     $  \beta_m  (\mu)     =      \wh{\beta} (\mu)      =      \beta (\mu)   $
                on $\[t,\t_\mu \[$,
   taking $(\t, A)    =    (\t_\mu, \es)$ in \eqref{eq:p611} yields that \pas
 \bea \label{eq:xxa614}
    X^{ \Th^m_\mu  }_s = X^{\wh{\Th}_\mu}_s = X^{\Th_\mu}_s \in \ol{O}_\d (x) , \q  \fa s \in [t,\t_\mu] .
 \eea
  Thus,  for any   $\eta \in \hL^p \big( \cF_{\t_\mu} \big)  $,
       the  BSDE$ \Big(t, \eta ,   f^{\Th^m_\mu}_{\t_\mu}    \Big)$
      and the BSDE$ \Big(t, \eta , f^{\Th_\mu}_{\t_\mu}    \Big)$
      are essentially the same. To wit,
   \bea
   \big(Y^{\Th^m_\mu} (\t_\mu, \eta),Z^{\Th^m_\mu} (\t_\mu, \eta)  \big)
         =    \big(Y^{\Th_\mu} (\t_\mu, \eta),Z^{\Th_\mu} (\t_\mu, \eta)  \big)    .   \label{eq:q104}
    \eea
  \if{0}
 By     \eqref{eq:xux611} and \eqref{eq:xxa614}, it holds     \pas ~ that
 \bea
 X^{\Th^m_\mu }_{  \t^m_\mu \land s} & \tneg \neg =&  \neg \tneg  X^{\Th^m_\mu }_{\t_\mu \land s} \neg +
  \dneg \int_{\t_\mu \land s}^{ \t^m_\mu \land s}  \dneg  b \Big(r, X^{\Th^m_\mu }_r,\mu_r, \big(\beta_m( \mu )\big)_r \Big)  dr
  \neg + \neg  \int_{\t_\mu \land s}^{ \t^m_\mu \land s}  \dneg  \si \Big(r, X^{\Th^m_\mu }_r,\mu_r, \big(\beta_m ( \mu)\big)_r  \Big)  dB_r,  \nonumber \\
 & \tneg \neg  =& \neg  \tneg  X^{\Th^m_\mu }_{\t_\mu \land s}  \neg + \dneg  \int_{\t_\mu \land s}^{ \t^m_\mu \land s}
   \dneg  b \Big(r, X^{\Th^m_\mu }_r,\mu_r, v_\sharp  \Big)  dr
  \neg + \neg  \int_{\t_\mu \land s}^{ \t^m_\mu \land s}  \dneg  \si \Big(r, X^{\Th^m_\mu }_r,\mu_r,  v_\sharp  \Big)  dB_r
   \neg = \neg  X^{\Th^m_\mu }_{\t_\mu \land s} \neg = \neg  X^{ \Th_\mu }_{\t_\mu \land s} , ~ s  \neg \in \neg  [t,T].   \qq \q   \label{eq:xux311b}
 \eea
 Since $\cZ^\eta = 0$ on $ \[\t_\mu, T\] $  by \eqref{eq:xux313},  we can deduce   from  \eqref{eq:xux614} that
  \beas
  \b1_{\{\t_\mu \le r < \t^m_\mu\}}
    f \Big( r, X^{\Th^m_\mu}_r,  \cY^\eta_r, \cZ^\eta_r, \mu_r, \big(\beta_m( \mu )\big)_r \Big)
    =  \b1_{\{\t_\mu \le r < \t^m_\mu\}}
    f \Big( r, X^{\Th^m_\mu}_r,  \cY^\eta_r, 0, \mu_r, v_\sharp \Big)
    =0,  \q r \in [t,T].
  \eeas
 Then it follows from \eqref{eq:xux311b}  that
 \beas
        \cY^\eta_s  & =&
       \eta  \neg + \neg    \int_s^T \neg \b1_{\{r < \t_\mu\}}
    f \Big( r, X^{\Th^m_\mu}_r,  \cY^\eta_r, \cZ^\eta_r, \mu_r, \big(\beta_m( \mu )\big)_r \Big)  \, dr
   \neg-\neg   \int_s^T \neg \cZ^\eta_r d B_r       \\
    & = &  \eta  \neg + \neg    \int_s^T \b1_{\{r < \t^m_\mu\}}
    f \Big(r, X^{\Th^m_\mu}_r,  \cY^\eta_r, \cZ^\eta_r, \mu_r, \big(\beta_m( \mu )\big)_r \Big)  \, dr
   \neg-\neg   \int_s^T \neg \cZ^\eta_r d B_r  , \q s \in [t,T]  .
    \eeas
    Thus $ \big( \cY^\eta,\cZ^\eta   \big) $
    also solve   BSDE$ \Big(t, \eta ,   f^{\Th^m_\mu}_{\t^m_\mu}    \Big)$, i.e.,
     \beas
   \big( \cY^\eta,\cZ^\eta   \big)
   \neg  = \neg   \big(Y^{\Th^m_\mu} (\t^m_\mu, \eta),Z^{\Th^m_\mu} (\t^m_\mu, \eta)  \big)    .
     \eeas
    In particular,
     \bea  \label{eq:xux417b}
      Y^{\Th^m_\mu}_{\t_\mu} (\t^m_\mu, \eta) = \cY^\eta_{\t_\mu}
      = Y^{\Th^m_\mu}_{\t_\mu} (\t_\mu, \eta) = \eta .
     \eea
  \fi

    Given $A \in \cF_t$, similar to \eqref{eq:xvx103},  we can deduce  from  \eqref{eq:xxa614}  that
    \if{0}
  \beas
\b1_A  X^{\Th^m_\mu }_{  \t^m_\mu  \land s} & \tneg \neg =&  \neg \tneg \b1_A   X^{\Th^m_\mu }_{\t_\mu \land s} \neg +
  \neg \b1_A \neg \int_{\t_\mu \land s}^{ \t^m_\mu  \land s}    b \Big(r, X^{\Th^m_\mu }_r,\mu_r, \big(\beta_m (\mu)\big)_r \Big)  dr
  \neg + \neg \b1_A  \neg \int_{\t_\mu \land s}^{ \t^m_\mu  \land s}
    \si \Big(r, X^{\Th^m_\mu }_r,\mu_r, \big(\beta_m ( \mu )\big)_r  \Big)  dB_r,    \\
 & \tneg \neg  =& \neg  \tneg \b1_A  X^{\Th_\mu }_{\t_\mu \land s}  \neg + \neg  \int_{\t_\mu \land s}^{ \t^m_\mu  \land s}
     \b1_A  b \Big(r, X^{\Th^m_\mu }_{\t^m_\mu  \land r}, \mu_r, \psi(r, \mu_r) \Big)     dr
  \neg + \neg \int_{\t_\mu \land s}^{ \t^m_\mu  \land s} \b1_A  \si \Big(r, X^{\Th^m_\mu }_{\t^m_\mu  \land r} , \mu_r,  \psi(r, \mu_r) \Big)  dB_r   , \q  s  \neg \in \neg  [t,T].
 \eeas
 It  follows that
 \fi
  \beas
 \b1_A \underset{r \in [t,s]}{\sup} \big| X^{\Th^m_\mu }_{  \t^m_\mu  \land r} \neg  - \neg  X^{\Th_\mu }_{\t_\mu \land r} \big|
 & \tneg \le & \tneg  \int_{\t_\mu \land s}^{ \t^m_\mu  \land s} \neg
  \b1_A  \big|   b \big(r, X^{\Th^m_\mu }_{\t^m_\mu  \land r}, \mu_r, \psi(r, \mu_r) \big) \big|    dr   \\
 &\tneg  & \tneg   + \neg  \underset{r \in [t,s]}{\sup} \left| \int_{\t_\mu \land r}^{ \t^m_\mu  \land r}  \neg \b1_A  \si \big(r', X^{\Th^m_\mu }_{\t^m_\mu  \land r'}, \mu_{r'},  \psi(r', \mu_{r'}) \big)     dB_{r'}   \right|  , \q   s  \neg \in \neg  [t,T].
 \eeas
 where we used the fact that $ \beta_m (\mu) =  \wh{\beta} (\mu)   =  \beta_\psi (\mu) $ on $\[\t_\mu, \t^m_\mu\[$.
  Let  $\wt{C}(\k,x,\d)$   denote a generic constant,  depending on $\k \neg + \neg |x| \neg + \neg \d$, $C^{\wt{\f}}_{x,\d} \dfnn
        \sup \big\{|\wt{\f} (s,\fx)| : (s,\fx) \in   \ol{O}_{\d + 3}(t,x) \cap ([t,T] \times \hR^k ) \big\} $,
   $T$, $\g $,  $p$ and $|g(0)|$,        whose form  may vary from line to line.
 Since $\t^m_\mu - \t_\mu \le \sum^{N_m}_{i=1} \b1_{ \cA^{\mu,m}_i } \, 2 \d^m_{s_i,x_i} < \frac{2}{m} $, \pas, using similar arguments to those that lead to \eqref{eq:xvx113} and using an analogous  decomposition and estimation to \eqref{eq:xvx055},
 we can deduce that
           \beas
 && \hspace{-1.2cm} E \bigg[ \b1_A  \underset{r \in [t,s]}{\sup} \big| X^{\Th^m_\mu }_{  \t^m_\mu  \land r} \neg  - \neg X^{\Th_\mu }_{\t_\mu \land r} \big|^2 \bigg]   \\
  & &  \le      \neg     4 E  \neg  \int_{\t_\mu \land s}^{ \t^m_\mu  \land s}  \neg
  \b1_A  \big|   b \Big(r, X^{\Th^m_\mu }_{\t^m_\mu  \land r}, \mu_r, \psi(r, \mu_r) \Big) \big|^2    dr
  \neg    +  \neg     8 E  \neg   \int_{\t_\mu \land s}^{ \t^m_\mu  \land s}  \neg  \b1_A
  \big|  \si \Big(r, X^{\Th^m_\mu }_{\t^m_\mu  \land r}, \mu_r,  \psi(r, \mu_r) \Big) \big|^2     dr        \\
 & &  \le   \neg    24 \g^2     \int_t^s    E \Big[ \b1_A \underset{r' \in [t,r]}{\sup} \big|   X^{\Th^m_\mu }_{\t^m_\mu  \land r'}   \neg -   \neg    X^{\Th_\mu }_{\t_\mu \land r'}  \big|^2 \Big] dr
 \neg + \neg  \frac{ \wt{C}(\k,x,\d)}{m} P(A) , \q \fa  s  \neg \in \neg  [t,T].
 \eeas
 \if{0}
   Then an application of    Gronwall's inequality yields that
  \beas
  E \bigg[ \b1_A \underset{r \in [t,s]}{\sup} \big| X^{\Th^m_\mu }_{  \t^m_\mu  \land r} -X^{\Th_\mu }_{\t_\mu \land r} \big|^2 \bigg]
  \le \frac{48}{m } \g^2   \big[  |x| \neg + \neg \d  \neg + \neg   (1  \neg + \neg  \k)^2 \, \big]^2 P  ( A )
     e^{24 \g^2  (s-t)}  , \q \fa  s  \neg \in \neg  [t,T].
  \eeas
  In particular,   $    E \bigg[ \b1_A \underset{r \in [t,T]}{\sup} \big| X^{\Th^m_\mu }_{  \t^m_\mu  \land r}  \neg - \neg X^{\Th_\mu }_{\t_\mu \land r} \big|^2 \bigg]
   \neg \le \neg  \frac{ C(\k,x,\d)}{m} P(A)$.  Here $C(\k,x,\d)$
   denotes a generic constant,  depending on $\k \neg + \neg |x| \neg + \neg \d$, $C^\f_{x,\d}  $,
   $T$, $\g $,  $p$ and $|g(0)|$,
        whose form  may vary from line to line.    Letting $A$ vary in $\cF_t$ yields that
  \bea   \label{eq:xvx043}
  E \bigg[  \underset{r \in [t,T]}{\sup} \big| X^{\Th^m_\mu }_{  \t^m_\mu  \land r}  \neg - \neg X^{\Th_\mu }_{\t_\mu \land r} \big|^2  \Big| \cF_t \bigg]
  \le  \frac{ C(\k,x,\d)}{m}, \q \pas
  \eea
  \fi
  Then similar to \eqref{eq:xvx043}, an application of    Gronwall's inequality leads to that
 \bea  \label{eq:xvx119}
  E \bigg[  \underset{r \in [t,T]}{\sup} \big| X^{\Th^m_\mu }_{  \t^m_\mu  \land r}  \neg - \neg X^{\Th_\mu }_{\t_\mu \land r} \big|^2  \Big| \cF_t \bigg]
  \le  \frac{ \wt{C}(\k,x,\d)}{m}, \q \pas
  \eea

 Let $i = 1, \cds \neg , N_m $ and  set  $\Th^{m,t_i}_\mu \dfnn \Big(t_i, X^{\Th^m_\mu}_{t_i}, [\mu]^{t_i},
 [\beta_m(\mu)]^{t_i}  \Big) $. Similar to \eqref{eq:xxd213}, it holds \pas ~ that
  \bea  \label{eq:xxd213b}
   Y^{\Th^m_\mu}_{t_i} \Big(T, g \Big(X^{\Th^m_\mu}_T\Big) \Big)  = J \big( \Th^{m,t_i}_\mu \big) .
  \eea
   Since    $ \big[ \beta_m(\mu)\big]^{t_i}_r (\o) = \big(\beta_m(\mu)\big)_r (\o)
 = \big( \beta^m_{ i } \big( [\mu ]^{t_i} \big) \big)_r (\o)$
    for any $  (r, \o) \in [t_i,T] \times \big( \cA^{\mu,m}_{i} \cap \cA^m_i \big)$,
   one can deduce from \eqref{eq:xxd213b}, \eqref{eq:xxa273} and \eqref{eq:s025} that  it holds \pas ~ on $ \cA^{\mu,m}_{i} \cap \cA^m_i     \in \cF_{t_i} $ that
   $$
    Y^{ \Th^m_\mu}_{\t^m_\mu } \neg  \Big(T, g \Big(X^{\Th^m_\mu}_T\Big) \neg \Big)
   \neg = \neg  Y^{ \Th^m_\mu}_{t_i} \neg  \Big(T, g \Big(X^{\Th^m_\mu}_T\Big) \neg \Big) 
   \neg  =  \neg   J \Big(  t_i, X^{ \Th^m_\mu}_{\t^m_\mu}, [\mu]^{t_i},  \beta^m_{ i } \big( [\mu]^{t_i} \big)  \Big)
    \neg   \le \neg   J \Big(  t_i, X^{\Th_\mu}_{\t_\mu}, [\mu]^{t_i}, \beta^m_{ i } \big( [\mu]^{t_i} \big)  \neg  \Big)
    \neg   +   c_0 \Big| X^{ \Th^m_\mu}_{\t^m_\mu} \neg -  X^{\Th_\mu}_{\t_\mu} \Big|^{2/p} .
   $$
   As $\dis \big| X^{\Th_\mu}_{\t_\mu} - x_i \big|^{2/p}
       <  ( \d^m_{s_i,x_i})^{2/p} < m^{-2/p} \le 1/m $ on $\cA^{\mu, m}_i$,
    we can also deduce from
 \eqref{eq:s025},     \eqref{eq:xxd083}, \eqref{eq:xux391} and the continuity of $\wt{\f}$
   that   it holds \pas ~ on $   \cA^{\mu,m}_{i} \cap \cA^m_i    $  that
   \beas
    J \Big(  t_i, X^{\Th_\mu}_{\t_\mu}, [\mu]^{t_i}, \beta^m_{ i } \big( [\mu]^{t_i} \big)  \Big)
     &    \dneg  \dneg   \le  &  \dneg  \dneg
         J \Big( t_i,x_i ,  [\mu]^{t_i}, \beta^m_{ i } \big( [\mu]^{t_i} \big)  \Big)
  \neg + \neg  \frac{c_0}{m}
     \neg   \le   \neg       I \big(t_i, x_i , \beta^m_{i} \big)   \neg +  \neg  \frac{c_0}{m}
    \neg  \le  \neg  \wt{\f}\big(t_i, x_i  \big)    \neg  +  \neg    \frac{c_0}{m} \\
    & \dneg  \dneg \le &  \dneg  \dneg   \wt{\f}\big(s_i, x_i  \big)    \neg  +  \neg    \frac{c_0}{m}
    \neg  \le   \neg     \wt{\f}\big(\t_\mu, X^{ \Th_\mu}_{\t_\mu}  \big)  \neg   +   \neg    \frac{c_0}{m}
    \dfnn \eta^m_\mu  \in \hL^\infty( \cF_{\t_\mu} )  .
   \eeas
 Thus it holds \pas ~ on $  \cup^{N_m}_{i=1} \big( \cA^{\mu,m}_{i} \cap \cA^m_i \big) $ that
   \beas
   Y^{\Th^m_\mu}_{\t^m_\mu} \Big(T, g \Big(X^{\Th^m_\mu}_T\Big)\Big)    \le
   \eta^m_\mu + c_0 \Big| X^{ \Th^m_\mu}_{\t^m_\mu} \neg -  X^{\Th_\mu}_{\t_\mu} \Big|^{2/p}
         \neg \dfnn \neg  \wt{\eta}^m_\mu  \in \hL^p \big( \cF_{\t^m_\mu} \big)    .
      \q   
  \eeas
  By \eqref{eq:xvx071}, it holds \pas ~ that
    \bea      \label{eq:j335}
     \Big|  Y^{\Th_\mu}_{t}  ({\t_\mu},  \eta^m_\mu  )
   \neg  - \neg  Y^{\Th_\mu}_{t}  \big({\t_\mu},    \wt{\f}\big({\t_\mu} , X^{\Th_\mu}_{{\t_\mu} }  \big)   \big)  \Big|^p
       \le   c_0    E \Big[ \big|   \eta^m_\mu      \neg  -  \neg
        \wt{\f}\big({\t_\mu} , X^{\Th_\mu}_{{\t_\mu} }  \big)     \big|^p   \Big|\cF_t\Big]
        \le  \frac{c_0}{m^p}          .
   \eea

     \if{0}
  Let $(Y^{m,\mu}, Z^{m,\mu}) \in \hG^p_\bF([0,T])$ be the unique solution of the following BSDE with zero generator:
   \beas
   Y^{m,\mu}_s = Y^{\Th^m_\mu}_{\t_\mu}  (\t^m_\mu, \eta^m_\mu )  - \int_s^T Z^{m,\mu}_r d B_r      , \q    s \in [0,T] .
   \eeas
     For any  $s \in [0,T]$, one can deduce that
   \beas
      Y^{m,\mu}_{{\t_\mu} \land s} = E[Y^{m,\mu}_{\t_\mu \land s} | \cF_{\t_\mu}] = E \bigg[ Y^{\Th^m_\mu}_{\t_\mu}  (\t^m_\mu, \eta^m_\mu )  - \int_{\t_\mu \land s}^T Z^{m,\mu}_r d B_r \Big| \cF_{\t_\mu}  \bigg]
   = Y^{\Th^m_\mu}_{\t_\mu}  (\t^m_\mu, \eta^m_\mu )  - \int_{{\t_\mu} \land s}^{\t_\mu} Z^{m,\mu}_r d B_r , \q  \pas ~
   \eeas
       By the continuity of process $Y^{m,\mu}$,  it holds \pas ~ that
    \bea   \label{eq:xvx031}
   Y^{m,\mu}_{{\t_\mu} \land s}  = Y^{\Th^m_\mu}_{\t_\mu}  (\t^m_\mu, \eta^m_\mu )  - \int_{{\t_\mu} \land s}^{\t_\mu} Z^{m,\mu}_r d B_r  = Y^{\Th^m_\mu}_{\t_\mu}  (\t^m_\mu, \eta^m_\mu )  - \int_s^T \b1_{\{r < \t_\mu\}} Z^{m,\mu}_r d B_r     , \q    s \in [0,T] .
   \eea
   Thus, we see that  $(Y^{m,\mu}_s, Z^{m,\mu}_s) \neg = \neg  \big(  Y^{m,\mu}_{{\t_\mu} \land s}, \b1_{\{s < \t_\mu\}} Z^{m,\mu}_s \big),
   s  \neg \in \neg  [0,T]   $. Also,
   taking $[\cd|\cF_{{\t_\mu} \land s}]$ in \eqref{eq:xvx031} shows that \pas
   \beas
    Y^{m,\mu}_s  \neg = \neg   Y^{m,\mu}_{{\t_\mu} \land s}  \neg = \neg  E\big[Y^{\Th^m_\mu}_{\t_\mu}  (\t^m_\mu,  \eta^m_\mu ) \big|\cF_{{\t_\mu} \land s}\big] ,
    \q   \fa  s  \neg \in \neg  [0,T] .
   \eeas

   On the other hand,
   let $(\wt{Y}^{m,\mu}, \wt{Z}^{m,\mu} ) \in \hG^p_\bF([0,T])$ be the unique solution of the following BSDE with zero generator:
   \bea  \label{eq:xvx036}
   \wt{Y}^{m,\mu}_s =    \eta^m_\mu     - \int_s^T \wt{Z}^{m,\mu}_r d B_r      , \q    s \in [0,T] .
   \eea
   Similar to $(Y^{m,\mu}, Z^{m,\mu}) $, it holds \pas ~ that
     \bea   \label{eq:xvx061}
     \big(\wt{Y}^{m,\mu}_s, \wt{Z}^{m,\mu}_s\big) \neg = \neg \big(  \wt{Y}^{m,\mu}_{{\t_\mu} \land s}, \b1_{\{s < \t_\mu\}} \wt{Z}^{m,\mu}_s \big)
       \q \hb{ and } \q      \wt{Y}^{m,\mu}_s   \neg  = \neg  E  [   \eta^m_\mu  |\cF_{{\t_\mu} \land s}  ] ,
       \q        \fa s \neg \in \neg  [0,T]          .
 \eea
   We   can deduce that      $( \cY^{m,\mu},  \cZ^{m,\mu} )
   \dfnn  \big\{ \big( \b1_{\{s < \t_\mu\}} Y^{m,\mu}_s \neg + \neg  \b1_{\{s  \ge  \t_\mu\}} Y^{\Th^m_\mu}_s  (\t^m_\mu, \eta^m_\mu ),
    \b1_{\{s < \t_\mu\}} Z^{m,\mu}_s  \neg + \neg  \b1_{\{s  \ge  \t_\mu\}} Z^{\Th^m_\mu}_s  (\t^m_\mu, \eta^m_\mu )     \big) \big\}_{s \in [0,T]} \in \hG^p_\bF([0,T]) $
   solves the following BSDE
    \bea
   \cY^{m,\mu}_s & \tneg =& \tneg  \b1_{\{s  \ge  \t_\mu\}} Y^{\Th^m_\mu}_s  (\t^m_\mu, \eta^m_\mu )
   \neg + \neg   \b1_{\{s < \t_\mu\}}Y^{\Th^m_\mu}_{\t_\mu}  (\t^m_\mu, \eta^m_\mu )   \neg - \neg  \b1_{\{s < \t_\mu\}} \int_s^T Z^{m,\mu}_r d B_r
   \neg  = \neg  Y^{\Th^m_\mu}_{{\t_\mu} \vee s}  (\t^m_\mu, \eta^m_\mu )  \neg - \neg \b1_{\{s < \t_\mu\}} \int_s^T \b1_{\{r < \t_\mu\}} Z^{m,\mu}_r    d B_r  \nonumber \\
   & \tneg =& \tneg  \eta^m_\mu  \neg + \neg  \int_{{\t_\mu} \vee s}^T f^{\Th^m_\mu}_{\t^m_\mu} \big(r,Y^{\Th^m_\mu}_r  (\t^m_\mu, \eta^m_\mu ), Z^{\Th^m_\mu}_r  (\t^m_\mu, \eta^m_\mu ) \big) dr
     \neg - \neg  \int_{{\t_\mu} \vee s}^T Z^{\Th^m_\mu}_r  (\t^m_\mu, \eta^m_\mu ) dB_r  \neg - \neg  \int_s^T \b1_{\{r < \t_\mu\}} Z^{m,\mu}_r    d B_r \nonumber   \\
   & \tneg =& \tneg  \eta^m_\mu  \neg + \neg  \int_s^T \b1_{\{r  \ge  {\t_\mu} \}} f^{\Th^m_\mu}_{\t^m_\mu} \big(r,\cY^{m,\mu}_r , \cZ^{m,\mu}_r \big) dr
     \neg - \neg  \int_s^T   \cZ^{m,\mu}_r   dB_r     , \q s \in [0,T] .   \label{eq:xvx034}
   \eea
   Set $(\wh{\cY}^{m,\mu}, \wh{\cZ}^{m,\mu}) \dfnn
      \big( \cY^{m,\mu} - \wt{Y}^{m,\mu}, \cZ^{m,\mu} - \wt{Z}^{m,\mu} \big) \in \hG^p_\bF([0,T])  $.
   Similar to \eqref{eq:xvx023}, we can deduce from \eqref{eq:xvx061} that
   \beas
     E \Big[      \big|   Y^{\Th^m_\mu}_{\t_\mu} (\t^m_\mu, \eta^m_\mu )   -\eta^m_\mu     \big|^p  \big| \cF_t  \Big]
    & \tneg = & \tneg    E \Big[      \big|     \wh{\cY}^{m,\mu}_{\t_\mu} \big|^p   \Big|\cF_t \Big]
        \le   E \Big[ \, \underset{s \in [0,T]}{\sup}  \big|    \wh{\cY}^{m,\mu}_s     \big|^p  \Big|\cF_t  \Big]
     \le   c_0 E \bigg[       \int_{\t_\mu}^T \neg  \big| f^{\Th^m_\mu}_{\t^m_\mu} \big(s,  \wt{Y}^{m,\mu}_s ,
    \wt{Z}^{m,\mu}_s   \big)  \big|^p ds   \Big|\cF_t   \bigg] \\
    & \tneg  = & \tneg  c_0 E \bigg[      \int_{\t_\mu}^{\t^m_\mu} \neg  \big| f  \big(s, X^{\Th^m_\mu}_{\t^m_\mu \land s}, \eta^m_\mu ,
    0 ,u_0, (\beta(\mu^m))_s  \big)  \big|^p ds  \Big| \cF_t  \bigg] , \q \pas
     \eeas
   \fi

    Similar to \eqref{eq:xvx114}, one can deduce that
   \beas
   E \bigg[      \Big|   Y^{\Th^m_\mu}_{\t_\mu} (\t^m_\mu, \eta^m_\mu )   -\eta^m_\mu     \Big|^p  \Big| \cF_t  \bigg]
       & \tneg   \le   & \tneg      c_0 E \bigg[       \int_{\t_\mu}^T \neg  \big| f^{\Th^m_\mu}_{\t^m_\mu} \big(s,  \wt{Y}^{m,\mu}_s ,
    \wt{Z}^{m,\mu}_s   \big)  \big|^p ds   \Big|\cF_t   \bigg]   \\
      & \tneg  = & \tneg    c_0 E \bigg[      \int_{\t_\mu}^{\t^m_\mu} \neg  \Big| f  \Big(s, X^{\Th^m_\mu}_{\t^m_\mu \land s}, \eta^m_\mu ,
    0 ,\mu_s, \psi(s, \mu_s)    \Big)  \Big|^p ds  \Big| \cF_t  \bigg]     ,    \q    \pas
    \eeas
    Using an  analogous decomposition and estimation to \eqref{eq:xvx055}, similar to \eqref{eq:xvx117}, we can deduce from
   \eqref{eq:xvx119} that
    \beas
  && \hspace{-1.5cm}  \Big| Y^{\Th^m_\mu}_{t} \big({\t_\mu} ,        Y^{\Th^m_\mu}_{\t_\mu} \big(\t^m_\mu,         \eta^m_\mu     \big)     \big)
    -  Y^{\Th^m_\mu}_{t} \big(\t_\mu,         \eta^m_\mu     \big)   \Big|^p     \le E \bigg[      \Big|   Y^{\Th^m_\mu}_{\t_\mu} (\t^m_\mu, \eta^m_\mu )   -\eta^m_\mu     \Big|^p  \Big| \cF_t  \bigg]\\
  &&        \le c_0 E \bigg[      \int_{\t_\mu}^{\t^m_\mu} \neg \Big(  \Big|X^{\Th^m_\mu}_{\t^m_\mu \land s}
    - X^{\Th_\mu}_{{\t_\mu} \land s} \Big|^2 + \big| X^{\Th_\mu}_{{\t_\mu} \land s}  \big|^2
    +  | \eta^m_\mu  |^p +  c_\k  \Big) ds  \Big| \cF_t  \bigg]
         \le   \frac{\wt{C}(\k,x,\d)}{m}   , \q \pas
     \eeas

 Applying  \eqref{eq:p677} with   $(\z, \t ,\eta ) = \big( \t_\mu,  \t^m_\mu,         \eta^m_\mu  \big)$,
           applying  \eqref{eq:q104} with $\eta = \eta^m_\mu$ and using \eqref{eq:j335}   yield that \pas
   \bea
 && \hspace{-1cm}   Y^{\Th^m_\mu}_{t}  \big( {\t^m_\mu},          \eta^m_\mu      \big)
      \neg =  \neg   Y^{\Th^m_\mu}_{t}  \Big( \t_\mu, Y^{\Th^m_\mu}_{\t_\mu}
     \big( {\t^m_\mu},    \eta^m_\mu  \big)  \Big)
      \le   \neg   Y^{\Th^m_\mu}_{t} \big({\t_\mu},         \eta^m_\mu     \big) + \frac{\wt{C}(\k,x,\d)}{m^{1/p}}
   \neg  = \neg   Y^{\Th_\mu}_{t} \big({\t_\mu},         \eta^m_\mu     \big) + \frac{\wt{C}(\k,x,\d)}{m^{1/p}} \nonumber  \\
    &&   \le   \neg   Y^{\Th_\mu}_{t}  \big({\t_\mu}, \wt{\f}\big({\t_\mu} , X^{\Th_\mu}_{{\t_\mu} } \big) \big)
        \neg +  \neg \frac{\wt{C}(\k,x,\d)}{m^{1/p}}
      \neg   \le   \neg   \esup{\mu \in \cU_t} \;  Y^{t,x,\mu, \beta (  \mu ) }_t \Big(\t_{\beta,\mu},
    \wt{\f}\big(\t_{\beta,\mu},  X^{t,x,\mu, \beta (  \mu ) }_{\t_{\beta,\mu}} \big) \Big)
    \neg  +  \neg   \frac{\wt{C}(\k,x,\d)}{m^{1/p}}   .    \qq       \label{eq:h041b}
    \eea
       As       $  \beta_m    ( \mu)   \neg = \neg    \wh{\beta}   ( \mu )        $, $ds \times dP -$a.s.
    on $\[t,\t^m_\mu \[$,   applying
   \eqref{eq:p611} with $(\t, A) = (\t^m_\mu, \es)$ yields that \pas
   \bea    \label{eq:q101c}
    X^{\Th^m_\mu  }_s = X^{\wh{\Th}_\mu}_s     , \q \fa s \in [t,\t^m_\mu] .
   \eea
     Given $i = 1, \cds, N_m$, \eqref{eq:q101c} shows that $ X^{\Th^m_\mu  }_{t_i} = X^{\wh{\Th}_\mu}_{t_i} $, \pas ~ on $\cA^{\mu,m}_i \backslash \cA^m_i$.
    Since     $  \big[\beta_m( \mu )\big]^{t_i}_r(\o)  \neg = \neg  \big( \beta_m( \mu )\big)_r(\o)
  \neg = \neg   \big(  \wh{\beta}  ( \mu   ) \big)_r(\o)    \neg = \neg     \big[\wh{\beta}  ( \mu   )\big]^{t_i}_r(\o)
     $ holds  $ds \times dP-$a.s.  on $       \[\t^m_\mu,T\]_{\cA^{\mu,m}_i \backslash \cA^m_i}  \neg = \neg    [t_i, T]  \neg \times \neg  (\cA^{\mu,m}_i \backslash \cA^m_i)  $. Then by   \eqref{eq:xxa273} and a similar argument to \eqref{eq:xxd213b},
 it holds \pas ~ on $\cA^{\mu,m}_i \backslash \cA^m_i$ that
   \bea   \label{eq:xvx073b}
   Y^{\Th^m_\mu}_{{\t^m_\mu}} \big(T, g \Big( X^{ \Th^m_\mu }_T \Big) \big)
   = Y^{\Th^m_\mu}_{t_i} \big(T, g \Big( X^{ \Th^m_\mu }_T \Big) \big)
   = J \big(  \Th^{m,t_i}_\mu \big)
    =       J \big( \wh{\Th}_\mu^{t_i} \big)
     =       Y^{\wh{\Th}_\mu}_{t_i} \big(T, g \big( X^{ \wh{\Th}_\mu }_T \big) \big)
   =       Y^{\wh{\Th}_\mu}_{{\t^m_\mu}} \big(T, g \big( X^{ \wh{\Th}_\mu }_T \big) \big) ,
   \eea
    where $ \wh{\Th}_\mu^{t_i}    \neg \dfnn \neg \Big(t_i, X^{\wh{\Th}_\mu}_{t_i}, [ \mu ]^{t_i},
 \big[\wh{\beta}  ( \mu   ) \big]^{t_i}  \Big) $.

  \ss   Given $A \in \cF_t$, one can deduce  that
  \beas
   \b1_A  X^{\wh{\Th}_\mu}_{{\t_\mu} \vee s}
 & \tneg \dneg = & \tneg \dneg  \b1_A  X^{\wh{\Th}_\mu}_{{\t_\mu} }  \neg +  \neg  \b1_A  \neg   \int_{\t_\mu}^{{\t_\mu} \vee s}  \neg   b \Big(r, X^{\wh{\Th}_\mu}_r,\mu_r, (\wh{\beta}(\mu))_r  \Big) dr  \neg + \neg
     \b1_A  \neg  \int_{\t_\mu}^{{\t_\mu} \vee s}  \neg    \si \Big(r, X^{\wh{\Th}_\mu}_r,\mu_r, (\wh{\beta}(\mu))_r  \Big) d B_r \\
  & \tneg \dneg    = & \tneg \dneg   \b1_A  X^{ \wh{\Th}_\mu }_{{\t_\mu} }  \neg + \neg \int_t^s    \b1_{\{r \ge {\t_\mu}\}} \b1_A  b \Big(r, X^{\wh{\Th}_\mu}_r, \mu_r, \psi(r,  \mu_r)  \Big) dr  \neg + \neg
   \int_t^s \b1_{\{r \ge {\t_\mu}\}} \b1_A   \si \Big(r, X^{\wh{\Th}_\mu}_r, \mu_r, \psi(r,  \mu_r)  \Big) d B_r ,   \q s \in \neg [t, T].
 \eeas
  It then follows from \eqref{eq:xxa614} that
  \beas
   \b1_A  \underset{r \in [t,s]}{\sup}  \Big| X^{\wh{\Th}_\mu}_{{\t_\mu} \vee r} \Big| & \tneg \le & \tneg
     \b1_A  (|x|  \neg + \neg  \d )   \neg + \dneg
    \int_t^s  \neg  \b1_{\{r \ge {\t_\mu}\}}  \b1_A  \Big| b \Big(r, X^{\wh{\Th}_\mu}_r, \mu_r,  \psi(r,  \mu_r)  \Big) \Big| dr \\
    & \tneg  & \tneg
 + \underset{r \in [t,s]}{\sup} \bigg| \int_t^r  \neg  \b1_{\{r' \ge {\t_\mu}\}}   \b1_A     \si \Big(r', X^{\wh{\Th}_\mu}_{r'},\mu_{r'}, \psi(r',  \mu_{r'})  \Big) d B_{r'} \bigg| ,   \q    s  \neg \in \neg  [t, T].
 \eeas
  Using an  analogous decomposition and estimation to \eqref{eq:xvx055},
  one can deduce from  H\"older's inequality,  Doob's martingale inequality, \eqref{b_si_linear_growth}, \eqref{b_si_Lip},
  \eqref{eq:xxa614} and Fubini's Theorem  that
 \beas
 && \hspace{-1cm}  E \bigg[  \b1_A  \underset{r \in [t,s]}{\sup}  \Big|X^{\wh{\Th}_\mu}_{{\t_\mu} \vee r}\Big|^2 \bigg]
      \neg    \le     \neg
      \wt{C}(\k,x,\d) P(A) \neg  + \neg
  c_0 E \neg \int_t^s  \neg  \b1_{\{r \ge {\t_\mu}\}}  \b1_A  \Big( \big| b \big(r, X^{\wh{\Th}_\mu}_r, \mu_r, \psi(r,  \mu_r)\big) \big|^2     \neg + \neg
     \big| \si \big(r, X^{\wh{\Th}_\mu}_r, \mu_r,  \psi(r,  \mu_r)  \big) \big|^2 \Big) dr  \\
 &  & \le        \wt{C}(\k,x,\d) P(A)  \neg  +  \neg
  c_0 E        \int_t^s  \neg  \b1_{\{r \ge {\t_\mu}\}} \b1_A \Big| X^{\wh{\Th}_\mu}_{{\t_\mu} \vee r} \Big|^2 \neg  dr       \neg  \le   \neg
    \wt{C}(\k,x,\d)  P(A)    \neg  + \neg
  c_0   \neg  \int_t^s  \neg  E \bigg[  \b1_A    \underset{r' \in [t,r]}{\sup}  \Big|X^{\wh{\Th}_\mu}_{{\t_\mu} \vee r'}\Big|^2    \bigg]  dr ,   ~     \fa  s  \neg \in \neg  [t, T] .
 \eeas
   Then an application of  Gronwall's inequality shows that
  $  E \bigg[  \b1_A  \underset{r \in [t,s]}{\sup}  \Big|X^{\wh{\Th}_\mu}_{{\t_\mu} \vee r}\Big|^2 \bigg]  \le \wt{C}(\k,x,\d) P(A)  e^{c_0(s-t)}$, $ s \in [t,T] $.
  In particular,
  $  E \bigg[  \b1_A  \underset{r \in [\t_\mu,T]}{\sup}  \Big|X^{\wh{\Th}_\mu}_{  r} \Big|^2 \bigg]
  =  E \bigg[  \b1_A  \underset{r \in [t,T]}{\sup}  \Big|X^{\wh{\Th}_\mu}_{{\t_\mu} \vee r}\Big|^2 \bigg]
  \le \wt{C}(\k,x,\d) P(A)   $.
    Letting $A$ vary in $\cF_t$ yields that
 \bea    \label{eq:xvx085}
  E \bigg[     \underset{r \in [\t_\mu,T]}{\sup}  \Big|X^{\wh{\Th}_\mu}_{  r}\Big|^2 \bigg| \cF_t \bigg]  \le \wt{C}(\k,x,\d) , \q   \pas
  \eea

 Let $(\wh{Y}^\mu, \wh{Z}^\mu) \in \hG^p_\bF([t,T])$ be the unique solution of the following BSDE with zero generator:
   \beas
   \wh{Y}^\mu_s = Y^{\wh{\Th}_\mu}_{\t_\mu} \Big(T, g \big( X^{ \wh{\Th}_\mu }_T \big) \Big)
     - \int_s^T \wh{Z}^\mu_r d B_r      , \q    s \in [t,T] .
   \eeas
   \if{0}
      For any  $s \in [t,T]$, one can deduce that
   \beas
  \wh{Y}^\mu_{{\t_\mu} \land s} = E[\wh{Y}^\mu_{{\t_\mu} \land s} | \cF_{\t_\mu}]
  = E \bigg[  Y^{\wh{\Th}_\mu}_{\t_\mu} \Big(T, g \big( X^{ \wh{\Th}_\mu }_T \big) \Big)
    - \int_{{\t_\mu} \land s}^T \wh{Z}^\mu_r d B_r \Big| \cF_{\t_\mu} \bigg]
   = Y^{\wh{\Th}_\mu}_{{\t_\mu}} \big(T, g \big( X^{ \wh{\Th}_\mu }_T \big) \big)   - \int_{{{\t_\mu}} \land s}^{\t_\mu} \wh{Z}^\mu_r d B_r , \q  \pas ~
   \eeas
       By the continuity of process $\wh{Y}^\mu$,  it holds \pas ~ that
    \bea   \label{eq:xvx031a}
   \wh{Y}^\mu_{{{\t_\mu}} \land s}  = Y^{\wh{\Th}_\mu}_{{\t_\mu}} \big(T, g \big( X^{ \wh{\Th}_\mu }_T \big) \big)  - \int_{{{\t_\mu}} \land s}^{{\t_\mu}} \wh{Z}^\mu_r d B_r  = Y^{\wh{\Th}_\mu}_{{\t_\mu}} \big(T, g \big( X^{ \wh{\Th}_\mu }_T \big) \big)  - \int_s^T \b1_{\{r < {{\t_\mu}}\}} \wh{Z}^\mu_r d B_r     , \q    s \in [t,T] .
   \eea
   Thus, we see that  $(\wh{Y}^\mu_s, \wh{Z}^\mu_s) \neg = \neg  \big(  \wh{Y}^\mu_{{\t_\mu} \land s}, \b1_{\{s < {\t_\mu}\}} \wh{Z}^\mu_s \big),
   s  \neg \in \neg  [t,T]   $. Also,
   taking $[\cd|\cF_{{\t_\mu} \land s}]$ in \eqref{eq:xvx031a} shows that \pas
   \beas
    \wh{Y}^\mu_s  \neg = \neg   \wh{Y}^\mu_{{\t_\mu} \land s}  \neg = \neg  E\big[Y^{\wh{\Th}_\mu}_{\t_\mu} \big(T, g \big( X^{ \wh{\Th}_\mu }_T \big) \big) \big|\cF_{{\t_\mu} \land s}\big] ,
    \q   \fa  s  \neg \in \neg  [t,T] .
   \eeas
   \fi
 Analogous to \eqref{eq:xvx034},  $( \wh{\cY}^\mu,  \wh{\cZ}^\mu )
  \neg \dfnn  \neg \Big\{ \Big( \b1_{\{s < {\t_\mu}\}} \wh{Y}^\mu_s   +    \b1_{\{s  \ge  {\t_\mu}\}} Y^{\wh{\Th}_\mu}_s  \Big(T, g \big( X^{ \wh{\Th}_\mu }_T \big) \Big),
    \b1_{\{s < {\t_\mu}\}} \wh{Z}^\mu_s    +    \b1_{\{s  \ge  {\t_\mu}\}} Z^{\wh{\Th}_\mu}_s  \Big(T, g \big( X^{ \wh{\Th}_\mu }_T \big) \Big)     \Big) \Big\}_{s \in [t,T]} \\ \in \hG^p_\bF([t,T]) $
   solves the following BSDE
   \if{0}
    \bea
   \wh{\cY}^\mu_s & \tneg \dneg =& \tneg \dneg   \b1_{\{s  \ge  {\t_\mu}\}} Y^{\wh{\Th}_\mu}_s  \big(T, g \big( X^{ \wh{\Th}_\mu }_T \big) \big)
   \neg + \neg   \b1_{\{s < {\t_\mu}\}}Y^{\wh{\Th}_\mu}_{\t_\mu} \big(T, g \big( X^{ \wh{\Th}_\mu }_T \big) \big)   \neg - \neg  \b1_{\{s < {\t_\mu}\}} \int_s^T  \neg  \wh{Z}^\mu_r d B_r
   \neg  = \neg  Y^{\wh{\Th}_\mu}_{{\t_\mu} \vee s}  \big(T, g \big( X^{ \wh{\Th}_\mu }_T \big) \big)  \neg - \neg \b1_{\{s < {\t_\mu}\}} \int_s^T \neg \b1_{\{r < {\t_\mu}\}} \wh{Z}^\mu_r    d B_r  \nonumber \\
   & \tneg  \dneg =& \tneg \dneg   g \big( X^{ \wh{\Th}_\mu }_T \big)  \neg + \neg  \int_{{\t_\mu} \vee s}^T f^{\wh{\Th}_\mu}_T \big(r,Y^{\wh{\Th}_\mu}_r  \big(T, g \big( X^{ \wh{\Th}_\mu }_T \big) \big) ,  Z^{\wh{\Th}_\mu}_r  \big(T, g \big( X^{ \wh{\Th}_\mu }_T \big) \big) \big) dr
     \neg - \neg  \int_{{\t_\mu} \vee s}^T  Z^{\wh{\Th}_\mu}_r  \big(T, g \big( X^{ \wh{\Th}_\mu }_T \big) \big)  dB_r  \neg - \neg  \int_s^T \b1_{\{r < {\t_\mu}\}} \wh{Z}^\mu_r    d B_r \nonumber   \\
   & \tneg \dneg  =& \tneg \dneg   g \big( X^{ \wh{\Th}_\mu }_T \big)  \neg + \neg  \int_s^T \b1_{\{r  \ge  {\t_\mu} \}} f^{\wh{\Th}_\mu}_T \big(r,\wh{\cY}^\mu_r , \wh{\cZ}^\mu_r \big) dr     \neg - \neg  \int_s^T    \wh{\cZ}^\mu_r   dB_r     , \q s \in [t,T] .
   \eea
   \fi
   \beas
   \wh{\cY}^\mu_s   =      g \Big( X^{ \wh{\Th}_\mu }_T \Big)  \neg + \neg  \int_s^T \b1_{\{r  \ge  {\t_\mu} \}} f^{\wh{\Th}_\mu}_T \big(r,\wh{\cY}^\mu_r , \wh{\cZ}^\mu_r \big) dr     \neg - \neg  \int_s^T    \wh{\cZ}^\mu_r   dB_r     , \q s \in [0,T] .
   \eeas
    Then \eqref{eq:s031}, \eqref{eq:xvx131} and H\"older's inequality imply that \pas
   \beas
     E \Big[     \underset{s \in [\t_\mu,T]}{\sup} \Big| Y^{\wh{\Th}_\mu}_s \Big(T, g \Big( X^{ \wh{\Th}_\mu }_T \Big) \Big)\Big|^p    \Big|\cF_t\Big]
       & \tneg \le & \tneg   E \Big[ \, \underset{s \in [t,T]}{\sup}  \big|    \wh{\cY}^\mu_s     \big|^p  \Big|\cF_t  \Big]
     \le   c_0 E \bigg[   \Big|  g \Big( X^{ \wh{\Th}_\mu }_T \Big) \Big|^p  \neg + \neg \int_{\t_\mu}^T \neg  \Big| f^{\wh{\Th}_\mu}_T \big(s,  0 ,    0   \big)  \Big|^p ds   \bigg|\cF_t   \bigg] \\
   & \tneg   = & \tneg  c_0 E \bigg[   \Big|  g \Big( X^{ \wh{\Th}_\mu }_T \Big) \Big|^p  \neg + \neg  \int_{\t_\mu}^T \neg  \Big| f  \big(s, X^{\wh{\Th}_\mu}_s, 0 ,
    0 , \mu_s, \psi(s,   \mu_s )  \big)  \Big|^p ds  \Big| \cF_t  \bigg]     .
     \eeas
  Using an  analogous decomposition and estimation to \eqref{eq:xvx055}, we can then deduce from  \eqref{f_linear_growth},
     \eqref{f_Lip}    and \eqref{eq:xvx085}      that
   \beas
        E \bigg[     \underset{s \in [\t_\mu,T]}{\sup} \Big| Y^{\wh{\Th}_\mu}_s \Big(T, g \Big( X^{ \wh{\Th}_\mu }_T \Big) \Big)\Big|^p    \bigg| \cF_t \bigg]
        \le      c_\k +  c_0  E \bigg[      \underset{s \in [{\t_\mu} ,T]}{\sup} \Big|X^{\wh{\Th}_\mu}_s   \Big|^2     \bigg| \cF_t  \bigg]
  \le \wt{C}(\k,x,\d)    , \q \pas
     \eeas

  \ss   Let $ \wh{\eta}^m_\mu  \neg \dfnn \neg  Y^{\Th^m_\mu}_{{\t^m_\mu}} \Big(T, g \Big( X^{ \Th^m_\mu }_T \Big) \Big)
   \neg \vee   \wt{\eta}^m_\mu    \neg \in \neg  \hL^p \big( \cF_{\t^m_\mu} \big)  $ and set $ \wt{\cA}^m_\mu \dfnn \Big\{ Y^{ \Th^m_\mu}_{\t^m_\mu} \Big(T, g \big(X^{ \Th^m_\mu}_T \big) \Big)  >
      \wt{\eta}^m_\mu \Big\} \in \cF_{\t^m_\mu} $. Clearly,  $\b1_{\wt{\cA}^m_\mu} \le    \b1_{\cA^m_\mu} $, \pas ~
   Applying \eqref{eq:xvx071} with $\wt{p} =\frac{1+p}{2}$, we can deduce from H\"older's inequality, \eqref{eq:xvx119}
    and \eqref{eq:xvx073b} that
    \bea
     && \hspace{-0.8cm}        \Big|  Y^{\Th^m_\mu}_{t}  ({\t^m_\mu},  \wh{\eta}^m_\mu  )
   \neg  - \neg  Y^{\Th^m_\mu}_{t}  \big({\t^m_\mu},      \eta^m_\mu   \big)  \Big|^{\wt{p}}
      \neg \le   \neg      c_0    E \Big[ \big|   \wh{\eta}^m_\mu      \neg  -  \neg
                \eta^m_\mu      \big|^{\wt{p}}   \Big|\cF_t\Big]
                   \neg   =   \neg    c_0    E \bigg[ \b1_{ (\wt{\cA}^m_\mu)^c   } \big|  \wt{\eta}^m_\mu       \neg  -  \neg
                \eta^m_\mu     \big|^{\wt{p}}  \neg + \neg  \b1_{ \wt{\cA}^m_\mu  } \Big|   Y^{\Th^m_\mu}_{{\t^m_\mu}} \Big(T, g \Big( X^{ \Th^m_\mu }_T \Big) \Big)      \neg  -  \neg
                \eta^m_\mu     \Big|^{\wt{p}}   \bigg|\cF_t\bigg]          \nonumber    \\
         & &  \le   \neg     c_0    E \Big[   \big|  X^{\Th^m_\mu}_{\t^m_\mu}       \dneg  -  \dneg
                 X^{\Th_\mu}_{\t_\mu}      \big|^{\frac{2 \wt{p}}{p}} \Big] \neg + \neg c_0
                  \bigg\{ E \Big[   \b1_{ \wt{\cA}^m_\mu  }     \Big|\cF_t\Big] \bigg\}^{\frac{p-\wt{p}}{p}}
                   \bigg\{ E \bigg[     \Big| Y^{\wh{\Th}_\mu}_{{\t^m_\mu}} \Big(T, g \big( X^{ \wh{\Th}_\mu }_T \big) \Big)         \neg  -  \neg   \eta^m_\mu     \Big|^p   \bigg|\cF_t\bigg] \bigg\}^{\frac{\wt{p}}{p}}   \nonumber    \\
                    & &  \le   \neg     c_0  \bigg\{  E \Big[   \big|  X^{\Th^m_\mu}_{\t^m_\mu}       \dneg  -  \dneg
                 X^{\Th_\mu}_{\t_\mu}      \big|^2   \Big] \bigg\}^{\frac{\wt{p}}{p}} \neg + \neg c_0
                  \bigg\{ E \Big[   \b1_{ \cA^m_\mu  }     \Big|\cF_t\Big] \bigg\}^{\frac{p-\wt{p}}{p}}
                   \bigg\{ E \bigg[     \underset{s \in [\t_\mu, T]}{\sup} \Big| Y^{\wh{\Th}_\mu}_s \Big(T, g \big( X^{ \wh{\Th}_\mu }_T \big) \Big)\Big|^p + \Big(C^{\wt{\f}}_{x,\d} + \frac{c_0}{m} \Big)^p     \bigg|\cF_t\bigg] \bigg\}^{\frac{\wt{p}}{p}}   \nonumber    \\
         & &     \le       \neg  \frac{\wt{C}(\k,x,\d)}{m^{\wt{p}/p}} \neg + \neg   \wt{C}(\k,x,\d)
    \bigg\{ E \Big[   \b1_{ \cup^{N_m}_{i=1} ( \cA^m_i)^c  }     \Big|\cF_t\Big] \bigg\}^{\frac{p-\wt{p}}{p}} , \q \pas
    \label{eq:xux711b}
        \eea

    Applying  \eqref{eq:p677} with $(\z, \t ,\eta ) \neg = \neg \Big({\t^m_\mu}, T,g \Big( X^{ \Th^m_\mu }_T \Big) \Big)$,
     we see from Proposition \ref{prop_BSDE_estimate_comparison} (2),  \eqref{eq:xux711b} and \eqref{eq:h041b}    that \pas
            \bea
  && \hspace{-1cm}     Y^{\Th^m_\mu}_{t}  \neg \Big(T,g \Big( X^{ \Th^m_\mu }_T \Big) \Big)
   \neg  =    \neg  Y^{\Th^m_\mu}_{t}  \neg \Big({\t^m_\mu},
    Y^{\Th^m_\mu}_{{\t^m_\mu}} \Big(T, g \Big( X^{ \Th^m_\mu }_T \Big) \Big)   \Big)
    \neg  \le \neg  Y^{\Th^m_\mu}_{t}  \neg  \big( {\t^m_\mu},         \wh{\eta}^m_\mu      \big)    \nonumber \\
    & &   \le   \neg  Y^{\Th^m_\mu}_{t}  \neg  \big( {\t^m_\mu},          \eta^m_\mu      \big) \neg + \neg   \frac{\wt{C}(\k,x,\d)}{m^{1/p}} \neg + \neg   \wt{C}(\k,x,\d)  \bigg\{   E \Big[ \b1_{\cup^{N_m}_{i=1} ( \cA^m_i)^c  }     \Big|\cF_t\Big]  \bigg\}^{\frac{p-\wt{p}}{p \wt{p}}}   \nonumber \\
  & &   \le   \neg   \esup{\mu \in \cU_t} \,  Y^{t,x,\mu, \beta (  \mu ) }_t \Big(\t_{\beta,\mu},
    \wt{\f}\big(\t_{\beta,\mu},  X^{t,x,\mu, \beta (  \mu ) }_{\t_{\beta,\mu}} \big) \Big)
     \neg + \neg   \frac{\wt{C}(\k,x,\d)}{m^{1/p}} \neg + \neg   \wt{C}(\k,x,\d)  \bigg\{   E \Big[ \b1_{\cup^{N_m}_{i=1} ( \cA^m_i)^c  }     \Big|\cF_t\Big]  \bigg\}^{\frac{p-\wt{p}}{p \wt{p}}}   . \label{eq:xux621}
    \eea

  Letting $      \wh{\cA}_m \dfnn \Big\{     E \Big[ \b1_{\cup^{N_m}_{i=1} ( \cA^m_i)^c  }     \Big|\cF_t\Big] >  m^{\frac{1+p}{1-p}}   \Big\} $, one can deduce that

                \beas
 P( \wh{\cA}_m  )
\le m^{\frac{1+p}{p-1}} E \Big[   E \Big[ \b1_{\cup^{N_m}_{i=1} ( \cA^m_i)^c  }     \Big|\cF_t\Big] \Big]
 =  m^{\frac{1+p}{p-1}} P \big(           \cup^{N_m}_{i=1} ( \cA^m_i)^c           \big)
 \le m^{\frac{1+p}{p-1}} \sum^{N_m}_{i=1} P \big(            ( \cA^m_i)^c           \big) \le m^{-p} .
 \eeas
       Multiplying $ \b1_{\wh{\cA}_m^c} $ to both sides of  \eqref{eq:xux621}     yields that
     \beas
     \b1_{\wh{\cA}_m^c} J \big( t , x , \mu ,  \beta_m (\mu)  \big)
     \le \b1_{\wh{\cA}_m^c}  \esup{\mu \in \cU_t} \;  Y^{t,x,\mu, \beta (  \mu ) }_t \Big(\t_{\beta,\mu},
    \wt{\f}\big(\t_{\beta,\mu},  X^{t,x,\mu, \beta (  \mu ) }_{\t_{\beta,\mu}} \big) \Big)
      \neg  +  \neg    \frac{\wt{C}(\k,x,\d)}{m^{1/p}}  ,  \q   \pas
     \eeas
   Since $\wh{\cA}_m $ does not depend on $ \mu $ nor on   $ \beta  $,
    taking essential supremum over $\mu \in \cU_t$ and applying Lemma \ref{lem_ess} (2) yield   that
      \beas
    \b1_{\wh{\cA}_m^c}  w_1(t,x)      \le  \b1_{\wh{\cA}_m^c}   I \big(t,x,   \beta_m  \big)
     \le       \b1_{\wh{\cA}_m^c}  \esup{\mu \in \cU_t} \;  Y^{t,x,\mu, \beta (  \mu ) }_t \Big(\t_{\beta,\mu},
    \wt{\f}\big(\t_{\beta,\mu},  X^{t,x,\mu, \beta (  \mu ) }_{\t_{\beta,\mu}} \big) \Big)
        \neg + \neg   \frac{\wt{C}(\k,x,\d)}{m^{1/p}}     ,  \q   \pas
      \eeas
        Then taking essential infimum over $\beta \in \fB_t$ and  using Lemma \ref{lem_ess} (2) again, we obtain
           \bea \label{eq:xxd231c}
    \b1_{\wh{\cA}_m^c}  w_1(t,x)
     \le       \b1_{\wh{\cA}_m^c}  \underset{\beta \in \fB_t }{\essinf} \; \esup{\mu \in \cU_t} \;  Y^{t,x,\mu, \beta (  \mu ) }_t \Big(\t_{\beta,\mu},
    \wt{\f}\big(\t_{\beta,\mu},  X^{t,x,\mu, \beta (  \mu ) }_{\t_{\beta,\mu}} \big) \Big)
        \neg + \neg   \frac{\wt{C}(\k,x,\d)}{m^{1/p}}     ,  \q   \pas
      \eea
      As  $ \dis   \sum_{m \in \hN} P \big( \wh{\cA}_m \big) \le \sum_{m \in \hN} m^{-p} < \infty $,
      similar to \eqref{eq:xux523},
       Borel-Cantelli theorem implies that      $    \lmt{m \to \infty} \b1_{\wh{\cA}_m} = 0$, \pas ~ \,
      Thus, letting $m \to \infty$ in \eqref{eq:xxd231c} yields \eqref{eq:xqxqx223}. \qed

    \ss \no   {\bf 2)}    For any $ ( t,x ,y,z, u,v  ) \in [0,T] \times \hR^k \times \hR  \times \hR^d \times \hU \times \hV$,  we   define
 \beas
  \fg ( x )    \dfnn    - g (x )
  \q \hb{and} \q \ff  ( t,x ,y,z, u,v  )
 \dfnn  -f  ( t,x , -y,-z, u,v  )   .
  \eeas
  Given $ (\mu, \nu) \in \cU_t \times \cV_t$, we  let $\Th$ stand for  $(t,x,\mu,\nu)$.
  For any   $\t \in \cS_{t,T}$ and any $\eta \in \hL^p(\cF_\t)$,    let
     $ \big(\cY^\Th (\t, \eta),\cZ^\Th (\t, \eta)  \big) $ denote the unique solution of the BSDE$\big(t,\eta,  \ff^\Th_\t    \big)$ in $ \hG^q_\bF  \big([t,T]\big) $, where
    \beas
     \ff^\Th_\t   (s,\o,y,z ) \dfnn \b1_{\{s < \t(\o)\}}  \ff  \Big(s,  X^\Th_s(\o) ,  y,  z, \mu_s (\o),\nu_s (\o) \Big) ,
    \q   \fa (s,\o,y,z  ) \in   [t,T] \times \O  \times \hR  \times \hR^d .
    \eeas
 Multiplying $-1$ in the BSDE$\big(t,\eta,  \ff^\Th_\t \big)$ shows that $\big(\neg - \neg  \cY^\Th (\t, \eta),$
  $ -\cZ^\Th (\t, \eta)   \big) \neg \in \neg \hG^q_\bF  \big([t,T]\big) $ solves the BSDE$\big(t,-\eta,  f^\Th_\t \big)$. To wit
   \bea  \label{eq:s304}
    \big(-\cY^\Th (\t, \eta),-\cZ^\Th (\t, \eta)   \big) = \big(Y^\Th (\t, -\eta), Z^\Th (\t, -\eta)   \big) .
   \eea

 Given $  (t,x ) \in [0,T] \times \hR^k$, let us consider the situation where    player II acts first by choosing
 a $\cV_t-$control to maximize
  $\cY^{t,x,  \a (\nu), \nu }_t   \Big(T,  \fg \Big(X^{t,x,\a (\nu), \nu}_T \Big)  \Big)$, where
  $\a \neg \in \neg \fA_t$ is player I's  strategic response.
  The corresponding priority value   of player II is
        $   \fw_2  (t,x   ) \dfnn \underset{\a \in \fA_t }{\essinf} \; \underset{\nu \in \cV_t}{\esssup} \; \cY^{t,x,  \a (\nu), \nu }_t   \Big(T,\fg \Big(X^{t,x,\a (\nu), \nu}_T \Big)  \Big)     $.
   We see    from   \eqref{eq:s304}  that
   $$
  ~    - \fw_2  (t,x )    \neg =  \neg   \underset{\a \in \fA_t }{\esssup} \; \underset{\nu \in \cV_t}{\essinf}  - \cY^{t,x,  \a (\nu), \nu }_t   \Big(T,\fg \Big(X^{t,x,\a (\nu), \nu}_T \Big)  \Big)
      \neg   =  \neg    \underset{\a \in \fA_t }{\esssup} \; \underset{\nu \in \cV_t}{\essinf} \;    Y^{t,x,  \a (\nu), \nu }_t   \Big(T, g \Big(X^{t,x,\a (\nu), \nu}_T \Big)  \Big)  \neg = \neg  w_2  (t,x ) .
    $$

   Let $t \in (0,T] $ and  let $\f , \wt{\f } \neg : [t,T]     \times
     \hR^k  \neg \to \neg  \hR $   be two continuous functions   satisfying  $     \f  (s,x) \le w_2 (s,x)
      \neg \le \neg  \wt{\f } (s,x)      $, $ (s,x)
        \neg \in \neg [t,T] \neg \times \neg \hR^k$. As
            $   - \wt{\f }(s,x) \le \fw_2 (s,x) \le - \f  (s,x)   $, $ (s,x)
        \neg \in \neg [t,T] \neg \times \neg \hR^k$,
   applying the weak dynamic programming principle of part (1) yields that for any $x \in \hR^k$ and $\d \in (0,T-t]$
   \beas
 \qq  && \hspace{-2.5cm} \underset{\a \in \fA_t}{\essinf} \; \underset{\nu \in \cV_t}{\esssup}    \;
    \cY^{t,x,\a (\nu), \nu }_t \Big(\t_{\a,\nu},
    - \wt{\f } \big(\t_{\a,\nu},  X^{t,x,\a (\nu), \nu}_{\t_{\a,\nu}} \big) \Big) \\
    &&  \le \fw_2 (t,x) \le
    \underset{\a \in \fA_t}{\essinf} \; \underset{\nu \in \cV_t}{\esssup}    \;
    \cY^{t,x,\a (\nu), \nu }_t \Big(\t_{\a,\nu},
    -  \f  \big(\t_{\a,\nu},  X^{t,x,\a (\nu), \nu}_{\t_{\a,\nu}} \big) \Big) , \q \pas
   \eeas
        Multiplying $-1$ above  and using \eqref{eq:s304}, we obtain
        the weak dynamic programming principle for $w_2$.          \qed

     \subsection{Proofs of Section \ref{sec:PDE}}

  \label{subsection:Proof_S3}

    We will prove that $\ul{w}_i$ and $\ol{w}_i$, $i=1,2$ are
 viscosity  solutions of   \eqref{eq:PDE} in a standard way: Assume oppositely that
 the corresponding inequality of  \eqref{eq:PDE} does not hold for some test function $\vf$.
 We decompose $\ul{H}_i$ or $\ol{H}_i $ with $\vf$ in  the reverse inequality until
  we reach a similar reverse inequality satisfied by
   a control $\wh{\mu}$   or a strategy $\wh{\beta}$. Then applying the comparison
  result of BSDE, i.e. Proposition  \ref{prop_BSDE_estimate_comparison} (2), to such an inequality
  leads to a contradiction to the weak dynamic programming principle.

  \ss \no {\bf Proof of Theorem \ref{thm_viscosity}:}
  We only need to prove for $\ul{w}_1$ and $\ol{w}_1$, then  the results of $\ol{w}_2$ and $\ul{w}_2$ follow   by
  a similar transformation  to that used in the   proof of Theorem \ref{thm_DPP}, part (2).

   \ss \no {\bf a)}   We first    show that $\ul{w}_1$  is a viscosity supersolution of   \eqref{eq:PDE} with Hamiltonian
 $\ul{H}_1$.
   Let $(t_0,x_0, \vf) \in (0,T) \times \hR^k \times \hC^{1,2}\big([0,T] \times \hR^k\big)$ be
   such that   $\ul{w}_1(t_0,x_0) = \vf (t_0,x_0)$
  and that  $\ul{w}_1-\vf$ attains a strict local minimum  at $(t_0,x_0)$, i.e., for some
     $\d_0 \in  \big(0,    t_0 \land  ( T-t_0 ) \big) $
   \bea     \label{eq:a071}
     (\ul{w}_1 - \vf) (t,x) >  (\ul{w}_1 - \vf) (t_0,x_0) = 0,
   \q \fa (t,x) \in O_{\d_0} (t_0,x_0) \big\backslash \big\{ (t_0,x_0) \big\} .
   \eea

    We simply  denote $\big(\vf(t_0   ,x_0),
     D_x \vf    (t_0  ,x_0),  D^2_x \vf    (t_0  ,x_0)\big) $ by $(y_0,z_0,\G_0 )$.  If  $\ul{H}_1 \big(t_0,x_0, y_0,z_0,\G_0\big) = - \infty$,
    then
  \beas           
     &&       \neg  -   \frac{\pa }{\pa t} \vf (t_0,x_0)
      \neg - \neg  \ul{H}_1 \big(t_0,x_0, y_0,z_0,\G_0\big)
      \ge  0    
     \eeas
     holds automatically.
   To make a contradiction,  we   assume    that when $\ul{H}_1 \big(t_0,x_0, y_0,z_0,\G_0\big) >    - \infty $,
  \bea  \label{eq:xvx315}
 \varrho \dfnn      \frac{\pa }{\pa t} \vf (t_0,x_0)
      \neg + \neg  \ul{H}_1 \big(t_0,x_0, y_0,z_0,\G_0\big)  > 0 .
 \eea

   For any $(t,x,y,z,\G, u,v ) \neg \in \neg  [0,T]  \neg \times \neg  \hR^k \neg \times  \neg  \hR
   \neg \times  \neg  \hR^d  \neg \times \neg \hS_k \neg \times \neg  \hU  \neg \times \neg  \hV   $, one can deduce from
   \eqref{b_si_linear_growth}$-$\eqref{f_Lip} that
 \bea
 && \hspace{-1cm} \big| H  (t,x,y,z,\G,u,v) \big|
   \neg \le  \neg  \frac14  |    \si \si^T(t,x,u,v)  |^2  \neg + \neg  \frac14 | \G |^2
   \neg  + \neg  \g | z | |b(t,x,u,v)|
     \neg + \neg  \g \Big(1  \neg + \neg |x|^{2/p}  \neg + \neg  |y|  \neg + \neg  \g |z||\si(t,x,u,v)|   \neg + \neg  [u]^{2/p}_{\overset{}{\hU}}  \neg + \neg  [v]^{2/p}_{\overset{}{\hV}} \Big) \nonumber  \\
  & & \le  \neg  \frac14 \g^2 \big( 1  \neg + \neg |x|  \neg + \neg [u]_{\overset{}{\hU}}
   \neg + \neg  [v]_{\overset{}{\hV}} \big)^2  \neg + \neg  \frac14 | \G |^2
     \neg + \neg  (\g \neg + \neg \g^2) | z | \big( 1  \neg + \neg |x|  \neg + \neg [u]_{\overset{}{\hU}}
      \neg + \neg  [v]_{\overset{}{\hV}} \big)
     \neg + \neg  \g \Big(1  \neg + \neg |x|^{2/p}  \neg + \neg  |y|   \neg + \neg  [u]^{2/p}_{\overset{}{\hU}}
      \neg + \neg  [v]^{2/p}_{\overset{}{\hV}} \Big) .     \label{eq:xvx311}
 \eea
 Set $C^0_\vf \dfnn |y_0|+|z_0|+|\G_0| = \big|\vf (t_0,x_0)\big| + \big| D_x \vf (t_0,x_0) \big| + \big| D^2_x \vf (t_0,x_0) \big| $,  and fix a $u_\sharp \in \pa O_\k (u_0) $.
  For any $u \notin O_\k(u_0)$,  we see from (A-u)  that  $ \psi(t_0, u)  \in \sO_u$,
   and it follows from  \eqref{eq:xvx311} that
 \bea
 && \hspace{-1cm} \underset{v \in \sO_u}{\inf}  \;  H (t_0,x_0, y_0,z_0,\G_0,u,v)
  \le \big| H (t_0,x_0, y_0,z_0,\G_0,u,\psi(t_0, u)) \big| \nonumber  \\
 && = \big| H (t_0,x_0, y_0,z_0,\G_0,u_\sharp,\psi(t_0, u_\sharp)) \big|
                \le \frac14 (C^0_\vf)^2 + C^0_\vf C(\k,x_0) + C(\k,x_0)  .   \label{eq:xvx313}
 \eea
 Here $C(\k,x_0)$ denotes a generic constant,  depending on $\k,    |x_0|  $,     $T$, $\g $,  $p$ and $|g(0)|$,
        whose form  may vary from line to line.

 \ss Similarly, it holds  for any $u \in O_\k(u_0)$ that
  \beas
  \underset{v \in \sO_u}{\inf}  \;  H (t_0,x_0, y_0,z_0,\G_0,u,v)
   \le    \big| H (t_0,x_0, y_0,z_0,\G_0,u,v_0) \big|
          \le    \frac14 (C^0_\vf)^2 + C^0_\vf C(\k,x_0) + C(\k,x_0) ,
 \eeas
 which together with \eqref{eq:xvx313} implies that
  \beas
     \ul{H}_1 \big(t_0,x_0, y_0,z_0,\G_0\big) \le     \underset{u \in \hU  }{\sup}  \;
     \underset{v \in \sO_u}{\inf}  \;  H (t_0,x_0, y_0,z_0,\G_0,u,v)
     \le \frac14 (C^0_\vf)^2 + C^0_\vf C(\k,x_0) + C(\k,x_0) < \infty   .
 \eeas
  Thus $\varrho < \infty$.

 \ss   As $\vf \in \hC^{1,2}\big([0,T] \times \hR^k\big) $, we see from \eqref{eq:xvx315} that for
  some   $\wh{u} \in \hU  $
 \beas
       \linf{(t,x)  \to  (t_0,x_0) }  \; \underset{v \in \sO_{\wh{u}} }{\inf}
 \;    H (t,x, \vf (t,x), D_x \vf (t,x),   D^2_x \vf (t,x),  \wh{u},v)
  \ge  \frac34 \varrho -  \frac{\pa }{\pa t} \vf (t_0,x_0)    .
 \eeas
 Moreover, there exists   a $ \d \in  (0, \d_0    )$ such that
 \bea  \label{eq:s423b}
    \underset{v \in \sO_{\wh{u}}}{\inf}
 \,  H (t,x, \vf (t,x), D_x \vf (t,x),   D^2_x \vf (t,x),  \wh{u},v)
 \neg  \ge  \neg  \frac12 \varrho  \neg - \neg   \frac{\pa }{\pa t} \vf (t,x), \q \fa (t,x)  \neg  \in  \neg  \ol{O}_{\d}  (t_0,x_0) .
 \eea

 Let  $     \wp    \dfnn  \inf \big\{ (  \ul{w}_1 \neg - \neg \vf    ) (t,x)  \neg : (t,x )  \neg \in \neg
     \ol{O}_\d  (t_0,x_0) \big\backslash O_{\frac{\d}{3}} (t_0,x_0) \big\}       $.
 Since the set $ \ol{O}_\d  (t_0,x_0) \big\backslash O_{\frac{\d}{3}} (t_0,x_0)$ is  compact, there exists
 a sequence $\{(t_n,x_n)\}_{n \in \hN} $ on $  \ol{O}_\d  (t_0,x_0) \big\backslash O_{\frac{\d}{3}} (t_0,x_0)$
 that  converges to some $(t_*,x_*) \in \ol{O}_\d  (t_0,x_0) \big\backslash O_{\frac{\d}{3}} (t_0,x_0) $
 and satisfies   $ \wp \neg = \neg  \lmtd{n \to \infty} (  \ul{w}_1 \neg - \neg \vf ) (t_n,x_n) $.
    The lower semicontinuity of $\ul{w}_1$ and the continuity of $\vf$ imply that
 $\ul{w}_1 - \vf$ is also lower semicontinuous.
 It follows that $ \wp \le ( \ul{w}_1 \neg - \neg \vf )(t_*,x_*) \le \lmtd { n \to \infty } ( \ul{w}_1 \neg - \neg \vf )(t_n,x_n ) = \wp  $,
 which together with \eqref{eq:a071} shows that
  \bea   \label{eq:s427}
  \wp  = \min \big\{ (\ul{w}_1 \neg - \neg \vf) (t,x): (t,x ) \in
     \ol{O}_\d  (t_0,x_0) \big\backslash O_{\frac{\d}{3}} (t_0,x_0) \big\} = (\ul{w}_1 \neg - \neg \vf)(t_*,x_*)  >  0 .
  \eea

  Then  we set     $\dis   \wt{\wp} \neg \dfnn  \neg   \frac{    \wp \land  \varrho}{ 2 (1 \vee \g) T }  \neg   > \neg  0$
  and let $ \big\{(t_j, x_j)\big\}_{j \in \hN}$ be a sequence of $O_{\frac{\d}{6}} (t_0,x_0) $ such   that
     \beas   
 \lmt{j \to \infty} (t_j,x_j) = (t_0,x_0)   \q \hb{and} \q  \lmt{j \to \infty} w_1 (t_j,x_j) = \ul{w}_1 (t_0,x_0) = \vf (t_0,x_0) = \lmt{j \to \infty}  \vf  (t_j,x_j) .
     \eeas
    So one can find a $j \in \hN$ such that
    \bea  \label{eq:s437}
  \big| w_1 (t_j,x_j)- \vf (t_j,x_j) \big| < \frac56 \wt{\wp}   t_0 .
   \eea

 Clearly,  $ \wh{\mu}_s \dfnn \wh{u}$, $s \in [t_j, T]$  is  a constant   $  \cU_{t_j}-$process.
  Fix $\beta \in \fB_{t_j}$.  We set $\Th \dfnn \big( t_j, x_j, \wh{\mu} ,    \beta (   \wh{\mu} )  \big)$
  and  define
 \beas
    \t = \t_{\beta,\wh{\mu}}  \dfnn  \inf \Big\{s  \neg \in \neg  ( t_j, T  ] \neg  : \big( s,X^{\Th}_s \big)
  \neg  \notin   \neg   O_{  \frac23 \d  } (t_j, x_j) \Big\}  \in \cS_{t_j,T}  .
 \eeas
    Since 
 $ \big| \big(T, X^{\Th}_T\big) \neg  - \neg (t_j,x_j) \big|  \neg  \ge  \neg  T \neg - \neg t_j
 \neg  \ge  \neg   T  \neg - \neg t_0  \neg - \neg  |t_j \neg - \neg t_0|  \neg > \neg  \d_0  \neg - \neg \frac{\d}{6}
  \neg >  \neg  \frac56 \d   \neg >  \neg  \frac23 \d   $,
  the continuity of $ X^{\Th} $ implies that \pas
   \bea
 && \qq \t  < T  ~ \;  \hb{and} ~ \; \big(\t \land  s , X^{\Th}_{\t \land  s}   \big)
 \in  \ol{O}_{  \frac23 \d }(t_j, x_j)
   \subset \ol{O}_{\frac56 \d}  (t_0, x_0) ,     \q \fa    s \in [t_j, T]  ;   \qq  \qq  \label{eq:a239b}   \\
 && \hb{in particular, }  ~
      \big( \t, X^{\Th}_{\t} \big)  \in  \pa O_{  \frac23 \d   }(t_j, x_j)
   \subset \ol{O}_{\frac56 \d}  (t_0,x_0) \big\backslash O_{\frac{\d}{2}} (t_0,x_0)  .  \label{eq:a041}
    \eea

    The continuity of $\vf$, $X^{\Th} $ and \eqref{eq:a239b} show that
  $\cY_s \dfnn \vf \big( \t \land  s, X^{\Th}_{\t \land  s} \big)
    + \wt{\wp} ( \t \land  s) $, $ s \in [t_j, T] $ defines a 
   bounded $\bF-$adapted  continuous process.  By   It\^o's formula,
    \bea
         \cY_s   & \tneg = &  \tneg   \cY_T     +   \int_s^T  \ff_r   dr
    -  \int_s^T  \cZ_r   d B_r ,   \q    s \in [t_j,T]  ,   \label{eq:s403b}
    \eea
where $\cZ_r  =        \b1_{\{r < \t\}}    D_x \vf    \big( r, X^{\Th}_r \big)
  \cd   \si \big( r, X^{\Th}_r, \wh{u},   \big(\beta (  \wh{\mu} )  \big)_r \big)$ and
 $$
     \ff_r   \neg   =   \neg    -  \b1_{\{ r  < \t  \}}
  \bigg\{ \wt{\wp}  \neg    +   \neg      \frac{\pa \vf }{\pa t} \big( r, X^{\Th}_r \big)
     \neg   +   \neg      D_x \vf    \big( r, X^{\Th}_r\big)   \cd    b \big( r, X^{\Th}_r \neg , \wh{u},
  (\beta (  \wh{\mu} )  )_r \big)
        \neg   +    \neg       \frac12  trace\Big( \si \si^T \big( r, X^{\Th}_r \neg , \wh{u},
  (\beta (  \wh{\mu} ) )_r \big) \cd    D^2_x \vf \big( r, X^{\Th}_r\big) \neg \Big) \neg \bigg\} .
     $$
 As $\vf  \neg \in  \neg   \hC^{1,2}\big([t,T]  \neg \times \neg  \hR^k\big)$,
  the measurability of $b$, $\si$, $X^{\Th}$, $\wh{u} $ and $\beta (  \wh{\mu} )$ implies   that  both $\cZ$ and $\ff$ are
     $\bF-$progressively    measurable.
  And one  can deduce from \eqref{b_si_linear_growth}, \eqref{b_si_Lip}, \eqref{eq:a239b} and H\"older's inequality  that
       \bea
 && \hspace{-1.5cm}  E \left[  \bigg( \int_{t_j}^T \neg |\cZ_s|^2   \,  ds  \bigg)^{ p /2 } \right]
    \le    (\g \wt{C}_\vf)^p \,  E \left[  \bigg( \int_{t_j}^{\t}
\Big( 1+ \big| X^{\Th}_s \big|  + [\wh{u}]_{\overset{}{\hU}}
  +   \big[ (\beta (  \wh{\mu} ) )_s \big]_{\overset{}{\hV}}   \Big)^2   ds  \bigg)^{ p /2 } \right] \nonumber \\
 &&    \le    c_0  \wt{C}^p_\vf  \bigg(   \big(1+ |x_0 | + \d + [\wh{u}]_{\overset{}{\hU}} \big)^p
+  \bigg\{ E    \int_{t_j}^T \neg \big[ (\beta (  \wh{\mu} ) )_s \big]_{\overset{}{\hV}}^2
  \,  ds    \bigg\}^{ p /2 } \bigg) < \infty , \q \hb{i.e. }
 \cZ \in \hH^{2,p} _{\bF  } \big( [t_j,T], \hR^d \big),  \qq \label{eq:s641}
  \eea
 where  $ \wt{C}_\vf \dfnn \underset{(t,x) \in \ol{O}_{\frac56 \d} (t_0,x_0) }{\sup} \big|D_x \vf (t,x)\big| < \infty $.
  Hence,
 $\big\{\big( \cY_s, \cZ_s \big) \big\}_{s \, \in  [t_j,T]}$
  solves the   BSDE$\big(  t_j  , \cY_T , \ff  \big)$.

 Let $\ell(x) \neg = \neg  c_\k  \neg + \neg  c_0 |x|^{2/p}$, $x  \neg \in \neg  \hR^k$ be the function appeared in Proposition \ref{prop_value_bounds}.
  Let $ \th_1  \neg : [0,T]  \neg \times \neg  \hR^k  \neg \to \neg  [0,1]$ be a continuous function    such that
    $\th_1  \neg \equiv \neg  0 $ on $ \ol{O}_{ \frac56 \d } (t_0,x_0) $ and  $\th_1  \neg \equiv \neg  1$
    on $ ( [0,T] \times \hR^k ) \big\backslash O_{ \d } (t_0,x_0)$.
  Also, let   $ \th_2   \neg  : [0,T]  \neg \times \neg  \hR^k  \neg \to \neg  [0,1]$ be another continuous function  such that
       $\th_2  \neg \equiv \neg 0 $ on $ \ol{O}_{\frac{\d}{3}} (t_0,x_0) $ and $\th_2  \neg \equiv \neg  1$
  on $ ([0,T]  \neg \times \neg  \hR^k) \big\backslash O_{ \frac{\d}{2} } (t_0,x_0)$.
    Define
    \bea   \label{eq:xvx151}
    \f (t,x) \dfnn - \th_1 (t,x) \ell (x) + \big(1-\th_1 (t,x)\big)
   \big(\vf(t,x) + \wp \th_2 (t,x) \big)        ,
   \q  \fa (t,x) \in [t_j,T] \times \hR^k ,
    \eea
    which is a continuous function satisfying  $ \f \le w_1$: given $(t,x) \in [t_j,T] \times \hR^k$,

 \no $\bullet $ if $(t,x) \in \ol{O}_{\frac{\d}{3}} (t_0,x_0) $,
  \eqref{eq:a071} shows that $\phi(t,x) = \vf (t,x)  \le \ul{w}_1 (t,x) \le w_1 (t,x)$;

 \no  $\bullet $ if $(t,x) \in O_\d  (t_0,x_0) \big\backslash \ol{O}_{\frac{\d}{3}} (t_0,x_0)   $,
 since $ \vf(t,x) \neg  + \neg   \wp \th_2 (t,x)  \neg  \le \neg   \vf (t,x) \neg  + \neg   \wp
  \neg  \le \neg    \ul{w}_1(t,x)  \neg  \le \neg   w_1(t,x) $
 by \eqref{eq:s427}, one can deduce   from Proposition   \ref{prop_value_bounds}  that
 $ \f (t,x) \le w_1(t,x) $;

 \no  $\bullet $ if $(t,x) \notin   O_{ \d } (t_0,x_0)   $,
 $ \f (t,x) = - \ell (x) \le w_1 (t, x) $.

   \ss 
 Then we can deduce from \eqref{eq:a041} that 
      \bea  \label{eq:xux811}
   \cY_T  =    \vf \big( \t, X^{\Th}_{\t} \big)    + \wt{\wp} T
 <  \vf \big( \t, X^{\Th}_{\t} \big) +\wp
     = \f \big( \t, X^{\Th}_{\t} \big) ,   \q   \pas
  \eea
 Since it holds $ds \times dP-$a.s. on $[t_j,T] \times \O$ that
  $
  \big[\big(\beta(\wh{\mu})\big)_s\big]_{\overset{}{\hV}} \le \k + C_\beta [\wh{\mu}_s]_{\overset{}{\hU}}
   = \k + C_\beta [\wh{u}]_{\overset{}{\hU}} \in \sO_{\wh{u}} $,
         \eqref{eq:a239b}, \eqref{eq:s423b}  and \eqref{f_Lip}  imply   that for $ds \times d P-$a.s.
  $ (s,\o)  \in [t_j,T] \times \O$
     \bea
     \ff_s (\o) & \tneg \dneg \le &  \tneg  \dneg   
 \b1_{\{ s \, < \t (\o) \}} \left\{ - \wt{\wp} - \frac12 \varrho   +    f \neg \left( s,\o, X^{\Th}_s(\o),
       \cY_s (\o)-\wt{\wp}s,    \cZ_s (\o), \wh{u},  ( \beta (  \wh{\mu} )  )_s (\o) \right)  \right\} \nonumber  \\
         &  \tneg \dneg  \le &  \tneg  \dneg  \b1_{\{ s \, < \t (\o) \}} \left\{ - \wt{\wp} - \frac12 \varrho  + \g \wt{\wp} T
          +    f \neg \left(s, \o, X^{\Th}_s (\o),
      \cY_s (\o),       \cZ_s (\o),\wh{u},  ( \beta (  \wh{\mu} )  )_s (\o) \right)  \right\}
  \neg  \le  \neg    f^{\Th }_{\t} \big(s , \o, \cY_s (\o),  \cZ_s (\o) \big)    . \q    \qq   \label{eq:xux813}
       \eea
   As   $f^{\Th}_{\t}$ is  Lipschitz continuous in $(y,z)$,
    Proposition   \ref{prop_BSDE_estimate_comparison} (2) implies that $P-$a.s.
   \beas
   \cY_s   \le Y^{\Th}_s \Big( \t,
   \f  \big( \t, X^{\Th}_{\t} \big)    \Big)  , \q  \fa  s \in [t_j,T] .
   \eeas
  Letting $s=t_j$ and using  the fact that $t_j \neg > \neg  t_0 -\frac16 \d  \neg > \neg  t_0 - \frac16 \d_0  \neg > \neg  \frac56 t_0 $, we obtain
 \beas
       \vf   (t_j   , x_j  )  +   \frac56    \wt{\wp}   t_0
  & < &  \vf   (t_j   , x_j  )   +   \wt{\wp}   t_j   =   \cY_{t_j}
     \le     Y^{t_j, x_j, \wh{\mu} ,    \beta (   \wh{\mu} ) }_{t_j} \Big( \t,
   \f  \Big( \t, X^{t_j, x_j, \wh{\mu} ,    \beta (   \wh{\mu} ) }_{\t} \Big)    \Big) \nonumber  \\
   & \le &   \underset{\mu \in \cU_{t_j}}{\esssup}
  Y^{t_j, x_j,  \mu, \beta (  \mu ) }_{t_j} \left( \t_{\beta,\mu},   \f  \big( \t_{\beta,\mu},
  X^{t_j, x_j, \mu, \beta (  \mu ) }_{\t_{\beta,\mu}} \big)    \right)     , \q \pas ,
      \eeas
 where $ \t_{\beta,\mu}   \dfnn  \inf \Big\{s  \neg \in \neg  ( t_j, T  ] \neg  :
  \big( s,X^{t_j, x_j, \mu, \beta (  \mu )}_s \big)
  \neg  \notin   \neg   O_{  \frac23 \d  } (t_j, x_j) \Big\} $, $\fa \mu \in \cU_{t_j}$.
 Taking essential  infimum over $\beta \in \fB_{t_j} $  and applying Theorem \ref{thm_DPP}
 with $(t,x,\d) = (t_j,x_j,\frac23 \d)$, we see from \eqref{eq:s437} that \pas
   \beas
   \vf  \big(t_j   , x_j     \big)  \neg + \neg  \frac56 \wt{\wp}   t_0
  \neg \le \neg   \underset{\beta \in \fB_{t_j}}{\essinf} \, \underset{\mu \in \cU_{t_j}}{\esssup}
  Y^{t_j, x_j,  \mu, \beta (  \mu ) }_{t_j} \Big( \t_{\beta,\mu},   \f  \Big( \t_{\beta,\mu},
  X^{t_j, x_j, \mu, \beta (  \mu ) }_{\t_{\beta,\mu}} \Big)    \Big)
   \neg   \le   \neg   w_1 (t_j,x_j)
  \neg < \neg  \vf (t_j,x_j)  \neg + \neg  \frac56 \wt{\wp}   t_0 \,.
 \eeas
   A contradiction appears.
 Therefore, $\ul{w}_1$  is a viscosity supersolution of   \eqref{eq:PDE} with Hamiltonian  $\ul{H}_1$.

      \ss \no {\bf b)}   Next,  we  show that
       $\ol{w}_1 $   is a viscosity subsolution    of   \eqref{eq:PDE} with   Hamiltonian $\ol{H}_1$.
 \if{0}
Clearly, $\cU_T = \hU$ and $\fB_T =  \fB_T  $  collects all functions $\beta: \hU \to \hV$ with  linear growth.
It follows that
  \beas
 \ol{w}_1(T,x) =    w_1(T,x)=  \underset{\beta \in \fB}{\inf} \, \underset{u \in \hU}{\sup}  \;
  Y^{T,x,u, \beta ( u )}_T \Big(T,g\big( X^{T,x,u, \beta ( u )}_T \big)  \Big)
  =  \underset{\beta \in \fB}{\inf} \, \underset{u \in \hU}{\sup}  \;  Y^{T,x,u, \beta ( u ) }_T \big(T,g(x)  \big)  = g(x), \q \fa x \in \hR^k
  \eeas
 \fi
  Let $(t_0,x_0, \vf) \in (0,T) \times \hR^k \times \hC^{1,2}\big([0,T] \times \hR^k\big)$
  be  such that $ \ol{w}_1  (t_0,x_0) = \vf (t_0,x_0)  $  and that  $\ol{w}_1-\vf$ attains a strict local maximum  at $(t_0,x_0)$, i.e., for some   $\d_0 \in  \big(0,   t_0 \land  ( T-t_0 ) \big) $
   \beas  
     (\ol{w}_1 - \vf) (t,x) <  (\ol{w}_1 - \vf) (t_0,x_0) = 0,
   \q \fa (t,x) \in O_{\d_0}  (t_0,x_0) \big\backslash \big\{ (t_0,x_0) \big\} .
   \eeas

    We still  denote $\big(\vf(t_0  ,x_0),  D_x \vf    (t_0  ,x_0),  D^2_x \vf    (t_0  ,x_0)\big) $
     by $(y_0,z_0,\G_0 )$.  If  $\ol{H}_1 \big(t_0,x_0, y_0,z_0,\G_0\big) =  \infty$,
    then
  \beas 
     &&       \neg  -   \frac{\pa }{\pa t} \vf (t_0,x_0)
      \neg - \neg  \ol{H}_1 \big(t_0,x_0, y_0,z_0,\G_0\big)
      \le  0    
     \eeas
     holds automatically.
   To make a contradiction,  we   assume    that when $\ol{H}_1 \big(t_0,x_0, y_0,z_0,\G_0\big) <     \infty $,
  \bea  \label{eq:s637}
 \varrho \dfnn   -   \frac{\pa }{\pa t} \vf (t_0,x_0)
      \neg - \neg  \ol{H}_1 \big(t_0,x_0, y_0,z_0,\G_0\big)  > 0 .
 \eea
  It is easy to see that
   \bea
   \ol{H}_1 \big(t_0,x_0, y_0,z_0,\G_0\big)  & \tneg \ge & \tneg  \lmtd{n \to \infty}  \;     \underset{u \in \hU  }{\sup} \;
 \underset{v \in \sO^n_u  }{\inf}   \;     H (t_0,x_0, y_0,z_0,\G_0,u,v)
    \ge \lmtd{n \to \infty}  \;
 \underset{v \in \sO^n_{u_0}  }{\inf}   \;  H (t_0,x_0, y_0,z_0,\G_0,u_0,v) \nonumber \\
 & \tneg  = & \tneg   \underset{v \in \ol{O}_\k(v_0)   }{\inf}       \;  H (t_0,x_0, y_0,z_0,\G_0,u_0,v) . \label{eq:xvx141}
   \eea
   For any $v \in \ol{O}_\k(v_0)$, one can deduce from  \eqref{eq:xvx311}  that
   $  | H (t_0,x_0, y_0,z_0,\G_0,u_0,v)| \le \frac14 (C^0_\vf)^2 + C^0_\vf C(\k,x_0) + C(\k,x_0)$, where
   $C^0_\vf 
   = \big|\vf (t_0,x_0)\big| + \big| D_x \vf (t_0,x_0) \big| + \big| D^2_x \vf (t_0,x_0) \big| $ as set in part (a).
      It then follows from \eqref{eq:xvx141} that
     \beas
   \ol{H}_1 \big(t_0,x_0, y_0,z_0,\G_0\big)
    \ge    \underset{v \in \ol{O}_\k(v_0)   }{\inf}       \;  H (t_0,x_0, y_0,z_0,\G_0,u_0,v)
    \ge - \frac14 (C^0_\vf)^2 - C^0_\vf C(\k,x_0) - C(\k,x_0) > - \infty .
   \eeas
    Thus $\varrho < \infty$.

 \ss  Then one can find an  $m   \in \hN$    such that
  \bea \label{eq:s629}
  \neg  -  \neg  \frac{\pa \vf }{\pa t}(t_0,x_0) - \frac78 \, \varrho     \ge
             \underset{u \in \hU  }{\sup} \;  \underset{v \in \sO^m_u  }{\inf}  ~\, \lsup{  u' \to u}\; \underset{(t,x,y,z,\G)
   \in O_{1/m} (t_0,x_0 ,y_0,z_0,\G_0)}{\sup} \;  H (t ,x  , y, z ,\G ,u',v) .
 \eea
As $\vf \in  \hC^{1,2}\big([0,T] \times \hR^k\big)$, there exists a $\d < \frac{1}{2m} \land  \d_0$ such that for any $ (t,x) \in \ol{O}_\d (t_0,x_0)$
   \bea
 &&   \Big| \frac{\pa \vf }{\pa t}(t,x)   -    \frac{\pa \vf }{\pa t}(t_0,x_0) \Big|    \le    \frac18 \,\varrho
             ~\;\;   \label{eq:s631} \\
  \hb{and} &&
 \big| \vf  (t,x)     -    \vf (t_0,x_0) \big|    \vee    \Big| D_x \vf    (t  ,x)
   -    D_x \vf    (t_0  ,x_0)\Big|
    \vee    \Big| D^2_x \vf    (t  ,x)    -    D^2_x \vf    (t_0  ,x_0)\Big|    \le    \frac{1}{2m} \,, \nonumber
      \eea
 the latter of which together with \eqref{eq:s629} implies that
  \beas
  \neg  -    \frac{\pa \vf }{\pa t}(t_0,x_0) - \frac78 \, \varrho     \ge
             \underset{u \in \hU  }{\sup} \;  \underset{v \in \sO^m_u  }{\inf}  ~\, \lsup{  u' \to u}\;
 \underset{(t,x )     \in \ol{O}_\d  (t_0,x_0)}{\sup} \;  H (t ,x , \vf(t   ,x),
     D_x \vf    (t  ,x),  D^2_x \vf    (t  ,x)  ,u',v)   .
 \eeas
  Then  for any $ u \in    \hU $,
 there exists a $ \fP_o (u) \in \sO^m_u  $ such that
 \beas
   \neg  -    \frac{\pa \vf }{\pa t}(t_0,x_0) - \frac34 \, \varrho  \ge
  \lsup{  u' \to u}\; \underset{(t,x )     \in \ol{O}_\d  (t_0,x_0)}{\sup} \;  H \big(t ,x , \vf(t   ,x),
     D_x \vf    (t  ,x),  D^2_x \vf    (t  ,x)  ,u', \fP_o(u) \big)   ,
 \eeas
 and we can find   a $ \l(u) \in  (0 , 1  )  $ such that
 for any $ u'   \in       O_{\l(u)}(u)  $
\bea \label{eq:s613}
   \neg  -  \neg  \frac{\pa \vf }{\pa t}(t_0,x_0) - \frac58 \, \varrho  \ge
    \underset{(t,x )     \in \ol{O}_\d  (t_0,x_0)}{\sup} \;  H \big(t ,x , \vf(t   ,x),
     D_x \vf    (t  ,x),  D^2_x \vf    (t  ,x)  ,u', \fP_o(u)\big)    .
 \eea

 Set $ \wt{\l}(u_0) = \l(u_0) $ and $ \wt{\l}(u) = \l(u) \land \big(\frac12 [u]_{\overset{}{\hU}} \big) $
 for any $u \in \hU \backslash \{u_0\}$.
 Since the separable metric space  $ \hU$ is   Lindel\"of,
  $\big\{\fO (u) \dfnn  O_{\wt{\l}(u)}(u) \big\}_{ u \in    \hU  } $ has a countable subcollection    $\{ \fO (u_i) \}_{i \in \hN
   }$ to cover   $\hU   $.       It is  clear that
 \beas
  \fP(u) \dfnn \sum_{ i \in \hN} \b1_{\big\{ u \in \fO (u_i) \backslash \,
  \underset{j<i}{\cup} \fO (u_j) \big\}} \fP_o (u_i) \in \hV  ,  \q  \fa     u  \in   \hU
 \eeas
   defines a $   \sB(\hU) /\sB(\hV)-$measurable function.

  \ss  Given $  u \neg \in   \neg   \hU  $,
    there exists an $i  \neg \in  \neg  \hN$ such that $u  \neg \in \neg  \fO (u_i) \backslash \,
  \underset{j<i}{\cup} \fO (u_j)$. If $u_i = u_0$,
    \bea    \label{eq:xux737a}
     \big[ \fP ( u) \big]_{\overset{}{\hV}} = \big[\fP_0( u_i)\big]_{\overset{}{\hV}}
       \le  \k    + m      [u_i]_{\overset{}{\hU}} =  \k    \le  \k    + m      [u]_{\overset{}{\hU}}  .
   \eea
   On the other hand, if $u_i \ne u_0$,  then
   $ [u_i]_{\overset{}{\hU}} \le [u]_{\overset{}{\hU}}
  \neg  + \neg     \rho_{\overset{}{\hU}} (u, u_i) \le [u]_{\overset{}{\hU}}
  \neg  + \neg    \wt{\l} ( u_i) \le  [u]_{\overset{}{\hU}}
  \neg  + \neg   \frac12 [u_i]_{\overset{}{\hU}} $,  and it follows that
  \bea   \label{eq:xux737}
     \big[ \fP ( u) \big]_{\overset{}{\hV}} = \big[\fP_0( u_i)\big]_{\overset{}{\hV}} \le  \k
     + m      [u_i]_{\overset{}{\hU}}
  \neg  \le  \neg \k+ 2 m       [u]_{\overset{}{\hU}}.
  \eea
   Also,  we see from \eqref{eq:s613} that
 \beas
  \neg  -    \frac{\pa \vf }{\pa t}(t_0,x_0) - \frac58 \, \varrho & \ge &
    \underset{(t,x )     \in \ol{O}_\d  (t_0,x_0)}{\sup} \;  H \big(t ,x , \vf(t   ,x),
     D_x \vf    (t  ,x),  D^2_x \vf    (t  ,x)  ,u, \fP_0(u_i) \big)  \\
  &=&  \underset{(t,x )     \in \ol{O}_\d  (t_0,x_0)}{\sup} \;  H (t ,x , \vf(t   ,x),
     D_x \vf    (t  ,x),  D^2_x \vf    (t  ,x)  ,u, \fP(u) )  ,
 \eeas
which together with \eqref{eq:s631} implies that
 \bea  \label{eq:s617}
  \neg  -  \neg  \frac{\pa \vf }{\pa t}(t ,x ) - \frac12 \, \varrho  \ge
          H (t ,x , \vf(t   ,x),
     D_x \vf    (t  ,x),  D^2_x \vf    (t  ,x)  ,u, \fP(u) )  , \q \fa (t,x) \in \ol{O}_\d(t_0,x_0) , ~ \fa u \in   \hU.
 \eea

  Similar to \eqref{eq:s427}, we set
 $        \wp    \dfnn   \min   \big\{ (\vf \neg - \neg \ol{w}_1      ) (t,x)  \neg : (t,x )  \neg \in \neg
     \ol{O}_\d  (t_0,x_0) \big\backslash O_{\frac{\d}{3}} (t_0,x_0) \big\}  > 0     $
   \if{0}
 Since the set $ \ol{O}_\d  (t_0,x_0) \big\backslash O_{\frac{\d}{3}} (t_0,x_0)$ is  compact, there exists
 a sequence $\{(t_n,x_n)\}_{n \in \hN} $ of $  \ol{O}_\d  (t_0,x_0) \big\backslash O_{\frac{\d}{3}} (t_0,x_0)$
 that  converges to some $(t_*,x_*) \in \ol{O}_\d  (t_0,x_0) \big\backslash O_{\frac{\d}{3}} (t_0,x_0) $
 and satisfies   $ \wp \neg = \neg  \lmtd{n \to \infty} (\ul{w}_1  \neg - \neg  \vf) (t_n,x_n) $.
 The lower semicontinuity of $\ul{w}_1$ and the continuity of $\vf$ imply that  $\ul{w}_1  \neg - \neg  \vf$ is also lower semicontinuous.
 Thus, it follows that $ \wp \le (\ul{w}_1-\vf)(t_*,x_*) \le \lmtd { n \to \infty } (\ul{w}_1-\vf)(t_n,x_n ) = \wp  $,
 which together with \eqref{eq:a071b} shows that
  \bea   \label{eq:s427b}
  \wp  = \min \big\{ (\ul{w}_1 - \vf) (t,x): (t,x ) \in
     \ol{O}_\d  (t_0,x_0) \big\backslash O_{\frac{\d}{3}} (t_0,x_0) \big\} = (\ul{w}_1-\vf)(t_*,x_*)  >  0 .
  \eea
  \fi
  and     $\dis   \wt{\wp} \neg \dfnn  \neg   \frac{    \wp \land  \varrho}{ 2 (1 \vee \g) T }  \neg   > \neg  0$.
    Let $ \big\{(t_j, x_j)\big\}_{j \in \hN}$ be a sequence of $O_{\frac{\d}{6}} (t_0,x_0) $ such   that
     \beas   
 \lmt{j \to \infty} (t_j,x_j) = (t_0,x_0)   \q \hb{and} \q  \lmt{j \to \infty} w_1 (t_j,x_j) = \ol{w}_1 (t_0,x_0) = \vf (t_0,x_0) = \lmt{j \to \infty} \vf (t_j,x_j) .
     \eeas
   So one can find a   $j \in \hN$ that
    \bea  \label{eq:s437b}
  \big| w_1 (t_j,x_j)- \vf (t_j,x_j) \big| < \frac56 \wt{\wp}   t_0 .
   \eea

  For any $\mu \in \cU_{t_j}$, the measurability of function $\fP$  shows that $ \big( \wh{\beta}(\mu) \big)_s \dfnn \fP(\mu_s)$,  $s \in [t_j,T]$ is a $\hV-$valued, $\bF-$progressively measurable process. By \eqref{eq:xux737a} and \eqref{eq:xux737},
  \beas   
  \big[ \big( \wh{\beta}(\mu) \big)_s \big]_{\overset{}{\hV}}
  = \big[   \fP (\mu_s )   \big]_{\overset{}{\hV}} \le \k + 2m [\mu_s]_{\overset{}{\hU}} , \q \fa s \in [t_j,T] .
  \eeas
  Let $E \int_{t_j}^T  [  \mu_s ]^q_{\overset{}{\hU}} \,   ds   <   \infty$ for some $q >2$.
   It then follows   that
        \beas
    E \int_{t_j}^T  \big[ \big( \wh{\beta}(\mu) \big)_s \big]^q_{\overset{}{\hV}} \,   ds
    \le 2^{q-1} \k^q T + 2^{2q-1} m^q  E \int_{t_j}^T  [  \mu_s ]^q_{\overset{}{\hU}} \,   ds   <   \infty .
    \eeas
       So  $ \wh{\beta}(\mu) \neg \in \neg  \cV_{t_j} $.
       Let    $\mu^1, \mu^2  \neg \in \neg  \cU_{t_j}$ such that
 $\mu^1  \neg = \neg   \mu^2 $,  $ds  \neg \times \neg  d P -$a.s.  on $ \[t_j,\t\[ \, \cup \, \[\t, T\]_A   $
 for some $ \t  \neg \in  \neg   \cS_{t_j,T}$  and $ A   \neg \in  \neg \cF_\t  $.
 Then it directly follows that $ \big( \wh{\beta}(\mu^1) \big)_s  \neg = \neg  \fP(\mu^1_s)
  \neg = \neg  \fP(\mu^2_s)  \neg = \neg  \big( \wh{\beta}(\mu^2) \big)_s    $, $ds  \neg \times \neg  d P -$a.s.
  on $ \[t_j,\t\[ \, \cup \, \[\t, T\]_A   $. Hence, $\wh{\beta}  \neg \in \neg  \fB_{t_j}$.

 \ss  Let   $\mu \in   \cU_{t_j} $.
     We set $\Th_\mu \dfnn \big( t_j, x_j, \mu ,    \wh{\beta} (   \mu )  \big)$
  and  define
 \beas
  \q  \t_\mu = \t_{\wh{\beta}, \mu}  \dfnn  \inf \Big\{s  \neg \in \neg  ( t_j, T  ] \neg  : \big( s,X^{\Th_\mu}_s \big)
  \neg  \notin   \neg   O_{ \frac23 \d} (t_j, x_j) \Big\} \in \cS_{t_j,T} .
 \eeas
    As 
 $ \big| \big(T, X^{\Th_\mu}_T\big) -(t_j,x_j) \big| \neg  \ge  \neg  T \neg - \neg t_j
 \neg  \ge  \neg   T  \neg - \neg t_0  \neg - \neg  |t_j \neg - \neg t_0|  \neg > \neg  \d_0  \neg - \neg \frac{\d}{6}
    \neg >  \neg  \frac23 \d $,
  the continuity of $ X^{\Th_\mu} $ implies that \pas
  \bea
  && \qq \t_\mu  < T  \q \hb{and}  \q
  \big(\t_\mu \land  s , X^{\Th_\mu}_{\t_\mu \land  s}   \big) \in    \ol{O}_{  \frac23 \d }(t_j, x_j)
   \subset \ol{O}_{\frac56 \d}  (t_0, x_0) ,  \q \fa    s \in [t_j, T]  .  \qq   \qq     \label{eq:a239} \\
  && \hb{In particular}, \q  \big( \t_\mu, X^{\Th_\mu}_{\t_\mu} \big)  \in   \pa O_{  \frac23 \d   }(t_j, x_j)
   \subset \ol{O}_{\frac56 \d}  (t_0,x_0) \big\backslash O_{\frac{\d}{2}} (t_0,x_0) .
     \label{eq:a041b}
  \eea

      The continuity of $\vf$, $X^{\Th_\mu} $ and \eqref{eq:a239} show that
  $\cY^\mu_s \dfnn \vf \big( \t_\mu \land  s, X^{\Th_\mu}_{\t_\mu \land  s} \big)
    - \wt{\wp} ( \t_\mu \land  s) $, $ s \in [t_j, T] $ defines a bounded $\bF-$adapted continuous process.
 Applying It\^o's formula yields that
    \bea
         \cY^\mu_s   & \tneg = &  \tneg   \cY^\mu_T     +   \int_s^T  \ff^\mu_r   dr
    -  \int_s^T  \cZ^\mu_r   d B_r ,   \q    s \in [t_j,T]  ,   \label{eq:s403}
    \eea
where $\cZ^\mu_r  \dfnn        \b1_{\{r < \t_\mu\}}    D_x \vf    \big( r, X^{\Th_\mu}_r \big)
  \cd   \si \big( r, X^{\Th_\mu}_r, \mu_r,   \big(\wh{\beta} (  \mu )  \big)_r \big)$ and
 $$
     \ff^\mu_r   \neg   \dfnn   \neg       \b1_{\{ r \, < \t_\mu  \}}
  \bigg\{ \wt{\wp}      -       \frac{\pa \vf }{\pa t} \big( r, X^{\Th_\mu}_r \big)
      -       D_x \vf    \big( r, X^{\Th_\mu}_r\big)   \cd    b \big( r, X^{\Th_\mu}_r \neg , \mu_r,
  (\wh{\beta} (  \mu )  )_r \big)
         -         \frac12  trace\Big( \si \si^T \big( r, X^{\Th_\mu}_r \neg , \mu_r,
  (\wh{\beta} (  \mu ) )_r \big) \cd    D^2_x \vf \big( r, X^{\Th_\mu}_r\big) \neg \Big) \neg \bigg\} .
     $$
 As $\vf  \neg \in  \neg   \hC^{1,2}\big([t,T] \neg  \times \neg  \hR^k\big)$,
  the measurability of $b$, $\si$, $X^{\Th_\mu}$, $\mu$ and $ \wh{\beta} (  \mu ) $ implies   that
   both $\cZ^\mu$ and $\ff^\mu$ are     $\bF-$progressively    measurable.
  Let $ \wt{C}_\vf \dfnn \underset{(t,x) \in \ol{O}_{\frac56 \d} (t_0,x_0) }{\sup} \big|D_x \vf (t,x)\big| < \infty $.
  Similar to \eqref{eq:s641}, we see from  \eqref{b_si_linear_growth}, \eqref{b_si_Lip} and \eqref{eq:a239}  that
  \bea
    E \left[  \bigg( \int_{t_j}^T \neg |\cZ^\mu_s|^2   \,  ds  \bigg)^{ p /2 } \right]
 \neg    \le \neg   c_0  \wt{C}^p_\vf  \bigg(   \big(1  \neg + \neg  |x_0 |  \neg + \neg  \d    \big)^p
  \neg + \neg   \bigg\{ E   \neg   \int_{t_j}^T \dneg [\mu_s]^2_{\overset{}{\hU}} \,     ds    \bigg\}^{ p /2 }
  \neg + \neg   \bigg\{ E    \neg  \int_{t_j}^T \dneg \big[ (\wh{\beta} (  \mu ) )_s \big]_{\overset{}{\hV}}^2 \,   ds
      \bigg\}^{ p /2 } \bigg)  \neg < \neg  \infty ,
  \eea
  i.e.  $\cZ^\mu \in \hH^{2,p} _{\bF  } \big( [t_j,T], \hR^d \big)$.
     Hence,  $\big\{\big( \cY^\mu_s, \cZ^\mu_s \big) \big\}_{s \, \in  [t_j,T]}$
  solves the   BSDE$\big( t_j , \cY^\mu_T , \ff^\mu  \big)$.

 Let $\ell$, $\th_1$ and $\th_2$ still  be the  continuous functions considered in part (a).
 Like $\f$ in \eqref{eq:xvx151},
    \beas
    \wt{\f} (t,x) \dfnn \th_1 (t,x) \ell (x) + \big(1-\th_1 (t,x)\big)
    (\vf(t,x) - \wp \th_2 (t,x) )   ,   \q  \fa (t,x) \in [t_j,T] \times \hR^k
    \eeas
     define a continuous function with $ \wt{\f}  \ge  w_1$.
  Similar to \eqref{eq:xux811} and \eqref{eq:xux813}, we can deduce from \eqref{eq:a041b},
  \eqref{eq:a239}, \eqref{eq:s617} and \eqref{f_Lip} that $ \cY^\mu_T \neg  \ge  \neg  \wt{\f} \big( \t_\mu, X^{\Th_\mu}_{\t_\mu} \big) $,
  \pas ~ and that  $\ff^\mu_s (\o)  \neg  \ge  \neg  f^{\Th_\mu }_{\t_\mu} \big(s ,\o, \cY^\mu_s  (\o),  \cZ^\mu_s (\o) \big) $
  for $ds  \neg \times \neg  d P-$a.s.  $ (s,\o)   \neg \in \neg  [t_j,T]  \neg \times \neg  \O$.
  \if{0}
 Since   $ \ol{w}_1(t,x)  \ge  w_1 (t,x)$
 for any $(t,x) \in [0,T] \times \hR^k$, we can deduce from  \eqref{eq:a041b} that
      \beas
   \cY^\mu_T   \ge     \vf \big( \t_\mu, X^{\Th_\mu}_{\t_\mu} \big)    -   \wt{\wp} T
 >  \vf \big( \t_\mu, X^{\Th_\mu}_{\t_\mu} \big) - \wp
  \ge   \wt{\f} \big( \t_\mu, X^{\Th_\mu}_{\t_\mu} \big)   .
  \eeas
  Also,      \eqref{eq:a239}, \eqref{eq:s617}  and \eqref{f_Lip}  show   that for $ds \times d P-$a.s.
  $ (s,\o)  \in [t_j,T] \times \O$
     \beas
     \ff^\mu_s (\o)  & \tneg  \ge  &  \tneg  
 \b1_{\{ s \, < \t_\mu (\o) \}} \left\{  \wt{\wp} + \frac12 \varrho   +    f \neg \left( s,\o, X^{\Th_\mu}_s (\o),
       \cY^\mu_s (\o) -\wt{\wp}s,    \cZ^\mu_s (\o), \mu_s (\o),  ( \wh{\beta} (  \mu )  )_s (\o) \right)  \right\}  \\
         &  \tneg  \ge  &  \tneg  \b1_{\{ s \, < \t_\mu (\o) \}} \left\{  \wt{\wp} + \frac12 \varrho  - \g \wt{\wp} T
          +    f \neg \left(s, \o, X^{\Th_\mu}_s (\o),
      \cY^\mu_s (\o),       \cZ^\mu_s  (\o), \mu_s (\o),  ( \wh{\beta} (  \mu )  )_s  (\o) \right)  \right\}   \\
   &  \tneg   \ge   &  \tneg   f^{\Th_\mu }_{\t_\mu} \big(s ,\o, \cY^\mu_s  (\o),  \cZ^\mu_s (\o) \big) .
       \eeas
       \fi
  As   $f^{\Th_\mu}_{\t_\mu}$ is  Lipschitz continuous in $(y,z)$, we know from Proposition
  \ref{prop_BSDE_estimate_comparison} (2) that $P-$a.s.
   \beas
   \cY^\mu_s    \ge  Y^{\Th_\mu}_s \left( \t_\mu,
   \wt{\f}  \big( \t_\mu, X^{\Th_\mu}_{\t_\mu} \big)    \right)  , \q  \fa  s \in [t_j,T] .
   \eeas
  Letting $s=t_j$ and using the fact that $t_j \neg > \neg  t_0 -\frac16 \d  \neg > \neg  t_0 - \frac16 \d_0  \neg > \neg  \frac56 t_0 $, we obtain
   \beas
       \vf   (t_j   , x_j  )  -   \frac56    \wt{\wp}   t_0
   >   \vf   (t_j   , x_j  )   -   \wt{\wp}   t_j   =   \cY^\mu_{t_j}
      \ge      Y^{t_j, x_j, \mu ,    \wh{\beta} (   \mu ) }_{t_j} \Big( \t_\mu,
   \wt{\f}  \Big( \t_\mu, X^{t_j, x_j, \mu ,    \wh{\beta} (   \mu ) }_{\t_\mu} \Big)    \Big)    , \q  \pas
      \eeas
   Taking essential supremum over $\mu \in \cU_{t_j} $ and applying Theorem \ref{thm_DPP}
 with $(t,x,\d) = (t_j,x_j,\frac23 \d)$, we see from \eqref{eq:s437b} that   \pas
 \beas
   \vf   (t_j   , x_j  ) \neg - \neg    \frac56    \wt{\wp}   t_0 & \tneg  \dneg  \ge  &  \tneg  \dneg  \underset{\mu \in \cU_{t_j}}{\esssup}
  Y^{t_j, x_j, \mu ,    \wh{\beta} (   \mu ) }_{t_j}   \Big(  \t_\mu ,
  \wt{\f} \Big(  \t_\mu , X^{t_j, x_j, \mu ,    \wh{\beta} (   \mu ) }_{ \t_\mu } \Big) \Big)
   \neg  \ge  \neg   \underset{\beta \in \fB_{t_j} }{\essinf} \, \underset{\mu \in \cU_{t_j}}{\esssup}
  Y^{t_j, x_j, \mu ,    \beta (   \mu ) }_{t_j}   \Big( \t_{\beta,\mu},
  \wt{\f} \Big( \t_{\beta,\mu}, X^{t_j, x_j, \mu ,    \beta (   \mu ) }_{\t_{\beta,\mu}} \Big) \Big)  \\
 &  \tneg  \dneg   \ge  &  \tneg  \dneg   w_1 (t_j,x_j)
  \neg > \neg  \vf (t_j,x_j)  \neg - \neg  \frac56 \wt{\wp}   t_0 \, ,
 \eeas
 where $ \t_{\beta,\mu}   \dfnn  \inf \Big\{s  \neg \in \neg  ( t_j, T  ] \neg  :
  \big( s,X^{t_j, x_j, \mu, \beta (  \mu )}_s \big)
  \neg  \notin   \neg   O_{  \frac23 \d  } (t_j, x_j) \Big\} $.  A contradiction appears.
 Therefore,     $\ol{w}_1$  is a viscosity supersolution of   \eqref{eq:PDE} with Hamiltonian $\ol{H}_1$.     \qed

       \medskip

\bibliographystyle{siam}
\bibliography{SDGVU_m}

 \end{document}